\documentclass{article}
\usepackage[T2A]{fontenc}
\usepackage[cp866]{inputenc}
\usepackage[english,russian]{babel}
\usepackage[tbtags]{amsmath}
\usepackage{amsfonts,amssymb,mathrsfs,amscd}

\let\No=\textnumero

\usepackage[hyper]{styl2e}

\numberwithin{equation}{section}
\theoremstyle{plain}
\newtheorem{theorem}{Теорема}[section]

\newtheorem{propos}{Предложение}[section]
\newtheorem{cor}{Следствие}[section]
\theoremstyle{definition}

\newtheorem{proof}{Доказательство}
\newtheorem{remark}{Замечание}[section]

\begin{document}
\begin{center}
\large{\textbf{Об операторах смежности локально конечных графов}}

\end{center}
\author[V.\,I.~Trofimov]{В.\,И.~Трофимов}
\address{Институт математики и механики им.~Н.\,Н.~Красовского УрО РАН, г.~Екатеринбург;\\
Уральский федеральный университет им. первого Президента России~Б.\,Н.~Ельцина, г.~Екатеринбург;\\
Уральский математический центр, г.~Екатеринбург}
\email{trofimov@imm.uran.ru} 

\date{}

\maketitle

\begin{fulltext}

\begin{abstract}

Граф $\Gamma$ называется локально конечным, если у графа $\Gamma$
для каждой вершины $v$
множество $\Gamma(v)$ смежных с ней вершин конечно.
Для произвольного локально конечного графа $\Gamma$ с множеством вершин $V(\Gamma)$ и произвольного поля $F$
на $F^{V(\Gamma)}$ (векторном пространстве над $F$ всех функций $V(\Gamma) \to F$ с естественными покомпонентными операциями) определен линейный оператор
$A^{({\rm alg})}_{\Gamma,F}: F^{V(\Gamma)} \to F^{V(\Gamma)}$, посредством
$(A^{({\rm alg})}_{\Gamma,F}(f))(v) = \sum_{u \in \Gamma(v)}f(u)$
для всех $f \in F^{V(\Gamma)}$, $v \in V(\Gamma)$.
В случае конечного \mbox{графа} $\Gamma$ отображение $A^{({\rm alg})}_{\Gamma,F}$ есть хорошо известный оператор, определяемый матрицей
смежности графа $\Gamma$ (над $F$), и теория собственных значений и собственных функций таких операторов
составляет (по крайней мере, \mbox{в случае} $F = \mathbb{C}$) хорошо разработанный раздел
теории конечных графов.
В настоящей работе разрабатывается теория собственных значений и собственных функций операторов
$A^{({\rm alg})}_{\Gamma,F}$
для бесконечных локально конечных графов $\Gamma$ (впрочем, отдельные ее результаты могут представлять интерес для конечных графов) и произвольных полей $F$,
хотя особый акцент делается на случай, когда $\Gamma$ --- связный граф с ограниченными в совокупности степенями вершин и $F = \mathbb{C}$.
Предпринимавшиеся ранее попытки в этом направлении не были, по мнению автора, вполне
удовлетворительными в том смысле, что ограничивались рассмотрением лишь собственных функций весьма специального вида
(и соответствующих им собственных значений).

Библиография: 18 наименований.
\end{abstract}


\markright{Об операторах смежности локально конечных графов}

\footnotetext[0]{Работа выполнена при поддержке Министерства науки и высшего образования
Российской Федерации в рамках проекта ``Уральский математический центр''
(соглашение \No~075-02-2022-877).}

\section{Введение}
\label{s1}

Под графом всюду далее понимается неориентированный граф без петель и без
кратных ребер. Если $\Gamma$ --- граф, то $V(\Gamma)$ и $E(\Gamma)$ --- соответственно его множество вершин
и множество ребер. Для $x \in V(\Gamma)$ через $\Gamma(x)$ обозначается множество всех
смежных с $x$ вершин графа $\Gamma$ (таким образом, $|\Gamma(x)|$ --- степень вершины $x$).
Граф $\Gamma$ называется локально конечным, если степени всех его вершин конечны.

Пусть $\Gamma$ --- локально конечный граф, $F$ --- поле и $F^{V(\Gamma)}$ --- векторное пространство
над $F$, образованное всеми функциями $V(\Gamma) \to F$ (с естественными покомпонентными операциями
сложения и умножения на скаляр).
Определим линейное отображение $A^{({\rm alg})}_{\Gamma,F}: F^{V(\Gamma)} \to F^{V(\Gamma)}$,
полагая для произвольных $f \in F^{V(\Gamma)}$
и $v \in V(\Gamma)$
$$(A^{({\rm alg})}_{\Gamma,F}(f))(v) := \sum_{u \in \Gamma(v)}f(u).$$
 Настоящая работа ориентирована, прежде всего, на исследование собственных значений и собственных векторов (собственных функций из $F^{V(\Gamma)}$)
линейного отображения  $A^{({\rm alg})}_{\Gamma,F}$.

Хотя собственные значения и собственные векторы так определенного линейного отображения
$A^{({\rm alg})}_{\Gamma,F}$ не зависят от наличия топологий на $F$ и $F^{V(\Gamma)}$,
нам будет удобно без потери общности предполагать в работе, что поле $F$ снабжено некоторым (возможно, тривиальным)
абсолютным значением $|.|_{\rm v}$, естественным образом превращающим $F$ в метрическое пространство (см.,
например,~\cite[с. 322]{Lang}),
а $F^{V(\Gamma)}$ снабжено соответствующей топологией произведения, превращающей $F^{V(\Gamma)}$
в топологическое векторное пространство. При этом линейное отображение $A^{({\rm alg})}_{\Gamma,F}$
превращается в
непрерывный линейный оператор в топологическом векторном пространстве $F^{V(\Gamma)}$, который
мы будем обозначать через $A_{\Gamma,F}$ и называть {\it оператором смежности} графа $\Gamma$ над
полем $F$ с абсолютным значением $|.|_{\rm v}$.

В случае конечного графа $\Gamma$ отображение $A^{({\rm alg})}_{\Gamma,F}$ есть хорошо известный оператор, определяемый матрицей
смежности графа $\Gamma$ (над $F$), и теория спектров таких операторов хорошо разработана (по крайней мере, в случае  $F = \mathbb{C}$), см.,
например,~\cite{C},~\cite{Bro}.
Предпринимались представляющие интерес попытки построения спектральной теории произвольных локально конечных графов
на основе их матриц смежности (при этом условие не более чем счетности и даже условие связности не являются сильными ограничениями
на графы, поскольку, по существу, достаточно рассматривать связные компоненты графов). Пожалуй, наиболее известная из них -
теория из~\cite{BM},~\cite{MW}. В связи с этим отметим следующее. Наше определение оператора смежности $A_{\Gamma,F}$, даже
в случае не более чем счетного графа $\Gamma$ и $F = \mathbb{C}$ с естественным абсолютным значением, кардинально отличается от определения
оператора смежности не более чем счетного локально конечного графа из~\cite{BM},~\cite{MW}, где оператором смежности такого графа называется ограничение введенного выше оператора
$A^{({\rm alg})}_{\Gamma,\mathbb{C}}$ на гильбертово пространство $l_2(V(\Gamma))$
квадратично суммируемых функций в случае, когда степени вершин  $\Gamma$ ограничены в совокупности, а в общем случае
называется оператор в $l_2(V(\Gamma))$, являющийся замыканием
ограничения $A^{({\rm alg})}_{\Gamma,\mathbb{C}}$
на пространство функций с конечными носителями. Не удивительно поэтому, что у нас получается кардинально другая теория
для собственных значений и собственных функций. (Так, для произвольного бесконечного
локально конечного связного графа
$\Gamma$ у $A^{({\rm alg})}_{\Gamma,\mathbb{C}}$ или, что то же, у
$A_{\Gamma,\mathbb{C}}$ с произвольным
абсолютным значением на $\mathbb{C}$
имеются собственные функции, соответствующие всем, за не более чем счетным исключением, комплексным
числам, см. теорему~\ref{t6.1}, но в случае, когда степени вершин $\Gamma$ не превосходят некоторого $d \in \mathbb{Z}_{\geq 1}$,
в $l_2(V(\Gamma))$ могут попасть лишь те из этих функций, которые соответствуют числам из $\mathbb{R}_{\geq - d} \cap \mathbb{R}_{\leq d}$, см. замечание~\ref{r6.8}.
Еще нагляднее: для ``бесконечной цепочки'' $\Gamma$ с
$V(\Gamma) = \{v_i : i \in \mathbb{Z}\}$, $E(\Gamma) = \{ \{v_i,v_{i+1}\} : i \in \mathbb{Z}\}$
для каждого $\lambda \in \mathbb{C}$
у оператора $A_{\Gamma,\mathbb{C}}$ имеются соответствующие собственному значению $\lambda$ собственные функции
$f_{c_1,c_2} \in \mathbb{C}^{V(\Gamma)}$, определяемые при $\lambda \not \in \{-2,+2\}$ посредством
$f_{c_1,c_2}(v_i) := c_1 a^i + c_2 b^i$ для всех $i \in \mathbb{Z}$,
где $a, b$ --- комплексные корни уравнения $x^2 - \lambda x + 1 = 0$, $c_1, c_2$ --- произвольные не равные  одновременно нулю комплексные числа, а при $\lambda = \varepsilon 2$, $\varepsilon \in \{-1,+1\}$, посредством
$f_{c_1,c_2}(v_i) := \varepsilon^i(c_1 + c_2 i)$ для всех $i \in \mathbb{Z}$, $c_1, c_2$ --- произвольные не равные  одновременно нулю комплексные числа;
но ни для какого $\lambda \in \mathbb{C}$
у $A_{\Gamma,\mathbb{C}}$ нет собственных функций в пространстве квадратично суммируемых
функций $l^2(V(\Gamma)$.)
По мнению автора, даже в случае, когда степени вершин связного графа $\Gamma$ ограничены в совокупности и $F = \mathbb{C}$,
пространство квадратично суммируемых функций $l_2(V(\Gamma))$ зачастую излишне мало, чтобы вместить представляющие интерес собственные
функции оператора $A^{({\rm alg})}_{\Gamma,\mathbb{C}}$ (см., однако, сказанное далее во введении
в связи с \S~\ref{s6}).
Качественно установку, реализуемую в настоящей работе, можно охарактеризовать как нахождение собственных значений и собственных функций
операторов смежности $A^{({\rm alg})}_{\Gamma,F}$
в максимально больших пространствах $F^{V(\Gamma)}$
(с возможным дальнейшим их анализом на предмет попадания в то или иное представляющее интерес подмножество
пространства $F^{V(\Gamma)}$). В последующей работе автор планирует предложить, по крайней
мере в случае, когда $\Gamma$ --- связный граф с ограниченными в совокупности степенями вершин и
$F = \mathbb{C}$, $A_{\Gamma,F}$-инвариантное подпространство пространства $F^{V(\Gamma)}$, которое значительно меньше $F^{V(\Gamma)}$,
но лишено указанного недостатка подпространства $l_2(V(\Gamma))$.

Рассматриваемые в работе вопросы соприкасаются с рядом разработанных направлений математики (помимо теории
спектров графов, см. выше, это анализ на графах и, в частности, исследования операторов Лапласа локально конечных графов,
теория топологических векторных пространств и функциональный анализ, теория бесконечных систем линейных уравнений
и бесконечных матриц, теория линейных клеточных автоматов,
но также теория разностных уравнений и теория скалярных полей на графах).
Однако приходится констатировать: ввиду специфики рассматриваемых в работе $F^{V(\Gamma)}$ и $A_{\Gamma,F}$ (для бесконечных $\Gamma$)
мало что из указанных областей удается
непосредственно и в достаточной общности использовать при рассмотрении поднятых в работе вопросов. Важное исключение
--- теорема Теплица (см.~\cite{T}, а также~\cite{A},~\cite[с. 107]{L}), утверждающая (в одной из эквивалентных формулировок), что
система линейных уравнений над полем разрешима, если разрешима каждая ее конечная подсистема. Кроме того, рассуждения
из доказательства теоремы~\ref{t5.2} близки к рассуждениям
из доказательства теоремы 3 работы~\cite{Wang} и отдельные хорошо известные результаты используются в \S~\ref{s6}.
В целом работа относится к алгебраической теории графов. Отчасти она стимулировалась~\cite{Tr} и~\cite[вопрос 15.89]{Ko}.

Приведем краткое описание содержания статьи. Конечно, следует
иметь в виду, что получаемые результаты о связных
локально конечных графах применимы к связным компонентам произвольных локально конечных графов.

Итак, пусть $\Gamma$ --- локально конечный граф и $F$ --- поле с некоторым (возможно, тривиальным)
абсолютным значением $|.|_{\rm v}$.

В \S~\ref{s2} содержатся важные для дальнейшего вспомогательные результаты. В нем приводится удобная для использования
переформулировка теоремы Теплица, выводятся следствия из нее и доказывается,
что для любого $\lambda \in F$ оператор $A_{\Gamma,F} - \lambda E$ (где $E$ --- единичный оператор
на $F^{V(\Gamma)}$) отображает замкнутые подпространства топологического
векторного пространства $F^{V(\Gamma)}$ в замкнутые же подпространства (см. следствие~\ref{c2.1}).

В  \S~\ref{s3}  устанавливается ряд общих свойств собственных функций оператора $A_{\Gamma,F}$.
В частности, показывается, что (см. предложения~\ref{p3.3},~\ref{p3.4} и замечание~\ref{r3.3}) для заданных $\lambda \in F$ и $v \in V(\Gamma)$
каждый содержащий $v$ носитель  собственной функции оператора $A_{\Gamma,F}$,
соответствующей собственному значению $\lambda$, содержит минимальный (по включению)
содержащий $v$ носитель собственной функции оператора $A_{\Gamma,F}$,
соответствующей собственному значению $\lambda$, причем
эту последнюю функцию можно выбрать принимающей значения только из $F_0(\lambda)$, где $F_0$ --- простое подполе поля $F$,
а этот минимальный носитель связен в графе $\Gamma^2$ (где $\Gamma^2$ --- граф с $V(\Gamma^2) = V(\Gamma)$ и с ребрами,
соединяющими всевозможные различные вершины, удаленные одна от другой на расстояние $\leq 2$ в $\Gamma$).

В  \S~\ref{s4} для  $\lambda \in F$ и  $v \in V(\Gamma)$ определяется $(F,\lambda)$-пропагатор графа $\Gamma$ относительно $v$
как такая функция $f \in F^{V(\Gamma)}$, что $(A_{\Gamma,F} - \lambda E)(f) = \delta_v$, где $\delta_v(u) = 0$ для всех $u \in V(\Gamma) \setminus \{v\}$
и $\delta_v(v) = 1$. ($(F,\lambda)$-пропагатор $\Gamma$ относительно $v$ есть, по существу, функция
Грина или фундаментальное решение относительно $v$ для оператора $A_{\Gamma,F} - \lambda E$.) $(F,\lambda)$-пропагаторы играют ключевую роль далее в этой работе.
В \S~\ref{s4} устанавливаются общие свойства
$(F,\lambda)$-пропагаторов. Доказывается, в частности, что (см. теоремы~\ref{t4.1} и~\ref{t4.2}) отсутствие
$(F,\lambda)$-пропагатора графа $\Gamma$ относительно вершины $v$ равносильно принадлежности $v$ конечному носителю собственной
функции оператора $A_{\Gamma,F}$,
соответствующей собственному значению $\lambda$, а наличие
$(F,\lambda)$-пропагатора графа $\Gamma$ относительно вершины $v$ с конечным носителем равносильно
отсутствию собственной
функции оператора $A_{\Gamma,F}$, соответствующей собственному значению $\lambda$,
носитель которой содержал бы $v$. Следствием
первого из этих результатов является то, что (см. следствие~\ref{c4.2}) отсутствие $(F,\lambda)$-пропагатора графа $\Gamma$ относительно вершины $v$ возможно
лишь при весьма специальных $\lambda$ (которые, в частности, должны быть собственными значениями матрицы смежности над $F$
индуцированного конечного подграфа графа $\Gamma$, что влечет их алгебраичность над простым подполем поля $F$ и вещественность в случае $F = \mathbb{C}$).

В  \S~\ref{s5} для ряда представляющих интерес свойств, включая ряд спектральных свойств,
оператора $A_{\Gamma,F}$ устанавливается (для дальнейшего использования) их эквивалентность определенным свойствам пропагаторов
графа $\Gamma$.
Доказывается, в частности, что (см. теорему~\ref{t5.3}) элемент $\lambda \in F$ тогда и только тогда не является
собственным значением оператора $A_{\Gamma,F}$ (другими словами, $A_{\Gamma,F} - \lambda E$ инъективен),
когда относительно каждой вершины $v$ у графа $\Gamma$
имеется $(F,\lambda)$-пропагатор с конечным носителем. Заметим, что сюръективность $A_{\Gamma,F} - \lambda E$
эквивалентна наличию $(F,\lambda)$-пропагатора графа $\Gamma$ относительно каждой вершины $v$ (см. теорему~\ref{t5.1}).
Доказывается также, что (см. теорему~\ref{t5.2}) наличие $(F,\lambda)$-пропагатора с бесконечным носителем графа
$\Gamma$ относительно некоторой вершины
влечет наличие у $A_{\Gamma,F}$ собственной функции с
бесконечным носителем, соответствующей собственному значению $\lambda$.

Особое значение
пропагаторов для проблематики настоящей работы обусловлено тем, что, наряду с возможностью выразить
и исследовать в их терминах представляющие
интерес свойства оператора смежности графа (см. \S~\ref{s4} и \S~\ref{s5}), они (в отличие от
собственных функций операторов смежности) зачастую допускают конструктивную
реализацию с использованием определенных сумм по путям графа и (в случае связного графа с ограниченными
в совокупности степенями вершин и $F = \mathbb{C}$) с использованием ограничения $A^{({\rm alg})}_{\Gamma,\mathbb{C}}$
на пространство квадратично суммируемых функций. Этим конструктивным реализациям пропагаторов посвящен \S~\ref{s6}.
Следствием использования определенных сумм по путям графа для
построения пропагаторов в сочетании с предшествующими результатами работы является то, что (см. теорему~\ref{t6.1})
в случае бесконечного локально конечного связного графа $\Gamma$ и поля $F$ характеристики $0$ каждый трансцендентный над простым
подполем поля $F$ элемент $\lambda \in F$ является собственным значением оператора $A_{\Gamma,F}$.
Кроме того, в случае, когда степени вершин бесконечного связного графа $\Gamma$ не превосходят некоторого числа $d \in \mathbb{Z}_{\geq 1}$  и  $F = \mathbb{C}$,
с использованием ограничения $A^{({\rm alg})}_{\Gamma,\mathbb{C}}$
на пространство квадратично суммируемых функций
становится возможным (см. утверждение 2) теоремы~\ref{t6.2}) для произвольной вершины $v$ графа $\Gamma$ и произвольного
$\lambda \in  \mathbb{C} \setminus (\mathbb{R}_{\geq - d} \cap \mathbb{R}_{\leq d})$ указать
вид квадратично суммируемого $(\mathbb{C},\lambda)$-пропагатора графа $\Gamma$ относительно вершины $v$,
причем однозначно определенного указанным
свойством. Как следствие этой единственности, ни для какого
$\lambda \in \mathbb{C} \setminus (\mathbb{R}_{\geq - d} \cap \mathbb{R}_{\leq d})$
ни одна из собственных функций оператора $A_{\Gamma, \mathbb{C}}$, соответствующих $\lambda$,
не является квадратично суммируемой функцией (хотя согласно теореме~\ref{t6.1}
собственные функции $A_{\Gamma, \mathbb{C}}$, соответствующие $\lambda$, существуют
по крайней мере для всех трансцендентных над $\mathbb{Q}$ комплексных чисел $\lambda$).
В связи с этим можно сказать, что (в случае бесконечного связного графа $\Gamma$ со степенями вершин $\leq d$)
для всех $\lambda \in \mathbb{C} \setminus (\mathbb{R}_{\geq - d} \cap \mathbb{R}_{\leq d})$
пространство квадратично суммируемых функций замечательным
образом подходит для реализаций $(\mathbb{C},\lambda)$-пропагаторов, но слишком мало
для реализаций соответствующих таким $\lambda$ собственных функций оператора $A_{\Gamma, \mathbb{C}}$ (которые, однако, могут быть исследованы
с привлечением указанных пропагаторов, используя результаты
предшествующих параграфов).

В \S~\ref{s7} как следствие предыдущих результатов работы показывается, что для локально конечного графа
$\Gamma$ при выборе $\lambda \in F$ в определенном смысле ``типична''
ситуация наличия $(F,\lambda)$-пропагаторов относительно всех вершин графа
(и, таким образом, сюръективность $A_{\Gamma,F} - \lambda E$). С другой стороны, в случае бесконечного локально конечного связного графа $\Gamma$ и
поля $F$ характеристики $0$ показывается в определенном смысле ``исключительность'' для $\lambda \in F$
ситуации, когда $\lambda$ не есть собственное значение $A_{\Gamma,F}$ или, другими словами, см. выше,
у $\Gamma$ относительно каждой вершины имеется $(F,\lambda)$-пропагатор с конечным носителем. Более того,
показывается, что (см. предложения~\ref{p7.3} и~\ref{p7.2}) в случае бесконечного локально конечного связного графа $\Gamma$ и
поля $F$ характеристики $0$ в определенном смысле ``исключительной'' для $\lambda \in F$
является ситуация наличия у $\Gamma$ относительно хотя бы одной вершины $(F,\lambda)$-пропагатора с конечным носителем.

В \S~\ref{s8} приводятся примеры бесконечных локально конечных связных графов $\Gamma$ и полей $F$
с представляющими интерес в контексте настоящей работы свойствами
операторов $A_{\Gamma,F}$.

Используемые в работе обозначения, касающиеся графов, в основном стандартны. Для графа $\Gamma$
через $d_{\Gamma}(.,.)$ обозначается обычное расстояние между вершинами (равное $\infty$ для вершин
из разных компонент связности графа); для $u \in V(\Gamma)$ и $r \in \mathbb{R}_{\geq 0}$ через $B_{\Gamma}(u,r)$
обозначается множество вершин графа $\Gamma$, удаленных от $u$ на расстояние $\leq r$.
Если $\Gamma$ --- граф и $\emptyset \not = X \subseteq V(\Gamma)$, то
$\langle X \rangle_{\Gamma}$ --- подграф графа $\Gamma$, порожденный $X$. $Aut(\Gamma)$ ---  группа автоморфизмов графа
$\Gamma$ (рассматриваемая как группа подстановок на $V(\Gamma)$).

Для поля $F$ и непустых множеств $X$, $Y$ через $M_{X\times Y}(F)$ обозначается
векторное пространство над $F$ всех матриц с элементами из $F$ со строками, проиндексированными элементами
множества $X$, и столбцами, проиндексированными элементами
множества $Y$. При этом через $M^{rf}_{X\times Y}(F)$ обозначается подмножество (векторное подпространство)
$M_{X\times Y}(F)$, состоящее из
всех матриц из $M_{X\times Y}(F)$, у которых в каждой строке имеется лишь конечное число ненулевых элементов;
через $M^{cf}_{X\times Y}(F)$ обозначается подмножество (векторное подпространство)
$M_{X\times Y}(F)$, состоящее из
всех матриц из $M_{X\times Y}(F)$, у которых в каждом столбце имеется лишь конечное число ненулевых элементов;
$M^{rcf}_{X\times Y}(F) := M^{rf}_{X\times Y}(F) \cap M^{cf}_{X\times Y}(F)$.
Если $F$ --- поле и $X$, $Y$, $Z$ --- непустые множества, то для каждой матрицы из $M_{Y\times Z}(F)$
определено (естественным образом)
умножение на любую матрицу из $M^{rf}_{X\times Y}(F)$ слева,
а для каждой матрицы из $M_{X\times Y}(F)$ определено (естественным образом)
умножение на любую матрицу из $M^{cf}_{Y\times Z}(F)$
справа. Ясно также, что для поля $F$ и непустого множества $X$ векторные подпространства $M^{rf}_{X\times X}(F)$,
$M^{cf}_{X\times X}(F)$ и $M^{rcf}_{X\times X}(F)$ векторного пространства $M_{X\times X}(F)$ являются (относительно
естественным образом определенных операций) ассоциативными $F$-алгебрами.
Для поля $F$, непустых множеств $X$, $Y$ и произвольной матрицы ${\bf M} \in M_{X\times Y}(F)$ через ${\bf M}^t$ обозначается
транспонированная к ${\bf M}$ матрица из $M_{Y\times X}(F)$ (определяемая естественным образом).

Для графа $\Gamma$ и поля $F$ через ${\bf A}_{\Gamma,F}$ обозначается матрица из
$M_{V(\Gamma)\times V(\Gamma)}(F)$, у которой для произвольных $u,v \in V(\Gamma)$ элемент
на пересечении строки, соответствующей $u$, и столбца, соответствующего $v$, равен $1$ при $\{u,v\} \in E(\Gamma)$ и
равен $0$ при $\{u,v\} \not \in E(\Gamma)$. Ясно, что ${\bf A}_{\Gamma,F}^t = {\bf A}_{\Gamma,F}$.
Ясно также, что $\Gamma$ локально конечен тогда и только тогда, когда
${\bf A}_{\Gamma,F} \in M^{rcf}_{V(\Gamma)\times V(\Gamma)}(F)$.

Если $\Gamma$ --- локально конечный граф и $F$ --- поле, то векторное пространство (над $F$) функций $F^{V(\Gamma)}$
можно естественным образом отождествить с векторным пространством (над $F$) матриц-столбцов
$M_{V(\Gamma)\times \{1\}}(F)$, сопоставляя каждой функции $f \in F^{V(\Gamma)}$
столбец ${\bf f} \in M_{V(\Gamma)\times \{1\}}(F)$, у которого элемент строки, соответствующей
вершине $u$, равен $f(u)$ для всех $u \in V(\Gamma)$. Ясно, что при этом для произвольных
$\lambda \in F$ и $f \in F^{V(\Gamma)}$ функции $(A^{({\rm alg})}_{\Gamma,F} - \lambda E)(f) \in F^{V(\Gamma)}$
сопоставляется столбец  $({\bf A}_{\Gamma,F} - \lambda {\bf E})({\bf f}) \in M_{V(\Gamma)\times \{1\}}(F)$
(где ${\bf E}$ --- единичная матрица из $M_{V(\Gamma)\times V(\Gamma)}(F)$).

\section{Терминология и вспомогательные результаты}
\label{s2}

В настоящем параграфе, содержащем вспомогательные результаты, предполагается, что $\Gamma$ --- локально конечный граф, $F$ --- поле
с (возможно, тривиальным) абсолютным значением $|.|_{\rm v}$. При этом предполагается (см. введение),
что поле $F$
снабжено топологией (метрикой), определяемой $|.|_{\rm v}$, а векторное пространство $F^X$ для произвольного
$\emptyset \not = X \subseteq V(\Gamma)$ --- соответствующей топологией произведения.
В случае связного (и, более общо,  не более, чем счетного) графа $\Gamma$ эта топология произведения на $F^{V(\Gamma)}$
может быть, как хорошо известно, задана метрикой $\rho$, определяемой следующим образом:
произвольным образом нумеруем все вершины графа
$\Gamma$, получая $V(\Gamma) = \{v_n | n \in {\cal N}\}$, где ${\cal N}$ $= \{ 1,...,|V(\Gamma)|\}$, если $\Gamma$ конечен, и
${\cal N} = \mathbb{Z}_{\geq 1}$
в противном случае, после чего для произвольных $f_1, f_2 \in F^{V(\Gamma)}$ полагаем
$$\rho(f_1,f_2) = \sum_{n \in \cal N} \left ( \frac{1}{2^n}\cdot \frac{|f_1(v_n) - f_2(v_n)|_{\rm v}}{1 + |f_1(v_n) -
f_2(v_n)|_{\rm v}}\right ).$$
Важно отметить, что для многих дальнейших результатов топология на $F$ (и на $F^{V(\Gamma)}$)
играет лишь вспомогательную роль. В этой связи напомним, что как свойство
элемента $\lambda \in F$
быть собственным значением оператора $A_{\Gamma,F}$, так и свойство функции $f \in F^{V(\Gamma)}$
быть собственной функцией оператора $A_{\Gamma,F}$, соответствующей собственному значению
$\lambda \in F$, не зависят от выбора абсолютного значения $|.|_{\rm v}$. (Однако наделение поля $F$
абсолютным значением  $|.|_{\rm v}$ может оказаться полезным, например, для нахождения для оператора $A_{\Gamma,F}$ и $\lambda \in F$ собственной функции, соответствующей собственному значению
$\lambda$, или $(F,\lambda)$-пропагатора, поскольку делает возможным привлечение топологических соображений, см. в связи с этим замечание~\ref{r3.1} и \S~\ref{s6}.)

Если не оговорено противное, подполе поля $F$ с абсолютным значением
$|.|_{\rm v}$ предполагается снабженным абсолютным значением,
являющимся ограничением на это подполе абсолютного значения $|.|_{\rm v}$.
Векторное подпространство топологического векторного пространства $F^{V(\Gamma)}$ предполагается, если не
оговорено противное, снабженным индуцированной топологией.

Векторные подпространства топологического векторного пространства $F^{V(\Gamma)}$ не предполагаются замкнутыми.
Через ${\rm Cl}$ обозначается операция взятия замыкания в $F^{V(\Gamma)}$.
Для произвольного $\lambda \in F$ ясно, что $A_{\Gamma,F} - \lambda E$ --- (всюду определенный) непрерывный линейный оператор из $F^{V(\Gamma)}$
в $F^{V(\Gamma)}$.
Если $X \subseteq V(\Gamma)$, то для $f \in F^{V(\Gamma)}$ и ${\cal F} \subseteq F^{V(\Gamma)}$
через $f|_X$ и ${\cal F}|_X$ обозначаются ограничения соответственно $f$ и $\cal F$ на $X$ (т. е.  образы $f$ и $\cal F$
при естественном проектировании
$F^{V(\Gamma)} \to F^X$).  Для произвольного $\emptyset \not = X \subseteq V(\Gamma)$ мы полагаем $\Gamma(X) := \cup_{u \in X}\Gamma(u)$,
а через $(A_{\Gamma,F} - \lambda E)_X$ для $\lambda \in F$ обозначаем линейный оператор из $F^{X \cup \Gamma(X)}$ в $F^X$,
отображающий функцию $f \in F^{X \cup \Gamma(X)}$ в функцию из $F^X$, значение которой
в произвольной вершине $u \in X$
равно $(\sum_{u' \in \Gamma(u)}f(u')) - \lambda f(u)$.

Через $F^{V(\Gamma)*}$ обозначается подпространство векторного пространства $F^{V(\Gamma)}$, состоящее из всех функций с конечными носителями.
Для $v \in V(\Gamma)$ через $\delta_v$ обозначается (как уже было сказано во введении) функция из $F^{V(\Gamma)}$
такая, что $\delta_v(u) = 0$ для всех $u \in V(\Gamma) \setminus \{v\}$ и $\delta_v(v) = 1$.
Кроме того, для $v \in V(\Gamma)$ полагаем
$$\alpha_v := (\sum_{v' \in \Gamma(v)} \delta_{v'}) - \lambda \delta_v.$$
Таким образом, $\alpha_v = (A_{\Gamma,F} - \lambda E)(\delta_v)$
для всех $v \in V(\Gamma)$, а $(A_{\Gamma,F} - \lambda E)(F^{V(\Gamma)*})$
совпадает с линейной оболочкой в $F^{V(\Gamma)}$ множества $\{\alpha_v : v \in V(\Gamma)\}$.

Далее в работе нам иногда будет удобно привлекать терминологию, используемую при рассмотрении
(систем) линейных уравнений.
Как обычно, для поля $F$ линейным уравнением над $F$ с неизвестными $x_j$, $j \in J \not = \emptyset,$
называется уравнение вида
$$\sum_{j \in J} c_j x_j = a,$$
где $a \in F$ и $c_j \in F$ для каждого $j \in J$, причем множество $\{j \in J : c_{j} \not = 0\}$ конечно.
В обычном смысле понимаются решения в $F$ такого уравнения.
В обычном же смысле понимаются системы (не обязательно конечные) линейных уравнений над $F$,
их решения в $F$, а также их разрешимость в $F$ (как наличие решения в $F$).

Важную роль в дальнейшем будет играть следующая теорема Теплица
(см.~\cite{T}, а также~\cite[гл. 2, \S~ 1.4]{L}):

\begin{propos}
\label{p2.1}

Пусть $F$ --- поле, $I$ и $J$ --- множества $($не обязательно конечные$)$,
причем $J \not = \emptyset$.
Для произвольных $i \in I$ и $j \in J$ пусть $c_{i,j} \in F$, причем для каждого $i \in I$
множество $\{j \in J : c_{i,j} \not = 0\}$ конечно. Пусть, наконец, $a_i \in F$ для каждого $i \in I$.
Тогда система линейных уравнений
$$\sum_{j \in J} c_{i,j} x_j = a_i,\ i \in I,$$
с неизвестными $x_j, j \in J,$ разрешима над $F$ тогда и только тогда, когда
для произвольного конечного подмножества $I'$ множества $I$ и произвольных
$b_{i'} \in F,$ $i' \in I',$ из
$$\sum_{i' \in I'}b_{i'}c_{i',j} = 0 \  \mbox{для каждого}\ j \in J$$
следует
$$\sum_{i' \in I'}b_{i'}a_{i'} = 0.$$

\end{propos}

Теорему Теплица чаще нам будет  удобно использовать в следующем виде (см.~\cite[гл. 2, \S~ 1.4]{L}, а также~\cite{A};
учтено, кроме того, что конечная система линейных уравнений над полем разрешима над ним, если
она разрешима над каким-либо его расширением).

\begin{propos}
\label{p2.2}

Система $($не обязательно конечная$)$ линейных уравнений над полем $F$
разрешима над $F$ тогда и только тогда, когда каждая ее конечная подсистема
разрешима над некоторым расширением поля $F$.

\end{propos}

Из теоремы Теплица в последней форме непосредственно следует

\begin{propos}
\label{p2.3}

Пусть $\Gamma$ --- локально конечный граф, $F$ --- поле с
абсолютным значением $|.|_{\rm v}$ и $\lambda \in F$.
Тогда для $f \in F^{V(\Gamma)}$ включение $f \in  (A_{\Gamma,F} - \lambda E)(F^{V(\Gamma)})$
эквивалентно включениям $f|_X \in  (A_{\Gamma,F} - \lambda E)_X(F^{X \cup \Gamma(X)})$
для всех конечных непустых подмножеств $X$ множества $V(\Gamma)$.

\end{propos}

Следующее предложение демонстрирует некоторые ``технические'' возможности теоремы Теплица,
которые полезно иметь в виду.

\begin{propos}
\label{p2.4}

Пусть $\Gamma$ --- локально конечный граф, $F$ --- поле с
абсолютным значением $|.|_{\rm v}$ и $\lambda \in F$.
Пусть, далее, $U \subseteq V(\Gamma)$, и для каждой вершины $u \in U$ пусть $c_u \in F$.
Пусть, кроме того, $U'$ --- конечное подмножество множества $V(\Gamma)$, $I$ --- некоторое множество
и для каждого $i \in I$ заданы $u_i \in U'$, $c'_{i} \in F$, $q_{i} \in \mathbb{R}_{\geq 0}$, а также $\tau_{i}$ --- одно из бинарных отношений
$\not =,<,\leq,\geq,>$ на $\mathbb{R}$.
Тогда следующие условия $(i)$ и $(ii)$ эквивалентны:

$(i)$ Найдется $f \in F^{V(\Gamma)}$ такая, что  $(A_{\Gamma,F} - \lambda E)(f) = 0$,
$f(u) = c_u$ для каждой вершины $u \in U$ и $|f(u_i) - c'_{i}|_{\rm v} \tau_{i} q_{i}$ для каждого $i \in I$;

$(ii)$ Для каждого конечного подмножества $X \subseteq V(\Gamma)$ найдется $f_X \in F^{V(\Gamma)}$ такая, что  $((A_{\Gamma,F} - \lambda E)(f_X))|_X = 0$,
$f_X(u) = c_u$ для каждой вершины $u \in U \cap X$  и $|f_X(u_i) - c'_{i}|_{\rm v} \tau_{i} q_{i}$ для каждого $i \in I$.

\end{propos}

\begin{proof}

Предположим, что имеет место $(ii)$, и докажем, что тогда выполняется $(i)$. Ограничение на $U'$ множества функций $f \in F^{V(\Gamma)}$ таких, что
$(A_{\Gamma,F} - \lambda E)(f) = 0$ и
$f(u) = c_u$ для каждой вершины $u \in U$, либо пусто, либо имеет вид $h + W$, где $h \in F^{U'}$ и $W$ --- подпространство (конечномерного) векторного пространства $F^{U'}$.
Но тогда из теоремы Теплица
(см. предложение~\ref{p2.2}) следует, что это ограничение совпадает для некоторого конечного подмножества $X$ множества $V(\Gamma)$
с ограничением на $U'$ множества функций $f \in F^{V(\Gamma)}$ таких, что
$(A_{\Gamma,F} - \lambda E)(f)|_X = 0$ и
$f(u) = c_u$ для каждой вершины $u \in U \cap X$.
(Действительно, для каждого конечного подмножества $X'$ множества $V(\Gamma)$ ограничение на $U'$
множества функций $f \in F^{V(\Gamma)}$ таких, что
$(A_{\Gamma,F} - \lambda E)(f)|_{X'} = 0$ и
$f(u) = c_u$ для каждой вершины $u \in U \cap X'$, имеет вид $h + W_{X'}$, где $W_{X'}$ --- содержащее
$W$ подпространство пространства $F^{U'}$. Поскольку по теореме Теплица пересечение $W_{X'}$
по всем конечным подмножествам $X' \subseteq V(\Gamma)$ совпадает с $W$, а векторное пространство $F^{U'}$
конечномерно, для некоторого конечного набора $X_1,...,X_k$ конечных подмножеств множества $V(\Gamma)$
имеем $W = \cap_{1 \leq j \leq k}W_{X_j}$. Остается заметить, что
$W_{X_1 \cup ... \cup X_k} = \cap_{1 \leq j \leq k}W_{X_j}$, и следовательно, в качестве $X$ можно
взять множество $X_1 \cup ... \cup X_k$.)
В силу $(ii)$ отсюда следует $(i)$.

Поскольку $(ii)$ тривиальным образом следует из $(i)$, предложение доказано.

\end{proof}

В заключительной части параграфа покажем, что для локально конечного графа $\Gamma$, поля $F$
с абсолютным значением $|.|_{\rm v}$ и $\lambda \in F$ оператор $A_{\Gamma,F} - \lambda E$
отображает замкнутые векторные подпространства топологического векторного пространства $F^{V(\Gamma)}$ в замкнутые же подпространства
(см. следствие~\ref{c2.1}). Для этого установим справедливость предложения~\ref{p2.5},
при доказательстве которого используются по существу стандартные рассуждения с проективными пределами, а также замкнутость произвольного
векторного подпространства векторного пространства $F^X$ с топологией произведения, где $F$ --- поле с топологией, определяемой
некоторым абсолютным значением $|.|_{\rm v}$, и $X$ --- конечное непустое множество. (Замкнутость произвольного векторного подпространства такого векторного
пространства $F^X$ следует из
возможности его представления как множества всех $f \in F^X$ таких, что $\sum_{u \in X}c_{i,u}f(u) = 0$, $i = 1,...,|X|$, для некоторых
фиксированных для данного подпространства элементов $c_{i,u}$ поля $F$, $1 \leq i \leq |X|$, $u \in X$.)

\begin{propos}
\label{p2.5}

Пусть  $\Gamma$ --- локально конечный связный граф, $F$ --- поле
с абсолютным значением $|.|_{\rm v}$ и $\lambda \in F$.
Пусть, кроме того, $(f_i)_{i \in \mathbb{Z}_{\geq 1}}$ --- последовательность векторов из  $F^{V(\Gamma)}$,
для которой последовательность $((A_{\Gamma,F} - \lambda E)(f_i))_{i \in \mathbb{Z}_{\geq 1}}$ сходится к некоторому
$f \in F^{V(\Gamma)}$.  Тогда найдется такой вектор $f'$, принадлежащий замыканию в $F^{V(\Gamma)}$
векторного подпространства, порожденного множеством векторов $\{f_i : i \in \mathbb{Z}_{\geq 1}\}$, что
$(A_{\Gamma,F} - \lambda E)(f') = f$.

\end{propos}

\begin{proof}

Пусть $v \in V(\Gamma)$ и $B_i := B_{\Gamma}(v,i)$ для всех $i \in \mathbb{Z}_{\geq 0}$.
Пусть, кроме того, $W$ --- векторное подпространство пространства
$F^{V(\Gamma)}$,  порожденное множеством векторов $\{f_i : i \in \mathbb{Z}_{\geq 1}\}$.
Тогда, поскольку согласно условию предложения $f \in {\rm Cl}((A_{\Gamma,F} - \lambda E)(W))$, для каждого $i \in  \mathbb{Z}_{\geq 0}$
имеем
$f|_{B_i} \in ({\rm Cl}((A_{\Gamma,F} - \lambda E)(W)))|_{B_i}$, что с учетом равенства
$({\rm Cl}((A_{\Gamma,F} - \lambda E)(W)))|_{B_i} = ((A_{\Gamma,F} - \lambda E)(W))|_{B_i}$, которое следует
из замкнутости векторного подпространства  $((A_{\Gamma,F} - \lambda E)(W))|_{B_i}$ топологического
векторного пространства $F^{B_i}$ (см. замечание непосредственно предшествующее формулировке предложения~\ref{p2.5}), влечет
существование $\tilde f_i \in W$ со свойством $((A_{\Gamma,F} - \lambda E)(\tilde f_{i}))|_{B_{i}} = f|_{B_i}$.

Для каждого $i \in  \mathbb{Z}_{\geq 0}$ положим $L_i :=  (A_{\Gamma,F} - \lambda E)_{B_i}$
и ${\cal F}_i := W|_{B_{i+1}} \cap (L_i)^{-1}(f|_{B_i})$ (где $(L_i)^{-1}(f|_{B_i})$ --- множество
всех функций из $F^{B_{i+1}}$, которые под действием $L_i$ отображаются в $f|_{B_i}$). Заметим,
что ${\cal F}_i$ --- непустое (в силу $\tilde f_{i+1}|_{B_{i+1}} \in {\cal F}_i$) аффинное подпространство аффинного пространства, ассоциированного
с векторным пространством $F^{B_{i+1}}$.
Кроме того, для произвольных целых неотрицательных чисел $i_1 \leq i_2$ отображение ограничения
$F^{B_{i_2+1}} \to F^{B_{i_1+1}}$ индуцирует аффинное отображение $\pi_{i_1,i_2}: {\cal F}_{i_2} \to {\cal F}_{i_1}$.
Получаемая таким образом проективная система множеств ${\cal F}_i$, $i \in \mathbb{Z}_{\geq 0}$,
и отображений $\pi_{i_1,i_2}$, $i_1,i_2 \in \mathbb{Z}_{\geq 0}$, $i_2 \geq i_1$,
имеет проективный предел (в силу наличия у нее так называемого свойства Миттаг-Леффлера:
для каждого $i \in \mathbb{Z}_{\geq 0}$ последовательность
$(\pi_{i,i+j}({\cal F}_{i+j}))_{j \in \mathbb{Z}_{\geq 0}}$ стабилизируется).
Каждый элемент этого проективного предела можно интерпретировать как такой вектор $f' \in  F^{V(\Gamma)}$, что
для любого $i \in \mathbb{Z}_{\geq 1}$ найдется вектор $f'_i \in W$, для которого $f'_i|_{B_i} \in  {\cal F}_{i-1}$ и $f'_i|_{B_i} = f'|_{B_i}$.
Ясно, что при этом последовательность  $(f'_i)_{i \in \mathbb{Z}_{\geq 1}}$ векторов из $W$ сходится к вектору $f'$ (который, следовательно,
принадлежит ${\rm Cl}(W)$) и $(A_{\Gamma,F} - \lambda E)(f') = f$.

\end{proof}

\begin{cor}
\label{c2.1}

Пусть  $\Gamma$ --- локально конечный граф, $F$ --- поле
с абсолютным значением $|.|_{\rm v}$ и $\lambda \in F$.
Если $W$ --- подпространство  векторного пространства $F^{V(\Gamma)}$, то
$${\rm Cl}((A_{\Gamma,F} - \lambda E)(W)) = (A_{\Gamma,F} - \lambda E)({\rm Cl}(W)).$$
Другими словами, $A_{\Gamma,F} - \lambda E$ отображает замкнутые подпространства
векторного пространства $F^{V(\Gamma)}$ в замкнутые же подпространства.

\end{cor}

\begin{proof}
 Включение
$$ (A_{\Gamma,F} - \lambda E)({\rm Cl}(W)) \subseteq {\rm Cl}((A_{\Gamma,F} - \lambda E)(W))$$
следует из непрерывности  $A_{\Gamma,F} - \lambda E$.
Обратное включение следует из предложения~\ref{p2.5}, примененного к связным компонентам графа $\Gamma$.

\end{proof}

\begin{remark}
\label{r2.1}

Утверждение следствия, вообще говоря,
не имеет места для подмножества $W$ векторного пространства $F^{V(\Gamma)}$,
не являющегося его векторным подпространством.

\end{remark}

\section{Некоторые общие свойства собственных функций операторов смежности
локально конечных графов}
\label{s3}

Настоящий параграф содержит предварительное рассмотрение собственных функций операторов смежности
локально конечных графов.
Результаты параграфа носят общий, но отчасти инструментарный характер.

Пусть $\Gamma$ --- локально конечный граф, $F$ --- поле
с абсолютным значением $|.|_{\rm v}$ и  $\lambda \in F$.
Как уже было отмечено, свойство функции из $V(\Gamma)$ в $F$ быть собственной функцией оператора смежности $A_{\Gamma,F}$, соответствующей собственному значению
$\lambda$, не зависит от выбора $|.|_{\rm v}$. Кроме того,
из предложения~\ref{p2.2} легко следует, что существование собственной функции оператора смежности $A_{\Gamma,F}$, соответствующей собственному значению
$\lambda$, равносильно существованию собственной функции оператора смежности $A_{\Gamma,F'}$, соответствующей собственному значению $\lambda$, для какого-либо
расширения $F'$ поля $F_0(\lambda)$, где $F_0$ --- простое подполе поля $F$. Мы зафиксируем эти факты (позволяющие для построения собственных функций
изменять удобным образом поля и снабжать их подходящими абсолютными значениями) в несколько уточненном виде
в следующем предложении.

\begin{propos}
\label{p3.1}

Пусть $\Gamma$ --- локально конечный граф, $F$ --- поле
с абсолютным значением $|.|_{\rm v}$, $F_0$ --- простое подполе поля $F$ и $\lambda \in F$.
Тогда справедливы следующие утверждения$:$

$1)$ Если $\varphi$ --- вложение поля $F$ в некоторое поле $F'$ с произвольным абсолютным
значением $|.|_{{\rm v}'}$
и $f$ --- собственная функция оператора смежности $A_{\Gamma,F'}$, соответствующая собственному значению $\varphi(\lambda)$, такая, что $f(u) \in \varphi(F)$ для всех
$u \in V(\Gamma)$, то $\varphi^{-1} f$, рассматриваемая как элемент
$F^{V(\Gamma)}$, является собственной функцией оператора смежности $A_{\Gamma,F}$, соответствующей собственному значению $\lambda$.

$2)$ Если $f$ --- собственная функция оператора смежности $A_{\Gamma,F}$, соответствующая собственному значению $\lambda$, и $v$ --- вершина графа $\Gamma$,
принадлежащая носителю $f$, то найдется собственная функция оператора смежности $A_{\Gamma,F_0(\lambda)}$, соответствующая собственному значению $\lambda$,
носитель которой содержится в носителе $f$ и содержит $v$.

\end{propos}

\begin{proof}

Утверждение $1)$ очевидно. Для доказательства утвер\-ждения $2)$
зададим следующим образом систему линейных уравнений над $F$ с неизвестными $x_u$, где $u \in V(\Gamma)$:
система состоит, во-первых,  из уравнений
$(\sum_{u' \in \Gamma(u)}x_{u'}) - \lambda x_u = 0$ для всех $u \in V(\Gamma)$,
во-вторых, из уравнений $x_u = 0$ для всех не принадлежащих носителю $f$ вершин $u$ графа $\Gamma$,
в-третьих, из уравнения  $x_v = 1$. Заметим, что все коэффициенты этой системы принадлежат $F_0(\lambda)$,
причем в силу наличия функции $f$ эта система совместна над $F$. В силу предложения~\ref{p2.2} отсюда следует, что
эта система совместна и над $F_0(\lambda)$.
Пусть $x'_u \in F_0(\lambda)$, где $u$ пробегает множество $V(\Gamma)$, --- некоторое
ее решение над $F_0(\lambda)$ (в том смысле, что при замене $x_u$ на $x'_u$ для всех $u \in V(\Gamma)$
каждое уравнение системы становится равенством в $F_0(\lambda)$). Тогда функция из $F_0(\lambda)^{V(\Gamma)}$, значение которой в каждой вершине $u \in V(\Gamma)$
равно $x'_u$, очевидным образом является требуемой
собственной функцией оператора смежности $A_{\Gamma,F_0(\lambda)}$, соответствующей собственному значению $\lambda$, носитель которой содержится
в носителе $f$ и содержит $v$.

\end{proof}

\begin{remark}
\label{r3.1}

Укажем на некоторые возможные использования предложения~\ref{p3.1}.
Предположим, что выполнены условия предложения~\ref{p3.1}, причем $\Gamma$ бесконечен и связен,
а характеристика поля $F$ равна $0$ (т. е.
$F_0 = \mathbb{Q}$), и нас интересует вопрос: является ли $\lambda$ собственным значением
оператора смежности $A_{\Gamma,F}$? В теореме~\ref{t6.1} будет доказано, что это так в случае, когда элемент $\lambda$
трансцендентен над $\mathbb{Q}$. Предположим поэтому, что $\lambda$ --- алгебраический над $\mathbb{Q}$
элемент. Тогда из предложения~\ref{p3.1} следует, что наличие собственной функции оператора смежности $A_{\Gamma,F}$,
соответствующей собственному значению $\lambda$,
равносильно наличию таковой для оператора $A_{\Gamma,\tilde F}$,
где в качестве $\tilde F$ можно взять любое расширение поля $\mathbb{Q}(\lambda)$.
Например, в качестве $\tilde F$ можно взять любое из следующих полей: $\mathbb{Q}(\lambda)$ (с произвольным
абсолютным значением), $\mathbb{C}$ (рассматриваемое как расширение $\mathbb{Q}(\lambda)$,
с произвольным абсолютным значением и, в частности, с обычным абсолютным значением на $\mathbb{C}$,
что позволяет привлекать для построения интересующих собственных функций аналитические методы, см. \S~\ref{s6}),
$\Omega_l$ для простого числа $l$ ($l$-адическое пополнение алгебраического замыкания поля
$l$-адических чисел $\mathbb{Q}_l$,
рассматриваемое как расширение $\mathbb{Q}(\lambda)$,
с произвольным абсолютным значением и, в частности, c продолжением $l$-адического абсолютного
значения $|.|_l$ на $\mathbb{Q}$, что также позволяет привлекать для построения собственных функций
аналитические методы). Кроме того, при рассмотрении, например, случая $\tilde F = \mathbb{C}$
с обычным абсолютным значением есть возможность (в силу утверждения 1) предложения~\ref{p3.1}) распорядиться
удобным образом выбором
вложения $\mathbb{Q}(\lambda)$ в $\mathbb{C}$.

\end{remark}

Далее, теорема Теплица (см. предложение~\ref{p2.2}) очевидным образом влечет

\begin{propos}
\label{p3.2}

Пусть $\Gamma$ --- локально конечный граф, $F$ --- поле
с абсолютным значением $|.|_{\rm v}$ и  $\lambda \in F$.
Пусть, кроме того, $v$ --- вершина графа $\Gamma$, не содержащаяся в носителях собственных функций оператора
смежности $A_{\Gamma,F}$, принадлежащих собственному значению $\lambda$.
Для каждого $r \in \mathbb{Z}_{\geq 0}$ определим систему $Syst_r$ линейных уравнений над $F$
относительно неизвестных $x_u$, где $u \in B_{\Gamma}(v,r + 1)$, включив в нее уравнение $x_v = 1$
и для каждого $u \in B_{\Gamma}(v,r)$ уравнение $(\sum_{u' \in \Gamma(u)}x_{u'}) - \lambda x_u = 0$.
Тогда найдется такое $r_{\lambda,v} \in \mathbb{Z}_{\geq 0}$ $($$``$радиус несовместности для $\lambda$ и $v$$"$$)$,
что система $Syst_r$ совместна $($над $F$$)$ тогда и только тогда, когда $r < r_{\lambda,v}$. При этом
число $r_{\lambda,v}$ не изменится, если заменить $F$ на $F_0(\lambda)$, где $F_0$ --- простое подполе
поля $F$.

\end{propos}

\begin{remark}
\label{r3.2}

Для локально конечного графа $\Gamma$, поля $F$ с абсолютным значением $|.|_{\rm v}$ и  $\lambda \in F$
множество всех
вершин графа $\Gamma$, не содержащихся ни в одном из носителей собственных функций оператора
смежности $A_{\Gamma,F}$, принадлежащих собственному значению $\lambda$, будет играть важную
роль в дальнейшем. В \S~\ref{s4} (см. теорему~\ref{t4.2}) будет показано, что это множество совпадает
с вводимым в \S~\ref{s4} множеством
$L_{F,\lambda}(\Gamma)$ всех вершин графа $\Gamma$, относительно которых у
$\Gamma$ имеются $(F,\lambda)$-пропагаторы с конечными носителями.

\end{remark}

Сделаем еще несколько наблюдений, касающихся носителей собственных функций операторов смежности
локально конечных графов.

\begin{propos}
\label{p3.3}

Пусть $\Gamma$ --- локально конечный граф, $F$ --- поле с
абсолютным значением $|.|_{\rm v}$, $\lambda \in F$ и $v \in V(\Gamma)$.
Если $X_i$, $i \in I$, --- непустое линейно упорядоченное по включению семейство подмножеств множества
$V(\Gamma)$, содержащих $v$ и являющихся носителями собственных функций оператора
$A_{\Gamma,F}$, соответствующих собственному значению $\lambda$,
то некоторое  подмножество множества $\cap_{i \in I}X_i$ является
содержащим $v$ носителем собственной функции оператора
$A_{\Gamma,F}$, соответствующей собственному значению $\lambda$.
Как следствие этого $($и леммы Цорна$)$, каждый содержащий $v$ носитель собственной функции оператора
$A_{\Gamma,F}$, соответствующей собственному значению $\lambda$, содержит
минимальный по включению среди содержащих $v$ носителей собственных функций оператора
$A_{\Gamma,F}$, соответствующих собственному значению $\lambda$.

\end{propos}

\begin{proof}

Утверждение следует из теоремы Теплица (см. предложение~\ref{p2.2}),
примененной к системе линейных уравнений над $F$ с неизвестными $x_u$, где $u \in V(\Gamma)$,
состоящей, во-первых, из уравнений
$(\sum_{u' \in \Gamma(u)}x_{u'}) - \lambda x_u = 0$ для всех $u \in V(\Gamma)$,
во-вторых, из уравнений $x_u = 0$ для всех $u \in V(\Gamma) \setminus \cap_{i \in I}X_i$,
в-третьих, из уравнения  $x_v = 1$.

\end{proof}

Если $\Gamma$ --- локально конечный граф, $F$ --- поле с
абсолютным значением $|.|_{\rm v}$, $\lambda \in F$ и $f$ --- собственная функция оператора
$A_{\Gamma,F}$, соответствующая собственному значению $\lambda$, то для $\emptyset \not = X \subseteq V(\Gamma)$, являющегося
объединением набора связных компонент носителя функции $f$, функция из $F^{V(\Gamma)}$, совпадающая с $f$ на $X$
и тождественно равная нулю на $V(\Gamma) \setminus X$, не является, вообще говоря, собственной функцией оператора
$A_{\Gamma,F}$, соответствующей собственному значению $\lambda$ (хотя ограничение $f$ на $X$ очевидным
образом является собственной функцией оператора
$A_{\langle X \rangle_{\Gamma},F}$, соответствующей собственному значению $\lambda$).
Однако такого рода эффект не возникает,
если вместо связных компонент носителя рассматривать вводимые ниже $\Gamma^2$-связные компоненты носителя.

Для произвольного графа $\Gamma$ будем обозначать через $\Gamma^2$ граф  с множеством вершин $V(\Gamma^2) = V(\Gamma)$
и множеством ребер, состоящим в точности из таких пар его различных вершин $\{u,u'\}$, что  $u$ и $u'$ соединены в графе $\Gamma$ путем длины $\leq 2$.
Пусть $X$ --- подмножество множества вершин графа $\Gamma$. Скажем, что $X$ является $\Gamma^2$-связным, если связен порожденный
$X$ подграф $\langle X \rangle_{\Gamma^2}$ графа $\Gamma^2$. {\it Компонентами $\Gamma^2$-связности} (или {\it $\Gamma^2$-связными компонентами}) $X$
будем называть множества вершин связных компонент
графа $\langle X \rangle_{\Gamma^2}$. Если $\Gamma$ --- локально конечный граф, $F$ --- поле
с абсолютным значением $|.|_{\rm v}$, $\lambda \in F$ и $f$ --- собственная функция оператора
$A_{\Gamma,F}$, соответствующая собственному значению $\lambda$, то для $\emptyset \not = X \subseteq V(\Gamma)$, являющегося
объединением произвольного набора  $\Gamma^2$-связных компонент носителя функции $f$, функция из $F^{V(\Gamma)}$, совпадающая с $f$ на $X$
и тождественно равная нулю на $V(\Gamma) \setminus X$, очевидным образом также является собственной функцией оператора
$A_{\Gamma,F}$, соответствующей собственному значению $\lambda$, но с носителем $X$. В частности,
для обсуждавшихся ранее (см. предложение~\ref{p3.3}) минимальных по включению среди содержащих $v$ носителей собственных функций оператора
$A_{\Gamma,F}$, соответствующих собственному значению $\lambda$, справедливо следующее утверждение.

\begin{propos}
\label{p3.4}

Пусть $\Gamma$ --- локально конечный граф, $F$ --- поле с
абсолютным значением $|.|_{\rm v}$, $\lambda \in F$ и $v \in V(\Gamma)$.
Тогда каждый минимальный по включению среди содержащих $v$ носителей собственных функций оператора
$A_{\Gamma,F}$, соответствующих собственному значению $\lambda$, $\Gamma^2$-связен.

\end{propos}

\begin{remark}
\label{r3.3}

С учетом утверждения 2) предложения~\ref{p3.1} для локально конечного графа $\Gamma$, поля $F$
с абсолютным значением $|.|_{\rm v}$, $\lambda \in F$ и $v \in V(\Gamma)$
каждый минимальный содержащий $v$ носитель собственной функции
оператора $A_{\Gamma,F}$, соответствующей собственному значению $\lambda$, не только $\Gamma^2$-связен
(см. предложение~\ref{p3.4}), но и является носителем содержащейся в $F_0(\lambda)^{V(\Gamma)}$, где $F_0$ --- простое подполе поля $F$, собственной функции оператора
$A_{\Gamma,F}$, соответствующей собственному значению $\lambda$
(или другими словами, является носителем собственной функции оператора $A_{\Gamma,F_0(\lambda)}$, соответствующей собственному значению $\lambda$).

\end{remark}

В дальнейшем нам потребуется следующее

\begin{propos}
\label{p3.5}

Пусть $\Gamma$ --- локально конечный граф, $F$ --- поле с
абсолютным значением $|.|_{\rm v}$ и $\lambda \in F$.
Предположим, что имеется такое семейство $X_i, i \in I,$ конечных подмножеств множества $V(\Gamma)$, что
$\cup_{i \in I} X_i$ бесконечно и для каждого $i \in I$ имеется собственная функция $f_i$ оператора $A_{\Gamma,F}$,
соответствующая собственному значению $\lambda$, с носителем $X_i$. Тогда для каждой вершины из $\cup_{i \in I} X_i$
найдется собственная функция оператора $A_{\Gamma,F}$, соответствующая собственному значению $\lambda$,
с бесконечным носителем, содержащим эту вершину и содержащимся в $\cup_{i \in I} X_i$.

\end{propos}

\begin{proof}

Для $i \in I$ компонента $\Gamma^2$-связности множества $X_i$ является носителем
собственной функции оператора $A_{\Gamma,F}$, соответствующей собственному значению $\lambda$, которая
на вершинах из этой компоненты принимает те же значения, что и $f_i$, а на всех остальных вершинах графа $\Gamma$
равна нулю. Поэтому, переходя в случае необходимости от семейства множеств $X_i, i \in I,$ к семейству
всех их компонент $\Gamma^2$-связности, будем, не теряя общности, предполагать $\Gamma^2$-связным каждое из множеств
$X_i, i \in I.$ Кроме того, будем, не теряя общности, предполагать, что все множества $X_i, i \in I,$
попарно различны.

Положим $X := \cup_{i \in I} X_i$. Для доказательства предложения достаточно доказать существование собственной
функции оператора $A_{\Gamma,F}$, соответствующей собственному значению $\lambda$,
с бесконечным носителем, содержащимся в $X$. Действительно, если $f$ --- такая функция, то для каждой
вершины из $X_i, i \in I,$ функция $f$ или функция $f + f_i$ обладает, очевидно, требуемым в предложении
свойством.

Так как конечные подмножества $X_i, i \in I,$ образуют покрытие бесконечного множества $X$,
то справедливо по меньшей мере одно из следующих двух утверждений:

1) для некоторого бесконечного подмножества $J \subseteq I$ множества $X_j,$ $j \in J,$ попарно дизъюнктны;

2) для некоторой вершины $v \in X$ множество $\{i \in I: v \in X_i\}$ бесконечно.

Если справедливо утверждение 1), то с учетом локальной конечности графа $\Gamma$
найдется такое бесконечное подмножество $J' \subseteq J$, что
для произвольных различных $j_1, j_2 \in J'$ имеем $(X_{j_1} \cup \Gamma(X_{j_1})) \cap (X_{j_2} \cup \Gamma(X_{j_2})) = \emptyset$.
Но тогда
функция $f \in F^{V(\Gamma)}$ такая, что $f(u) = 0$ при
$u \in V(\Gamma) \setminus \cup_{j \in J'}X_j$ и $f(u) = f_{j}(u)$ при $u \in X_{j}$, где $j \in J'$,
очевидным образом является собственной функцией оператора $A_{\Gamma,F}$,
соответствующей собственному значению $\lambda$, с бесконечным носителем $\cup_{j \in J'}X_j$, содержащимся в
$X$.

Предположим поэтому, что справедливо утверждение 2).
Тогда с учетом локальной конечности графа $\Gamma$ и $\Gamma^2$-связности каждого из множеств
$X_i, i \in I,$ существует такой бесконечный путь $u_0 = v, u_1, u_2,...$ графа $\Gamma^2$ с попарно различными вершинами,
что для любого целого положительного числа $n$ найдется $i \in I$ со свойством $\{u_0,...,u_n\} \subseteq X_{i}$.
Ввиду конечности множеств $X_i, i \in I,$ очевидно, далее, наличие такой бесконечной возрастающей последовательности $n_0 = 0,n_1,n_2,...$ целых неотрицательных чисел,
что для каждого $k \in \mathbb{Z}_{\geq 1}$ найдется $i_k \in I$, для которого
$u_{n_{k-1}} \in X_{i_k}$ и $\{u_{n_k},u_{n_{k+1}},u_{n_{k+2}},...\} \cap X_{i_k} = \emptyset$.

Зададим следующим образом систему линейных уравнений над полем $F$ с неизвестными $x_u$, где $u \in V(\Gamma)$:
система состоит, во-первых, из уравнений
$(\sum_{u' \in \Gamma(u)}x_{u'}) - \lambda x_u = 0$ для всех $u \in V(\Gamma)$,
во-вторых, из уравнений $x_u = 0$ для всех $u \in V(\Gamma) \setminus X$,
в-третьих, из уравнений  $x_{u_{n_k}} = 1$ для всех $k \in \mathbb{Z}_{\geq 0}$.
В силу наличия функций $f_{i_1},f_{i_2},f_{i_3},...$ каждая конечная подсистема этой системы совместна над $F$.
(В качестве решения можно взять значения в вершинах подходящей линейной комбинации
функций $f_{i_1},f_{i_2},f_{i_3},...$ Заметим, что для каждого $k \in \mathbb{Z}_{\geq 1}$
по выбору $f_{i_k}$ имеем $(\sum_{u' \in \Gamma(u)}f_{i_k}(u')) - \lambda f_{i_k}(u) = 0$ для всех $u \in V(\Gamma)$,
$f_{i_k}(u) = 0$ для всех $u \in V(\Gamma) \setminus X$, $f_{i_k}(u_{n_{k-1}}) \not = 0$ и
$f_{i_k}(u_{n_l}) = 0$ для всех $l \in \mathbb{Z}_{\geq k}$.) Следовательно, по теореме Теплица
(см. предложение~\ref{p2.2}) совместна над $F$ и вся система.
Пусть $x'_u \in F$, где $u$ пробегает множество $V(\Gamma)$, --- некоторое
ее решение над $F$. Тогда, определяя $f \in F^{V(\Gamma)}$ посредством $f(u) = x'_u$ для всех  $u \in V(\Gamma)$, мы получает
собственную функцию оператора $A_{\Gamma,F}$, соответствующую собственному значению $\lambda$, с бесконечным носителем,
содержащимся в $X$, что завершает доказательство предложения~\ref{p3.5}.

\end{proof}

Еще одним следствием теоремы Теплица является следующее

\begin{propos}
\label{p3.6}

Пусть $\Gamma$ --- локально конечный граф и $F$ --- поле с таким абсолютным значением
$|.|_{\rm v}$, что определяемая этим абсолютным значением метрика на $F$ полна и не является дискретной.
Тогда для $\lambda \in F$ следующие условия $(i)$ и $(ii)$ эквивалентны$:$

$(i)$ Найдется собственная функция оператора смежности $A_{\Gamma,F}$, соответствующая собственному значению $\lambda$,
носитель которой совпадает с $V(\Gamma)$.

$(ii)$ Для каждых конечного подмножества $X \subseteq V(\Gamma)$ и вершины $u \in X$ найдется $f_{X,u} \in F^{V(\Gamma)}$ такая, что
$((A_{\Gamma,F} - \lambda E)(f_{X,u}))|_X = 0$
и $f_{X,u}(u)  \not = 0$.

\end{propos}

\begin{proof}

Доказывая предложение~\ref{p3.6}, мы будем без потери общности предполагать граф $\Gamma$ связным
(поскольку очевидно, что в общем случае предложение достаточно доказать для компонент связности графа $\Gamma$)
и бесконечным (поскольку для конечного графа $\Gamma$ в $(ii)$ можно положить $X = V(\Gamma)$, а
векторное пространство $W$ собственных функций оператора смежности $A_{\Gamma,F}$, соответствующих собственному значению $\lambda$, в силу бесконечности поля $F$ не может быть теоретико-множественным объединением
по всем $u \in V(\Gamma)$
собственных подпространств $W_{u}$, состоящих для каждой вершины $u \in V(\Gamma)$
из равных $0$ в $u$ собственных функций оператора $A_{\Gamma,F}$, соответствующих собственному значению $\lambda$).

Пусть имеет место $(ii)$.
Занумеруем все вершины графа $\Gamma$ положительными целыми числами, полагая
$V(\Gamma) = \{v_1,v_2,...\}$, и определим (см. ниже) функции
$f_i \in F^{V(\Gamma)}$, $i \in \mathbb{Z}_{\geq 1}$, обладающие следующими двумя свойствами:

\smallskip

$1)$  $(A_{\Gamma,F} - \lambda E)(f_i) = 0$ для каждого  $i \in \mathbb{Z}_{\geq 1}$;

$2)$ для каждого $j \in \mathbb{Z}_{\geq 1}$ последовательность
$f_1(v_j), f_2(v_j), ...$
сходится (в метрике, определяемой на $F$ абсолютным значением $|.|_{\rm v}$) к некоторому ненулевому элементу $c_j$ поля $F$.

\smallskip

\noindent
Ясно, что тогда функция $f \in F^{V(\Gamma)}$, определяемая посредством
$f(v_j) = c_j$ для всех $j \in \mathbb{Z}_{\geq 1}$, будет обладать требуемыми в $(i)$ свойствами.

Функции $f_i$, $i \in \mathbb{Z}_{\geq 1}$, определим индукцией по $i$. Пусть $i \in \mathbb{Z}_{\geq 1}$ и функции $f_{i'}$ уже определены
для всех $i' \in \mathbb{Z}_{\geq 1}$ со свойством $i' < i$. Из $(ii)$ и теоремы Теплица
(см. предложение~\ref{p2.2}) следует существование $\tilde f_i \in F^{V(\Gamma)}$ такой,
что $(A_{\Gamma,F} - \lambda E)(\tilde f_i) = 0$ и $\tilde f_i (v_i) = 1$ (поскольку выполнение этих условий есть разрешимость
системы линейных уравнений, каждая конечная подсистема которой разрешима в силу $(ii)$).
При $i = 1$ полагаем $f_1 = \tilde f_1$. При $i > 1$ полагаем $f_i := f_{i-1} + \varepsilon_i \tilde f_i$, где
$\varepsilon_i \not = 0$ таково, что, во-первых, $|\varepsilon_i|_{\rm v} |\tilde f_i(v_{i'})|_{\rm v} \leq 3^{-i}|f_{i-1}(v_{i'})|_{\rm v}$ для всех $i' < i$
и, во-вторых,  $|\varepsilon_i|_{\rm v} \leq 3^{-i}|f_{i-1}(v_{i})|_{\rm v}$ при  $f_{i-1}(v_{i}) \not = 0.$
Легко убедиться, что это определение корректно, а так определенные функции $f_i \in F^{V(\Gamma)}$,
$i \in \mathbb{Z}_{\geq 1},$ обладают свойствами $1)$ и $2)$.

Поскольку $(ii)$ тривиальным образом следует из $(i)$, предложение доказано.

\end{proof}

\section{Пропагаторы}
\label{s4}

Если $\Gamma$ --- локально конечный граф, $F$ --- поле с абсолютным значением $|.|_{\rm v}$, $\lambda \in F$ и
$v \in V(\Gamma)$, то {\it $(F,\lambda)$-пропагатором}
графа $\Gamma$ относительно вершины $v$ мы называем любую функцию $f \in F^{V(\Gamma)}$ такую, что
$(A_{\Gamma,F} - \lambda E)(f) = \delta_v$, где, напомним, $\delta_v \in F^{V(\Gamma)}$, $\delta_v(v) = 1$ и $\delta_v(u) = 0$ для всех
$u \in V(\Gamma) \setminus \{v\}$.
Заметим, что свойство функции $f \in F^{V(\Gamma)}$
быть $(F,\lambda)$-пропагатором
графа $\Gamma$ относительно вершины $v$ не зависит от выбора абсолютного значения $|.|_{\rm v}$.
Очевидно, что $(F,\lambda)$-пропагаторы
локально конечного графа $\Gamma$ относительно вершины $v$ --- это в точности те $f \in F^{V(\Gamma)}$,
для которых выполняется условие: если $X$ --- связная компонента графа $\Gamma$, то
$f|_X$ является $(F,\lambda)$-пропагатором графа $\langle X \rangle_{\Gamma}$ относительно вершины $v$
в случае $v \in X$ и является либо нулевой функцией, либо собственной функцией
оператора $A_{\langle X \rangle_{\Gamma},F}$, соответствующей собственному значению $\lambda$, в случае
$v \not \in X$.

\begin{remark}
\label{r4.1}

$(F,\lambda)$-пропагатор графа $\Gamma$ относительно вершины $v$ есть, по существу, функция
Грина (или фундаментальное решение) относительно $v$ для оператора $A_{\Gamma,F} - \lambda E$  и как таковая обладает рядом ``ожидаемых'' свойств.
Однако специфика пространства,
на котором задан оператор $A_{\Gamma,F} - \lambda E$, не позволяет использовать здесь многие ``стандартные'' приемы работы с функциями Грина.

\end{remark}

В этом и особенно следующем параграфах будет показано, что ряд представляющих интерес свойств
оператора $A_{\Gamma,F} - \lambda E$
(где $\Gamma$ --- локально конечный граф, $F$ --- поле с абсолютным значением $|.|_{\rm v}$, $\lambda \in F$)
и, в частности, свойство инъективности (или, другими словами, свойство $\lambda$ не быть собственным значением
$A_{\Gamma,F}$) могут быть сформулированы в терминах $(F,\lambda)$-пропагаторов графа $\Gamma$.
При этом, как будет видно из дальнейшего (см. особенно \S~\ref{s6}), исследовать $(F,\lambda)$-пропага\-торы графа $\Gamma$
зачастую проще, чем исследовать непосредственно сами эти представляющие интерес свойства и (что особо значимо для тематики настоящей работы) чем
исследовать непосредственно собственные функции оператора $A_{\Gamma,F}$,
соответствующие собственному значению $\lambda$.

Пусть $\Gamma$ --- локально конечный граф и $F$ --- поле  с абсолютным значением $|.|_{\rm v}$.
Вообще говоря, для заданных $\lambda \in F$ и $v \in V(\Gamma)$ у графа $\Gamma$ относительно $v$ может не быть $(F,\lambda)$-пропагатора (как, например,
для $\lambda = 0$ в случае $V(\Gamma) = \{v\}$), а в случае наличия такового он может не быть единственным (соответствующие многочисленные
примеры будут возникать далее в работе). Заметим, что для $\lambda \in F$
разность двух различных $(F,\lambda)$-пропагаторов графа $\Gamma$
относительно одной и той же вершины есть собственная функция оператора $A_{\Gamma,F}$, соответствующая собственному значению $\lambda$.
Обратно, для $\lambda \in F$ сумма $(F,\lambda)$-пропагатора графа $\Gamma$
относительно некоторой его вершины и собственной функции оператора $A_{\Gamma,F}$, соответствующей собственному значению $\lambda$,
есть другой $(F,\lambda)$-пропагатор графа $\Gamma$ относительно той же вершины.

Для локально конечного графа $\Gamma$, поля $F$ с абсолютным значением $|.|_{\rm v}$,
$\lambda \in F$  и  $v \in V(\Gamma)$, как уже было отмечено, свойство функции из $V(\Gamma)$ в $F$ быть
$(F,\lambda)$-пропагатором графа $\Gamma$ относительно $v$
не зависит от выбора $|.|_{\rm v}$. Кроме того (ср. начало \S~\ref{s3}),
используя предложение~\ref{p2.2}, легко доказать, что существование
$(F,\lambda)$-пропагатора графа $\Gamma$ относительно $v$
равносильно существованию
$(F',\lambda)$-пропагатора графа $\Gamma$ относительно $v$
для какого-либо расширения $F'$ поля $F_0(\lambda)$, где $F_0$ --- простое подполе поля $F$.
Мы зафиксируем эти факты (позволяющие для построения пропагаторов изменять
удобным образом поля и снабжать их подходящими абсолютными значениями)
в несколько уточненном виде
в предложении~\ref{p4.1} (являющемся аналогом для пропагаторов предложения~\ref{p3.1} для собственных функций).

Если $\Gamma$ --- локально конечный граф, $F$ --- поле  с абсолютным значением $|.|_{\rm v}$, $\lambda \in F$, $v \in V(\Gamma)$ и $f$ --- $(F,\lambda)$-пропагатор
графа $\Gamma$ относительно вершины $v$, то {\it расширенным носителем} $(F,\lambda)$-пропагатора $f$ будем называть объединение
носителя $f$ и $\{v\}$ (вообще говоря, $v$ может не принадлежать носителю $f$).

\begin{propos}
\label{p4.1}

Пусть $\Gamma$ --- локально конечный граф, $F$ --- поле
с абсолютным значением $|.|_{\rm v}$, $F_0$ --- простое подполе поля $F$ и  $\lambda \in F$.
Тогда для $v \in V(\Gamma)$ справедливы следующие утверждения$:$

$1)$ Если $\varphi$ --- вложение поля $F$ в некоторое поле $F'$ с произвольным абсолютным
значением $|.|_{{\rm v}'}$
и $f$ --- $(F',\varphi(\lambda))$-пропагатор графа $\Gamma$ относительно вершины $v$ такой,
что $f(u) \in \varphi(F)$ для всех $u \in V(\Gamma)$, то функция $\varphi^{-1} f \in F^{V(\Gamma)}$
является $(F,\lambda)$-пропагатором графа $\Gamma$ относительно вершины $v$.

$2)$ Если $f$ --- $(F,\lambda)$-пропагатор графа $\Gamma$ относительно вершины $v$, то найдется $(F_0(\lambda),\lambda)$-пропагатор графа $\Gamma$
относительно вершины $v$, расширенный носитель которого содержится в расширенном носителе $f$.

\end{propos}

\begin{proof}

Утверждение $1)$ очевидно. Доказательство утверждения $2)$
аналогично доказательству утверждения $2)$ предложения~\ref{p3.1}, но система линейных уравнений над $F$ с неизвестными $x_u$, где $u \in V(\Gamma)$,
теперь состоит, во-первых,  из уравнений
$(\sum_{u' \in \Gamma(u)}x_{u'}) - \lambda x_u = 0$ для всех $u \in V(\Gamma) \setminus \{v\}$,
во-вторых, из уравнений $x_u = 0$ для всех не принадлежащих расширенному носителю $f$ вершин $u$ графа $\Gamma$,
в-третьих, из уравнения  $(\sum_{u' \in \Gamma(v)}x_{u'}) - \lambda x_v = 1$.

\end{proof}

Пусть $\Gamma$ --- локально конечный граф, $F$ --- поле  с абсолютным значением $|.|_{\rm v}$
и $\lambda \in F$.
Если у графа $\Gamma$ относительно каждой его вершины $u$ имеется некоторый $(F,\lambda)$-пропагатор, скажем $p_u$,
то матрица ${\bf P}$ из $M_{V(\Gamma)\times V(\Gamma)}(F)$, у которой
для $u_1, u_2 \in V(\Gamma)$ элемент на пересечении строки, соответствующей $u_1$, и столбца, соответствующего $u_2$,
равен $p_{u_2}(u_1)$, является правой обратной для матрицы ${\bf A}_{\Gamma,F} - \lambda {\bf E}$ (а с учетом симметричности
${\bf A}_{\Gamma,F} - \lambda {\bf E}$ матрица ${\bf P}^t$ является левой обратной для ${\bf A}_{\Gamma,F} - \lambda {\bf E}$).
Обратно, если ${\bf P}$ --- матрица из $M_{V(\Gamma)\times V(\Gamma)}(F)$, которая является правой обратной
для ${\bf A}_{\Gamma,F} - \lambda {\bf E}$ (последнее эквивалентно тому, что ${\bf P}^t$
является левой обратной для ${\bf A}_{\Gamma,F} - \lambda {\bf E}$),
то для каждого $u \in V(\Gamma)$ функция $p_u \in F^{V(\Gamma)}$ такая,
что ее значение $p_u(w)$ для произвольного $w \in V(\Gamma)$ равно элементу матрицы ${\bf P}$ на пересечении строки, соответствующей $w$, и столбца,
соответствующего $u$,
является $(F,\lambda)$-пропагатором графа $\Gamma$ относительно $u$.
Наконец, заметим, что для $\lambda \in F$ с учетом следствия~\ref{c2.1} {\it сюръективность оператора $A_{\Gamma,F} - \lambda E$ эквивалентна наличию
у графа $\Gamma$ относительно каждой его вершины $(F,\lambda)$-пропагатора}.

Для локально конечного графа $\Gamma$, поля $F$  с абсолютным значением $|.|_{\rm v}$
и $\lambda \in F$ вершину $v \in V(\Gamma)$ назовем
{\it $(F,\lambda)$-сингулярной}, если у графа $\Gamma$ отсутствуют $(F,\lambda)$-пропагаторы относительно $v$.
Через $S_{F,\lambda}(\Gamma)$ будем обозначать
множество всех $(F,\lambda)$-сингулярных вершин графа $\Gamma$.
(Таким образом, с учетом замеченного выше, условие
сюръективности $A_{\Gamma,F} - \lambda E$, также как и условие наличия в
$M_{V(\Gamma)\times V(\Gamma)}(F)$ правой обратной
для ${\bf A}_{\Gamma,F} - \lambda {\bf E}$ матрицы, эквивалентно условию $S_{F,\lambda}(\Gamma) = \emptyset$.)
Ясно, что множество $S_{F,\lambda}(\Gamma)$ является $Aut(\Gamma)$-инвариантным, а из предложения~\ref{p4.1}
следует, что $S_{F,\lambda}(\Gamma) = S_{F',\lambda}(\Gamma)$ для произвольного расширения $F'$
(с любым абсолютным значением) поля $F_0(\lambda)$, где $F_0$ --- простое подполе поля $F$.

В случае конечного графа $\Gamma$ для элемента $\lambda$ поля $F$
условие $S_{F,\lambda}(\Gamma) \not = \emptyset$, как легко убедиться, равносильно условию
принадлежности $\lambda$ спектру графа $\Gamma$ над $F$, причем для $v \in V(\Gamma)$ условие
$v \in S_{F,\lambda}(\Gamma)$ равносильно условию принадлежности $v$ носителю некоторой
собственной функции матрицы смежности графа $\Gamma$ над $F$, соответствующей собственному значению
$\lambda$.
Ниже (см. теорему~\ref{t4.1}) будет доказано, что и в общем случае произвольного локально конечного
графа $\Gamma$ множество $S_{F,\lambda}(\Gamma)$ совпадает с объединением всех
конечных носителей собственных
функций оператора $A_{\Gamma,F}$, соответствующих собственному значению $\lambda$.

Из теоремы Теплица (см. предложение~\ref{p2.2}) очевидным образом вытекает справедливость следующиго утверждения.

\begin{propos}
\label{p4.2}

Пусть $\Gamma$ --- локально конечный граф, $F$ --- поле
с абсолютным значением $|.|_{\rm v}$ и  $\lambda \in F$.
Тогда для $v \in V(\Gamma)$ следующие условия
$(i)$ и $(ii)$ равносильны$:$

$(i)$ $v \in S_{F,\lambda}(\Gamma)$.

$(ii)$ Для некоторого конечного подмножества $X$ множества $V(\Gamma)$ не существует
$f' \in F^{V(\Gamma)}$ такой, что
$((A_{\Gamma,F} - \lambda E)(f'))(v) = 1$ и $((A_{\Gamma,F} - \lambda E)(f'))(u) = 0$ для всех $u \in X \setminus \{v\}$. $($Ясно, что при этом тем же свойством
обладает и любое содержащее $X$ конечное подмножество
множества $V(\Gamma)$.$)$

\end{propos}

Следствием предложений~\ref{p4.2} и~\ref{p4.1} является

\begin{propos}
\label{p4.3}

Пусть $\Gamma$ --- локально конечный граф, $F$ --- поле
с абсолютным значением $|.|_{\rm v}$ и  $\lambda \in F$.
Пусть, кроме того, $v \in S_{F,\lambda}(\Gamma)$.
Для каждого $r \in \mathbb{Z}_{\geq 0}$ определим систему $PrSyst_r$ линейных уравнений над $F$
относительно неизвестных $x_u$, где $u \in B_{\Gamma}(v,r + 1)$, включив в нее уравнение
$(\sum_{u' \in \Gamma(v)}x_{u'}) - \lambda x_v = 1$
и для каждого $u \in B_{\Gamma}(v,r) \setminus \{v\}$ уравнение $(\sum_{u' \in \Gamma(u)}x_{u'}) - \lambda x_u = 0$.
Тогда найдется такое $\tilde{r}_{\lambda,v} \in \mathbb{Z}_{\geq 0}$
$($$``$радиус пропагаторной несовместности для $\lambda$ и $v$$"$$)$,
что система $PrSyst_r$ совместна $($над $F$$)$ тогда и только тогда, когда $r < \tilde{r}_{\lambda,v}$. При этом
число $\tilde{r}_{\lambda,v}$ не изменится, если заменить $F$ на $F_0(\lambda)$, где $F_0$ --- простое подполе
поля $F$.

\end{propos}

Для локально конечного связного графа $\Gamma$ и элемента $\lambda$ поля $F$
вполне может оказаться, что $S_{F,\lambda}(\Gamma) \not = \emptyset$ и даже
$S_{F,\lambda}(\Gamma) = V(\Gamma)$
(таковыми являются, например, конечные регулярные степени $d \in \mathbb{Z}_{\geq 1}$ связные графы
для $\lambda = d\cdot 1_{F}$ или конечные связные вершинно-симметрические (т. е. с вершинно-транзитивными
группами автоморфизмов)
графы для $\lambda$, принадлежащих их спектрам над  $F$;
примеры бесконечных локально конечных связных вершинно-симметрических графов $\Gamma$
с $S_{\mathbb{C},\lambda}(\Gamma) = V(\Gamma)$ для некоторых $\lambda$
см. в \S~\ref{s8}, раздел~\ref{s8.5}).
Однако, как будет показано далее (см., например, следствие~\ref{c4.2}), для
локально конечного графа $\Gamma$ и элемента $\lambda$ поля $F$ в определенном смысле
типичной все же является
ситуация, когда у $\Gamma$ относительно каждой его вершины имеется $(F,\lambda)$-пропагатор,
т. е. когда $S_{F,\lambda}(\Gamma) = \emptyset$.

\medskip

Следующие предложения~\ref{p4.4} и~\ref{p4.5}
являются аналогами для пропагаторов предложений~\ref{p3.3} и~\ref{p3.4} для собственных функций операторов смежности
локально конечных графов.

\begin{propos}
\label{p4.4}

Пусть $\Gamma$ --- локально конечный граф, $F$ --- поле с
абсолютным значением $|.|_{\rm v}$, $\lambda \in F$ и $v \in V(\Gamma)$.
Если $X_i$, $i \in I$, --- непустое линейно упорядоченное по включению семейство подмножеств множества
$V(\Gamma)$, являющихся расширенными носителями $(F,\lambda)$-пропагаторов графа $\Gamma$
относительно $v$,
то некоторое  подмножество множества $\cap_{i \in I}X_i$ является
расширенным носителем $(F,\lambda)$-пропагатора графа $\Gamma$
относительно $v$.
В частности, каждый расширенный носитель $(F,\lambda)$-пропагатора графа $\Gamma$
относительно $v$ содержит
минимальный по включению среди расширенных носителей $(F,\lambda)$-про\-пагаторов графа $\Gamma$
относительно $v$.

\end{propos}

\begin{proof}

Утверждение следует из теоремы Теплица (см. предложение~\ref{p2.2}),
примененной к системе линейных уравнений над $F$ с неизвестными $x_u$, где $u \in V(\Gamma)$,
состоящей, во-первых, из уравнений
$(\sum_{u' \in \Gamma(u)}x_{u'}) - \lambda x_u = 0$ для всех $u \in V(\Gamma) \setminus \{v\}$,
во-вторых, из уравнений $x_u = 0$ для всех $u \in V(\Gamma) \setminus \cap_{i \in I}X_i$,
в-третьих, из уравнения $(\sum_{v' \in \Gamma(v)}x_{v'}) - \lambda x_v = 1$.

\end{proof}

Если $\Gamma$ --- локально конечный граф, $F$ --- поле  с абсолютным значением $|.|_{\rm v}$, $\lambda \in F$,
$v \in V(\Gamma)$ и $p_v$ ---
$(F,\lambda)$-пропагатор графа $\Gamma$ относительно вершины $v$, то для $X \subseteq V(\Gamma)$,
содержащего $v$ и являющегося
объединением набора связных компонент расширенного носителя $p_v$, функция из $F^{V(\Gamma)}$,
совпадающая с $p_v$ на $X$
и тождественно равная нулю на $V(\Gamma) \setminus X$, не является, вообще говоря,
$(F,\lambda)$-пропагатором графа $\Gamma$ относительно
вершины $v$, хотя ограничение $p_v$ на $X$ очевидным
образом является $(F,\lambda)$-пропагатором графа $\langle X \rangle_{\Gamma}$ относительно
вершины $v$.
Однако (аналогично тому, как это имеет место для собственных функций, см. \S~\ref{s3}) такого рода эффект не возникает,
если вместо связных компонент расширенного носителя $p_v$ рассматривать его $\Gamma^2$-связные компоненты.

Действительно, если $\Gamma$ --- локально конечный граф, $F$ --- поле  с абсолютным значением $|.|_{\rm v}$, $\lambda \in F$, $v \in V(\Gamma)$ и
$p_v$ --- $(F,\lambda)$-пропагатор графа $\Gamma$ относительно
вершины $v$, то для $X \subseteq V(\Gamma)$, содержащего $v$ и являющегося
объединением некоторого набора ${\Gamma^2}$-связных компонент расширенного носителя $p_v$, функция из $F^{V(\Gamma)}$, совпадающая с $p_v$ на $X$
и тождественно равная нулю на $V(\Gamma) \setminus X$, очевидным образом также является $(F,\lambda)$-пропагатором графа $\Gamma$ относительно
вершины $v$. В частности, содержащая $v$ компонента  ${\Gamma^2}$-связности расширенного носителя $p_v$, которую уместно назвать
{\it основной компонентой ${\Gamma^2}$-связности} расширенного носителя $p_v$,
является расширенным носителем $(F,\lambda)$-пропагатора графа $\Gamma$ относительно вершины $v$.
Как следствие,
для обсуждавшихся ранее (см. предложение~\ref{p4.4}) минимальных по включению
среди расширенных носителей $(F,\lambda)$-про\-пагаторов графа $\Gamma$
относительно $v$ справедливо следующее утверждение.

\begin{propos}
\label{p4.5}

Пусть $\Gamma$ --- локально конечный граф, $F$ --- поле с
абсолютным значением $|.|_{\rm v}$, $\lambda \in F$ и $v \in V(\Gamma)$.
Тогда каждый минимальный по включению среди расширенных носителей $(F,\lambda)$-про\-пагаторов графа $\Gamma$
относительно $v$ является $\Gamma^2$-связным.

\end{propos}

\begin{remark}
\label{r4.2}

В силу утверждения 2) предложения~\ref{p4.1} для каждого минимального подмножества из предыдущего предложения
$(F,\lambda)$-пропагатор графа $\Gamma$ относительно вершины $v$, для которого это подмножество
есть расширенный носитель, может быть выбран из $F_0(\lambda)^{V(\Gamma)}$, т. е. может быть выбран являющимся
$(F_0(\lambda),\lambda)$-пропагатором.

\end{remark}

В дальнейшем нам потребуется следующее

\begin{propos}
\label{p4.6}

Пусть $\Gamma$ --- локально конечный граф, $F$ --- поле
с абсолютным значением $|.|_{\rm v}$ и $\lambda \in F$.
Пусть, кроме того, для вершины $v$ графа $\Gamma$ имеется такая последовательность
$(v_i)_{i \in \mathbb{Z}_{\geq 1}}$
попарно различных вершин графа $\Gamma$, что для каждого $i \in \mathbb{Z}_{\geq 1}$
существует $(F,\lambda)$-пропагатор графа $\Gamma$ относительно вершины $v_i$ с конечным $\Gamma^2$-связным расширенным носителем $X_i$,
содержащим $v$. Тогда найдется собственная функция оператора $A_{\Gamma,F}$, соответствующая собственному значению $\lambda$,
с бесконечным носителем, содержащим вершину $v$ и содержащимся в $\cup_{i \in \mathbb{Z}_{\geq 1}} X_i$.

\end{propos}

\begin{proof}

По условию для каждого $i \in \mathbb{Z}_{\geq 1}$
существует $(F,\lambda)$-про\-пагатор $p_i$ графа $\Gamma$ относительно вершины $v_i$ с расширенным носителем $X_i$.
Положим $X:= \cup_{i \in \mathbb{Z}_{\geq 1}} X_i$.

В силу локальной конечности графа $\Gamma$ и $\Gamma^2$-связности каждого из множеств
$X_i, i \in I,$ существует такой бесконечный путь $u_0 = v, u_1, u_2,...$ графа $\Gamma^2$ с попарно различными вершинами,
что для любого целого положительного числа $n$ найдется $i \in I$ со свойством $\{u_0,...,u_n\} \subseteq X_{i}$.

Пусть $k \in \mathbb{Z}_{\geq 0}$.
Зададим следующим образом систему линейных уравнений над полем $F$ с неизвестными $x_u$, где $u \in V(\Gamma)$:
система состоит, во-первых, из уравнений
$(\sum_{u' \in \Gamma(u)}x_{u'}) - \lambda x_u = 0$ для всех $u \in V(\Gamma)$,
во-вторых, из уравнений $x_u = 0$ для всех $u \in V(\Gamma) \setminus X$,
в-третьих, из уравнения  $x_{u_k} = 1$.
В силу наличия функций $p_{1},p_{2},p_{3},...$ каждая конечная подсистема этой системы совместна над $F$.
Следовательно, по теореме Теплица (см. предложение~\ref{p2.2}) совместна над $F$ и вся система.
Пусть $x'_u \in F$, где $u$ пробегает множество $V(\Gamma)$, --- некоторое
ее решение над $F$. Тогда, определяя $f_k \in F^{V(\Gamma)}$ посредством $f_k(u) = x'_u$ для всех  $u \in V(\Gamma)$, мы получает
собственную функцию оператора $A_{\Gamma,F}$, соответствующую собственному значению $\lambda$, с носителем,
содержащимся в $X$, и такую, что $f_k(u_k) = 1$.

Если носитель функции $f_0$ бесконечен, то она обладает требуемыми в предложении~\ref{p4.6}
свойствами. Предположим, что носитель функции $f_0$ конечен.
Если для некоторого $k \in \mathbb{Z}_{\geq 1}$ носитель функции $f_k$ бесконечен, то функция
$f_k$ или функция $f_k + f_0$ обладает требуемыми в предложении~\ref{p4.6}
свойствами. Предположим поэтому, что для каждого $k \in \mathbb{Z}_{\geq 0}$
носитель функции $f_k$ конечен. Но в этом случае справедливость предложения~\ref{p4.6}
вытекает из предложения~\ref{p3.5} (в котором следует положить $I = \mathbb{Z}_{\geq 0}$,
а в качестве $X_i$, $i \in I$, взять носитель функции $f_i$).

\end{proof}

Для локально конечного графа $\Gamma$, поля $F$  с абсолютным значением $|.|_{\rm v}$
и $\lambda \in F$ вершину $v \in V(\Gamma)$ назовем
{\it $(F,\lambda)$-локальной}, если у графа $\Gamma$ имеется $(F,\lambda)$-пропагатор относительно $v$
с конечным носителем или, другими словами, $\delta_v \in (A_{\Gamma,F} - \lambda E)(F^{V(\Gamma)*})$.
Через $L_{F,\lambda}(\Gamma)$ будем обозначать
множество всех $(F,\lambda)$-локальных вершин графа $\Gamma$.
(Таким образом, условие, что $A_{\Gamma,F} - \lambda E$
отображает $F^{V(\Gamma)*}$ сюръективно на себя, также как и условие наличия в
$M^{rcf}_{V(\Gamma)\times V(\Gamma)}(F)$ обратной
для ${\bf A}_{\Gamma,F} - \lambda {\bf E}$ матрицы, эквивалентно условию $L_{F,\lambda}(\Gamma) = V(\Gamma)$.)
Ясно, что множество $L_{F,\lambda}(\Gamma)$ является $Aut(\Gamma)$-инвариантным, а из утверждения 2) предложения~\ref{p4.1}
следует, что $L_{F,\lambda}(\Gamma) = L_{F',\lambda}(\Gamma)$ для произвольного расширения $F'$
(с произвольным абсолютным значением) поля $F_0(\lambda)$, где $F_0$ --- простое подполе поля $F$.
Ниже (см. теорему~\ref{t4.2}) будет доказано, что $L_{F,\lambda}(\Gamma)$ есть в точности дополнение в $V(\Gamma)$ объединения всех
носителей собственных
функций оператора $A_{\Gamma,F}$, соответствующих собственному значению $\lambda$. (Таким образом, окажется,
что $L_{F,\lambda}(\Gamma)$ есть в точности множество тех вершин $v$ графа $\Gamma$, для которых
в соответствии с предложением~\ref{p3.2} определено число $r_{\lambda,v}$.)

В случае конечного графа $\Gamma$ очевидным образом справедливо разбиение
$V(\Gamma) = S_{F,\lambda}(\Gamma) \cup L_{F,\lambda}(\Gamma)$,
причем $L_{F,\lambda}(\Gamma) \not = V(\Gamma)$
тогда и только тогда, когда $\lambda$ принадлежит спектру графа $\Gamma$ над $F$.

\medskip

Для дальнейшего полезно сделать несколько по существу тривиальных наблюдений относительно собственных
функций матрично записанного оператора смежности локально конечного графа. Пусть $\Gamma$ --- локально конечный граф,
$F$ --- поле  с абсолютным значением $|.|_{\rm v}$ и $\lambda \in F$.
Для $u \in V(\Gamma)$ (в соответствии с \S~\ref{s2}) полагаем $\alpha_u := (\sum_{u' \in \Gamma(u)}\delta_{u'}) - \lambda \delta_u$.
Пусть, кроме того, $\boldsymbol{\alpha}_u$, где $u \in V(\Gamma)$, --- столбец матрицы
${\bf A}_{\Gamma,F} - \lambda {\bf E}$, соответствующий $u$.
(Таким образом, $\boldsymbol{\alpha}_u^t$ --- строка матрицы ${\bf A}_{\Gamma,F} - \lambda {\bf E}$,
соответствующая $u$.) Ясно, что элемент столбца $\boldsymbol{\alpha}_u$, стоящий в строке, соответствующей
$v \in V(\Gamma)$, равен $\alpha_u(v)$.
Далее, в силу ${\bf A}_{\Gamma,F} - \lambda {\bf E} \in M^{rcf}_{V(\Gamma)\times V(\Gamma)}(F)$
для произвольной функции $f \in F^{V(\Gamma)}$ естественным образом корректно определен столбец $\sum_{u \in V(\Gamma)} f(u)\boldsymbol{\alpha}_u$.
При этом $f$ тогда и только тогда является собственной функцией оператора $A_{\Gamma,F}$, соответствующей собственному значению $\lambda$,
когда $f$ --- ненулевая функция со свойством $\sum_{u \in V(\Gamma)} f(u)\boldsymbol{\alpha}_u = {\bf 0},$ где
${\bf 0}$ --- нулевой столбец.
(Аналогично, $f$ тогда и только тогда является $(F,\lambda)$-пропагатором графа $\Gamma$ относительно его вершины $v$,
когда $\sum_{u \in V(\Gamma)} f(u)\boldsymbol{\alpha}_u = \boldsymbol{\delta}_v$, где $ \boldsymbol{\delta}_v$ --- столбец, у которого на месте,
соответствующем вершине $v$, стоит 1, а на остальных местах стоят нули.)
В частности, если $X$ --- конечное подмножество множества $V(\Gamma)$ и $v$ --- некоторая вершина из $X$, то наличие у оператора
$A_{\Gamma,F}$ собственной функции, соответствующей собственному значению $\lambda$, носитель которой содержится в $X$ и содержит $v$,
эквивалентно тому, что столбец $\boldsymbol{\alpha}_v$ равен линейной комбинации (с коэффициентами из $F$)
столбцов $\boldsymbol{\alpha}_u$, $u \in X \setminus \{v\}$, или, равносильно, что функция $\alpha_v$
равна линейной комбинации (с коэффициентами из $F$)
функций $\alpha_u$, $u \in X \setminus \{v\}$.

\begin{propos}
\label{p4.7}

Пусть $\Gamma$ --- локально конечный граф, $F$ --- поле  с абсолютным значением $|.|_{\rm v}$ и $\lambda \in F$.
Тогда для конечного подмножества $X$ множества $V(\Gamma)$ и $v \in X$ следующие условия
$(i)$ и $(ii)$ равносильны$:$

$(i)$ У оператора $A_{\Gamma,F}$ имеется такая собственная функция $f$, соответствующая собственному значению $\lambda$, что носитель $f$
содержится в $X$ и содержит $v$.

$(ii)$ Не существует $f' \in F^{V(\Gamma)}$ такой, что
$((A_{\Gamma,F} - \lambda E)(f'))(v) = 1$ и $((A_{\Gamma,F} - \lambda E)(f'))(u) = 0$ для всех $u \in X \setminus \{v\}$.

\end{propos}

\begin{proof}
Для произвольных $u_1,u_2 \in V(\Gamma)$ пусть
$\alpha_{u_1,u_2} = \alpha_{u_2}(u_1)$ --- элемент
матрицы ${\bf A}_{\Gamma,F} - \lambda {\bf E}$, стоящий на пересечении строки, соответствующей $u_1$,
и столбца, соответствующего $u_2$. (При этом $\alpha_{u_1,u_2} = \alpha_{u_2,u_1}$ в силу
симметричности матрицы ${\bf A}_{\Gamma,F} - \lambda {\bf E}$.)

Предположим, что имеет место $(i)$, но существует функция $f' \in F^{V(\Gamma)}$ такая, что
$((A_{\Gamma,F} - \lambda E)(f'))(v) = 1$ и $((A_{\Gamma,F} - \lambda E)(f'))(u) = 0$ для всех $u \in X \setminus \{v\}$.
Не теряя общности, будем считать, что носитель $f'$ содержится в $X \cup \Gamma(X)$
(и следовательно, носители функций $f$ и $f'$ конечны).
Тогда, с одной стороны,
\begin{equation}
\label{eq4.1}
\sum_{u_1,u_2 \in V(\Gamma)}f'(u_1)\alpha_{u_1,u_2}f(u_2) = \sum_{u_1 \in V(\Gamma)}f'(u_1)
\big(\sum_{u_2 \in V(\Gamma)}\alpha_{u_1,u_2}f(u_2)\big) = 0.
\end{equation}
Но, с другой стороны,
$$\sum_{u_1,u_2 \in V(\Gamma)}f'(u_1)\alpha_{u_1,u_2}f(u_2) = \sum_{u_2 \in V(\Gamma)}f(u_2)
\big(\sum_{u_1 \in V(\Gamma)}\alpha_{u_1,u_2}f'(u_1)\big)$$
$$= \sum_{u_2 \in X}f(u_2)
\big(\sum_{u_1 \in X \cup \Gamma(X)}\alpha_{u_1,u_2}f'(u_1)\big) = \sum_{u_2 \in X}f(u_2)
\big(\sum_{u_1 \in X \cup \Gamma(X)}\alpha_{u_2,u_1}f'(u_1)\big)$$
$$= f(v) \not = 0.$$
Таким образом, $(i)$ влечет $(ii)$.

Докажем, что $(ii)$ влечет $(i)$. Если функция $\alpha_v \in F^{V(\Gamma)}$
принадлежит линейной оболочке в $F^{V(\Gamma)}$ функций $\alpha_u$, $u \in X\setminus \{v\}$, то,
как было отмечено непосредственно перед формулировкой предложения, имеет место $(i)$.
Предположим поэтому, что имеет место $(ii)$, но функция $\alpha_v$ не принадлежит линейной оболочке
в $F^{V(\Gamma)}$
функций $\alpha_u$, $u \in X \setminus \{v\}$. Тогда на содержащем $\alpha_u, u \in X,$ подпространстве
$\langle \delta_u : u \in X \cup \Gamma (X)\rangle$ векторного пространства $F^{V(\Gamma)}$ найдется
такой линейный функционал $\chi$, что $\chi(\alpha_v) = 1$ и $\chi(\alpha_u) = 0$ для всех
$u \in X \setminus \{v\}$.
Определим $f'' \in F^{V(\Gamma)}$, полагая $f''(u) = \chi(\delta_u)$ для всех $u \in X \cup \Gamma (X)$
и полагая $f''(u) = 0$ для всех  $u \in V(\Gamma) \setminus (X \cup \Gamma (X))$.
Тогда, как легко видеть, для функции $f''$ имеем
$((A_{\Gamma,F} - \lambda E)(f''))(v) = 1$ и $((A_{\Gamma,F} - \lambda E)(f''))(u) = 0$ для всех $u \in X \setminus \{v\}$.
Полученное противоречие со сделанным предположением о том, что имеет место $(ii)$, завершает доказательство
предложения.

\end{proof}

Следующее предложение в определенном смысле дуально предложению~\ref{p4.7}.

\begin{propos}
\label{p4.8}

Пусть $\Gamma$ --- локально конечный граф, $F$ --- поле  с абсолютным значением $|.|_{\rm v}$ и $\lambda \in F$.
Тогда для конечного подмножества $X$ множества $V(\Gamma)$ и $v \in X$ следующие условия
$(i)$ и $(ii)$ равносильны$:$

$(i)$ Не существует $f \in F^{V(\Gamma)}$ такой, что $v$ принадлежит носителю $f$ и
$((A_{\Gamma,F} - \lambda E)(f))(u) = 0$ для всех $u \in X$.

$(ii)$ У графа $\Gamma$ относительно вершины $v$ имеется такой $(F,\lambda)$-пропа\-гатор $f'$, что носитель
$f'$ содержится в $X$.

\end{propos}

\begin{proof}
Как и при доказательстве предложения~\ref{p4.7}, для произвольных $u_1,u_2 \in V(\Gamma)$ пусть
$\alpha_{u_1,u_2} = \alpha_{u_2}(u_1)$ --- элемент
матрицы ${\bf A}_{\Gamma,F} - \lambda {\bf E}$, стоящий на пересечении строки, соответствующей $u_1$,
и столбца, соответствующего $u_2$.

Предположим, что имеет место $(ii)$, но существует функция $f \in F^{V(\Gamma)}$ такая, что $f(v) \not = 0$ и
$((A_{\Gamma,F} - \lambda E)(f))(u) = 0$ для всех $u \in X$.
При этом, не теряя общности, будем считать, что носитель $f$ содержится в $X \cup \Gamma(X)$
(и следовательно, носители функций $f$ и $f'$ конечны).
Тогда, как легко убедиться, для $f$ и $f'$, с одной стороны, сохраняет силу~\eqref{eq4.1}, а с другой стороны,
$$\sum_{u_1,u_2 \in V(\Gamma)}f'(u_1)\alpha_{u_1,u_2}f(u_2) =
\sum_{u_1 \in X,u_2 \in X \cup \Gamma(X)}f'(u_1)\alpha_{u_1,u_2}f(u_2)$$
$$= \sum_{u_2 \in X \cup \Gamma(X)}f(u_2)\big(\sum_{u_1 \in X}\alpha_{u_1,u_2}f'(u_1)\big)
= \sum_{u_2 \in X \cup \Gamma(X)}f(u_2)\big(\sum_{u_1 \in V(\Gamma)}\alpha_{u_1,u_2}f'(u_1)\big)$$
$$= \sum_{u_2 \in X \cup \Gamma(X)}f(u_2)\big(\sum_{u_1 \in V(\Gamma)}\alpha_{u_2,u_1}f'(u_1)\big)= f(v) \not = 0.$$
Полученное противоречие означает, что $(ii)$ влечет $(i)$.

Докажем, что $(i)$ влечет $(ii)$. Условие $(ii)$
эквивалентно условию, что $\delta_v$ принадлежит
линейной оболочке в
$F^{V(\Gamma)}$ функций $\alpha_u$, $u \in X$. Предположим, что
имеет место $(i)$, но функция $\delta_v$ не принадлежит линейной оболочке
в $F^{V(\Gamma)}$
функций $\alpha_u$, $u \in X$.
Тогда на содержащем $\alpha_u, u \in X,$ подпространстве
$\langle \delta_u : u \in X \cup \Gamma (X)\rangle$ векторного пространства $F^{V(\Gamma)}$ найдется
такой линейный функционал $\chi$, что $\chi(\delta_v) = 1$ и $\chi(\alpha_u) = 0$ для всех
$u \in X$.
Определим $f'' \in F^{V(\Gamma)}$, полагая $f''(u) = \chi(\delta_u)$ для всех $u \in X \cup \Gamma (X)$
и полагая $f''(u) = 0$ для всех  $u \in V(\Gamma) \setminus (X \cup \Gamma (X))$.
Тогда, как легко видеть, для функции $f''$ имеем
$f''(v) = 1$ и $((A_{\Gamma,F} - \lambda E)(f''))(u) = 0$ для всех $u \in X$.
Полученное противоречие со сделанным предположением о том, что имеет место $(i)$, завершает доказательство
предложения.

\end{proof}

Из предложения~\ref{p4.7} и теоремы Теплица (см. предложение~\ref{p2.2}) вытекает справедливость следующего утверждения.

\begin{theorem}
\label{t4.1}

Пусть $\Gamma$ --- локально конечный граф, $F$ --- поле  с абсолютным значением $|.|_{\rm v}$ и $\lambda \in F$.
Тогда для $v \in V(\Gamma)$ следующие условия
$(i)$ и $(ii)$ эквивалентны$:$

$(i)$ У оператора $A_{\Gamma,F}$ имеется собственная функция, соответствующая собственному значению $\lambda$,  носитель которой конечен и содержит $v$.

$(ii)$ У графа $\Gamma$ отсутствуют $(F,\lambda)$-пропагаторы относительно вершины $v$.
Другими словами, $v \in S_{F,\lambda}(\Gamma)$.

\end{theorem}

\begin{cor}
\label{c4.1}

Пусть $\Gamma$ --- локально конечный граф, $F$ --- поле  с абсолютным значением $|.|_{\rm v}$ и $\lambda \in F$.
Тогда множество $S_{F,\lambda}(\Gamma)$
совпадает с объединением всех конечных носителей собственных функций оператора $A_{\Gamma,F}$, соответствующих собственному значению $\lambda$.

\end{cor}

\begin{cor}
\label{c4.2}

Пусть $\Gamma$ --- локально конечный граф, $F$ --- поле с абсолютным значением $|.|_{\rm v}$, $F_0$ --- простое подполе поля $F$  и $\lambda \in F$.
Пусть, кроме того,
$v \in S_{F,\lambda}(\Gamma)$.
Тогда справедливы следующие утверждения$:$

$1)$ Найдется такое конечное содержащее $v$ подмножество $X$ множества $V(\Gamma)$, что порожденный $X$ подграф
$\langle X \rangle_{\Gamma}$ графа $\Gamma$ связен, а $\lambda$ является собственным значением матрицы смежности
графа $\langle X \rangle_{\Gamma}$ над $F$, причем имеется соответствующая $\lambda$ собственная функция
матрицы смежности графа $\langle X \rangle_{\Gamma}$ над $F$, содержащаяся в
$(F_0(\lambda))^X$, носитель которой равен $X$. Кроме того, найдется такое конечное содержащее $X$ подмножество $Y$
множества $V(\Gamma)$, что порожденный $Y$ подграф
$\langle Y \rangle_{\Gamma^2}$ графа $\Gamma^2$ связен, а $\lambda$ является собственным значением матрицы смежности
графа $\langle Z \rangle_{\Gamma}$ над $F$ для любого содержащего $Y$ конечного подмножества $Z$
множества $V(\Gamma)$, причем имеется соответствующая $\lambda$ собственная функция
матрицы смежности графа $\langle Z \rangle_{\Gamma}$ над $F$, содержащаяся в
$(F_0(\lambda))^Z$, носитель которой равен $Y$. В частности, $\lambda$ --- алгебраический элемент над $F_0$ $($и поле $F_0(\lambda)$ конечно в случае поля $F$
положительной характеристики$)$.

$2)$ Если $F$ - поле характеристики $0$ $($т. е. $F_0 = \mathbb{Q}$$)$, а $X$ --- такое, как в $1)$,
то $\lambda$ --- корень унитарного многочлена из $\mathbb{Z}[x]$ степени $|X|$, все корни которого вещественны
и не превосходят по $($обычной на $\mathbb{R}$$)$
абсолютной величине максимума степеней вершин графа $\langle X \rangle_{\Gamma}$
$($и не превосходят максимума степеней в графе $\Gamma$ вершин из $X$, причем достигает этого
максимума лишь в случае, когда $X$ --- содержащая $v$ связная компонента графа $\Gamma$
и $\langle X \rangle_{\Gamma}$ --- $($конечный$)$ регулярный граф$)$.

\end{cor}

\begin{proof}
Легко убедиться, что в качестве $X$ и $Y$ в $1)$ можно взять содержащую $v$ соответственно связную компоненту и $\Gamma^2$-связную компоненту
содержащего $v$ конечного носителя собственной функции
оператора $A_{\Gamma,F_0(\lambda)}$, соответствующей собственному значению $\lambda$, существование которого вытекает из теоремы~\ref{t4.1} (в которой
в качестве $F$ следует взять $F_0(\lambda)$). Утверждение 2) следствия вытекает из утверждения 1)
и хорошо известных свойств спектров конечных графов.

\end{proof}

\begin{remark}
\label{r4.3}

Пусть $\Gamma$ --- локально конечный граф, $F$ --- поле положительной характеристики с абсолютным значением $|.|_{\rm v}$, $F_0$ --- простое подполе
поля $F$ и $v \in V(\Gamma)$. Тогда для $\lambda \in F$
отсутствие $(F,\lambda)$-пропагатора графа $\Gamma$ относительно $v$ (влекущее согласно утверждению 1)
следствия~\ref{c4.2} алгебраичность $\lambda$ над $F_0$)
эквивалентно (см. утверждение 2) предложения~\ref{p4.1}) отсутствию $(F_0(\lambda),\lambda)$-пропагатора графа $\Gamma$ относительно $v$ (где
$F_0(\lambda)$ --- {\it конечное} поле), а также эквивалентно (см. теорему~\ref{t4.1})
принадлежности $v$ конечному носителю собственной функции оператора смежности
$A_{\Gamma,F_0(\lambda)}$ (где $F_0(\lambda)$ --- {\it конечное} поле), соответствующей собственному значению $\lambda$.

\end{remark}

Подобно тому, как предложение~\ref{p4.7} и теорема Теплица  (см. предложение~\ref{p2.2}) влекут теорему~\ref{t4.1},
предложение~\ref{p4.8} и теорема Теплица  (см. предложение~\ref{p2.2}) влекут следующую важную для дальнейшего теорему
(в определенном смысле дуальную теореме~\ref{t4.1}).

\begin{theorem}
\label{t4.2}

Пусть $\Gamma$ --- локально конечный граф, $F$ --- поле  с абсолютным значением $|.|_{\rm v}$ и $\lambda \in F$.
Тогда для $v \in V(\Gamma)$ следующие условия
$(i)$ и $(ii)$ равносильны$:$

$(i)$ Вершина $v$ не принадлежит носителю никакой собственной функции оператора $A_{\Gamma,F}$, соответствующей собственному значению $\lambda$.

$(ii)$ У графа $\Gamma$ относительно вершины $v$ имеется $(F,\lambda)$-пропагатор с конечным носителем.
Другими словами, $v \in L_{F,\lambda}(\Gamma)$.

\end{theorem}

\section{Инъективность и ряд других свойств операторов $A_{\Gamma,F} - \lambda E$ в терминах пропагаторов}
\label{s5}

В этом параграфе будет показано, что для локально конечного графа $\Gamma$, поля $F$ с абсолютным значением $|.|_{\rm v}$ и $\lambda \in F$
ряд важных свойств оператора $A_{\Gamma,F} - \lambda E$
может быть естественным образом сформулировано в терминах $(F,\lambda)$-пропагаторов графа $\Gamma$
относительно его вершин.  Что важно, в этом и следующих параграфах будет показана продуктивность
такой переформулировки.

Согласно теореме~\ref{t4.1} для локально конечного графа $\Gamma$, поля $F$ с некоторым абсолютным значением
$|.|_{\rm v}$ и $\lambda \in F$
условие $S_{F,\lambda}(\Gamma) \not = \emptyset $ равносильно условию существования у оператора $A_{\Gamma,F}$ собственной функции,
соответствующей собственному значению $\lambda$, с конечным носителем или, другими словами, условию, что ограничение оператора
$A_{\Gamma,F} - \lambda E$ на $F^{V(\Gamma)*}$ не является инъективным оператором.
Таким образом, с учетом сказанного в начале предыдущего параграфа (после доказательства предложения~\ref{p4.1}),
а также следствия~\ref{c2.1} получаем, что справедлива

\begin{theorem}
\label{t5.1}

Пусть $\Gamma$ --- локально конечный граф, $F$ --- поле  с абсолютным значением $|.|_{\rm v}$ и $\lambda \in F$.
Тогда следующие условия $(i)-(v)$ равносильны$:$

$(i)$ Оператор $A_{\Gamma,F} - \lambda E$ сюръективен.

$(ii)$ Матрица ${\bf A}_{\Gamma,F} - \lambda {\bf E}$ имеет в $M_{V(\Gamma)\times V(\Gamma)}(F)$ правую обратную матрицу или левую
обратную матрицу $($а следовательно, имеет также взаимно транспонированные правую обратную матрицу и левую обратную матрицу$)$.

$(iii)$ $S_{F,\lambda}(\Gamma) = \emptyset $.

$(iv)$ Оператор $A_{\Gamma,F} - \lambda E$ отображает подпространство $F^{V(\Gamma)*}$ пространства $F^{V(\Gamma)}$
инъективно в себя $($или, другими словами, у оператора $A_{\Gamma,F}$ отсутствуют собственные функции, соответствующие собственному значению
$\lambda$, с конечным носителем$)$.

$(v)$ $F^{V(\Gamma)*} \subseteq (A_{\Gamma,F} - \lambda E)(F^{V(\Gamma)})$.

\end{theorem}

Далее, очевидным следствием теоремы~\ref{t4.1} является

\begin{propos}
\label{p5.1}

Пусть $\Gamma$ --- локально конечный граф, $F$ --- поле  с абсолютным значением $|.|_{\rm v}$ и $\lambda \in F$.
Множество $S_{F,\lambda}(\Gamma)$ конечно тогда и только тогда, когда собственные функции оператора $A_{\Gamma,F}$,
соответствующие собственному значению $\lambda$ и имеющие конечный носитель, образуют конечномерное подпространство
пространства $F^{V(\Gamma)}$.

\end{propos}

Кроме того, из теоремы~\ref{t4.1} и предложения~\ref{p3.5} вытекает

\begin{propos}
\label{p5.2}

Пусть $\Gamma$ --- локально конечный граф, $F$ --- поле  с абсолютным значением $|.|_{\rm v}$ и $\lambda \in F$.
Если множество $S_{F,\lambda}(\Gamma)$ бесконечно, то у оператора $A_{\Gamma,F}$ имеется
собственная функция, соответствующая собственному значению $\lambda$, с бесконечным носителем,
содержащимся в $S_{F,\lambda}(\Gamma)$.

\end{propos}

\begin{remark}
\label{r5.1}

Из теоремы~\ref{t4.1} и предложения~\ref{p5.2} легко следует, что для локально конечного графа $\Gamma$,
поля $F$ с абсолютным значением $|.|_{\rm v}$ и $\lambda \in F$ бесконечность $S_{F,\lambda}(\Gamma)$ влечет принадлежность каждой вершины из $S_{F,\lambda}(\Gamma)$
бесконечному содержащемуся в $S_{F,\lambda}(\Gamma)$ носителю собственной функции оператора $A_{\Gamma,F}$, соответствующей собственному значению $\lambda$. (Действительно, согласно теореме 4.1 каждая вершина $v$ из $S_{F,\lambda}(\Gamma)$ принадлежит
содержащемуся в $S_{F,\lambda}(\Gamma)$ конечному носителю некоторой собственной функции $f_1$ оператора
$A_{\Gamma,F}$, соответствующей собственному
значению $\lambda$. Кроме того, согласно предложению 5.2 имеется собственная функция $f_2$ оператора
$A_{\Gamma,F}$, соответствующая собственному
значению $\lambda$, с бесконечным содержащимся в $S_{F,\lambda}(\Gamma)$ носителем.
Если $v$ принадлежит носителю $f_2$, то $f_2$ - требуемая функция. Если же $v$ не принадлежит
носителю $f_2$, то $f_1 + f_2$ --- требуемая функция.)

\end{remark}

\begin{remark}
\label{r5.2}

Отметим еще одно следствие теоремы~\ref{t4.1}: если $\Gamma$ --- локально конечный граф, $F$ --- поле  с абсолютным значением $|.|_{\rm v}$ и $\lambda \in F$,
то $S_{F,\lambda}(\Gamma)$ не имеет одноэлементных компонент связности при $\lambda \not = 0$ и не имеет одноэлементных компонент $\Gamma^2$-связности
при $\lambda = 0$ и отсутствии у $\Gamma$ одноэлементных компонент связности.

\end{remark}

Доказательство следующей теоремы во многом близко к доказательству теоремы 3 из~\cite{Wang}.

\begin{theorem}
\label{t5.2}

Пусть $\Gamma$ --- локально конечный граф, $F$ --- поле  с абсолютным значением $|.|_{\rm v}$ и $\lambda \in F$.
Если у графа $\Gamma$ относительно некоторой его вершины $v$ имеется $(F,\lambda)$-пропагатор с бесконечным носителем,
то у оператора $A_{\Gamma,F}$ имеется собственная функция $f$, соответствующая собственному значению $\lambda$,
с бесконечным носителем. Если, кроме того, у $\Gamma$ относительно $v$ отсутствуют $(F,\lambda)$-пропагаторы с конечным носителем
$($или, другими словами, $v \not \in L_{F,\lambda}(\Gamma)$$)$, то $f$ можно выбрать с дополнительным свойством, что носитель $f$ содержит $v$.

\end{theorem}

\begin{proof}
Заметим, что второе утверждение теоремы следует из ее первого утверждения и теоремы~\ref{t4.2}.
Действительно, в предположениях второго утверждения теоремы при справедливости ее первого утверждения
у оператора $A_{\Gamma,F}$ имеется собственная функция $f_1$, соответствующая собственному значению $\lambda$,
с бесконечным носителем, а из $v \not \in L_{F,\lambda}(\Gamma)$ согласно теореме~\ref{t4.2} следует
существование у оператора $A_{\Gamma,F}$ собственной функции $f_2$, соответствующей собственному значению $\lambda$,
носитель которой содержит $v$. Ясно, что среди функций $f_1$, $f_2$, $f_1 + f_2$
имеется функция с бесконечным носителем, содержащим $v$, которую и можно взять в качестве $f$
во втором утверждении теоремы.

Итак, остается доказать первое утверждение теоремы.
Пусть у графа $\Gamma$ относительно его вершины $v$ имеется $(F,\lambda)$-пропагатор $p_v$ с бесконечным носителем.
Если у графа $\Gamma$ относительно вершины $v$ имеется также $(F,\lambda)$-пропагатор с конечным носителем, то
разность $p_v$ и этого пропагатора есть, очевидно, собственная функция оператора $A_{\Gamma,F}$, соответствующая собственному значению
$\lambda$, с бесконечным носителем. Будем поэтому предполагать, что у графа $\Gamma$ относительно вершины
$v$ отсутствуют $(F,\lambda)$-пропагаторы с конечным носителем. Обозначая для каждой вершины $u$ графа $\Gamma$ через $\boldsymbol{\alpha}_u$ столбец матрицы
${\bf A}_{\Gamma,F} - \lambda {\bf E}$, соответствующий $u$, последнее условие, как легко понять,
можно переформулировать следующим образом:
для произвольного конечного подмножества $U$ множества $V(\Gamma)$ линейная оболочка множества $\{\boldsymbol{\alpha}_u: u \in U\}$
в векторном пространстве над $F$ столбцов
не содержит столбец $\boldsymbol{\delta}_v$ (у которого на месте, соответствующем вершине $v$ стоит $1$, а на остальных местах стоят нули).
Другими словами, для произвольного конечного подмножества $U$ множества $V(\Gamma)$
система линейных уравнений (над $F$ с неизвестными $x_u, u \in U$), определяемая посредством
$$\sum_{u \in U} x_u\boldsymbol{\alpha}_u = \boldsymbol{\delta}_v,$$
не имеет решений в $F$, что согласно теореме Теплица (см. предложение~\ref{p2.1}) означает наличие
такого конечного подмножества $U'$ множества $V(\Gamma)$  и таких $b_{u'} \in F$, где $u' \in U'$,
что, с одной стороны,
\begin{equation}
\label{eq5.1}
\sum_{u' \in U'}b_{u'}\alpha_{u',u} = 0
\end{equation}
для каждой вершины $u \in U$, где $\alpha_{u',u}, u' \in U',$ --- элемент столбца $\boldsymbol{\alpha}_u$ в строке, соответствующей вершине $u'$,
но, с другой стороны,
\begin{equation}
\label{eq5.2}
\sum_{u' \in U'}b_{u'}\delta_{u',v} \not = 0,
\end{equation}
где $\delta_{u',v} = 0$ при $u' \not = v$ и $\delta_{v,v} = 1$ при $v \in U'$.
В силу~\eqref{eq5.2} имеем $v \in U'$ и $b_v \not = 0$, а потому согласно~\eqref{eq5.1}
строка, соответствующая вершине $v$, $V(\Gamma) \times U$-подматрицы матрицы
${\bf A}_{\Gamma,F} - \lambda {\bf E}$
является линейной комбинацией других строк этой подматрицы.
Ввиду симметричности матрицы ${\bf A}_{\Gamma,F} - \lambda {\bf E}$ это
также означает, что столбец, соответствующий вершине $v$, $U \times V(\Gamma)$-подматрицы матрицы
${\bf A}_{\Gamma,F} - \lambda {\bf E}$
является линейной комбинацией других столбцов этой подматрицы. С учетом произвольности выбора конечного подмножества $U$ множества $V(\Gamma)$
отсюда следует, что каждая конечная подсистема системы линейных уравнений (над $F$ с неизвестными $y_u$, $u \in V(\Gamma) \setminus \{v\})$, определяемой посредством
\begin{equation}
\label{eq5.3}
\sum_{u \in V(\Gamma)\setminus \{v\}} y_u\boldsymbol{\alpha}_u = \boldsymbol{\alpha}_v,
\end{equation}
разрешима в $F$, откуда согласно теореме Теплица (см. предложение~\ref{p2.2}) вытекает
разрешимость в $F$ этой системы линейных уравнений.

Пусть $y_u' \in F, u \in V(\Gamma)\setminus \{v\}$, --- некоторое решение в $F$ системы линейных уравнений,
задаваемой посредством~\eqref{eq5.3}. Тогда функция $f \in F^{V(\Gamma)}$, определяемая посредством
$f(v) = 1$ и $f(u) = - y_u'$ для всех $u \in V(\Gamma)\setminus \{v\}$, есть собственная функция оператора $A_{\Gamma,F}$,
соответствующая собственному значению $\lambda$. Для завершения доказательства теоремы~\ref{t5.2}
остается заметить, что носитель функции $f$ бесконечен в силу теоремы~\ref{t4.1},
так как  у графа $\Gamma$ относительно вершины $v$ имеется $(F,\lambda)$-пропагатор $p_v$.

\end{proof}

Поскольку для
локально конечного графа $\Gamma$, поля $F$  с абсолютным значением $|.|_{\rm v}$ и $\lambda \in F$
из наличия у оператора $A_{\Gamma,F}$ собственной функции, соответствующей собственному значению $\lambda$, с бесконечным носителем
вытекает наличие у $\Gamma$ относительно произвольной содержащейся в $V(\Gamma) \setminus S_{F,\lambda}(\Gamma)$
вершины $(F,\lambda)$-пропагатора с бесконечным носителем, то теорема~\ref{t5.2} с учетом следствия~\ref{c4.1} влечет

\begin{cor}
\label{c5.1}

Пусть $\Gamma$ --- локально конечный граф, $F$ --- поле  с абсолютным значением $|.|_{\rm v}$ и $\lambda \in F$.
Если все $(F,\lambda)$-пропагаторы графа $\Gamma$ относительно некоторой вершины из $V(\Gamma) \setminus S_{F,\lambda}(\Gamma)$
имеют конечные носители, то все $(F,\lambda)$-пропагаторы графа $\Gamma$ относительно любой вершины из $V(\Gamma) \setminus S_{F,\lambda}(\Gamma)$
имеют конечные носители,
причем все $(F,\lambda)$-пропагаторы графа $\Gamma$ относительно любой фиксированной вершины из $V(\Gamma) \setminus S_{F,\lambda}(\Gamma)$ имеют
одно и то же ограничение на  $V(\Gamma) \setminus S_{F,\lambda}(\Gamma)$.

\end{cor}

Следующая теорема дает для локально конечного графа $\Gamma$, поля $F$
с произвольным абсолютным значением $|.|_{\rm v}$ и $\lambda \in F$
несколько характеризаций, включая характеризации в терминах $(F,\lambda)$-пропагаторов графа $\Gamma$,
свойства инъективности оператора $A_{\Gamma,F} - \lambda E$ или, другими словами, свойства
$\lambda$ не быть собственным значением оператора смежности $A_{\Gamma,F}$.

\begin{theorem}
\label{t5.3}

Пусть $\Gamma$ --- локально конечный граф, $F$ --- поле  с абсолютным значением $|.|_{\rm v}$ и $\lambda \in F$.
Тогда следующие условия $(i)-(viii)$ равносильны$:$

$(i)$ Элемент $\lambda$ не является собственным значением оператора смежности $A_{\Gamma,F}$ графа $\Gamma$,
или, другими словами, оператор $A_{\Gamma,F} - \lambda E$ инъективен.

$(ii)$ Оператор $A_{\Gamma,F} - \lambda E$ биективен.

$(iii)$ У графа $\Gamma$ относительно некоторой его вершины имеется единственный $(F,\lambda)$-пропагатор.

$(iv)$ $L_{F,\lambda}(\Gamma) = V(\Gamma)$.

$(v)$ У графа $\Gamma$ относительно каждой его вершины $v$ имеется единственный $(F,\lambda)$-пропагатор $p_v$
и этот пропагатор имеет конечный носитель, причем для произвольных $v_1,v_2 \in V(\Gamma)$ имеет место равенство
$p_{v_1}(v_2) = p_{v_2}(v_1)$ $($в частности, каждая вершина графа $\Gamma$ принадлежит лишь конечному
числу носителей пропагаторов $p_u$, $u \in V(\Gamma)$$)$.

$(vi)$ Матрица ${\bf A}_{\Gamma,F} - \lambda {\bf E}$ имеет в алгебре $M^{rcf}_{V(\Gamma)\times V(\Gamma)}(F)$
обратную $($симмет\-рическую$)$ матрицу.

$(vii)$ Оператор $A_{\Gamma,F} - \lambda E$ отображает подпространство $F^{V(\Gamma)*}$ пространства $F^{V(\Gamma)}$
сюръективно на себя.

$(viii)$ Оператор $A_{\Gamma,F} - \lambda E$ отображает подпространство $F^{V(\Gamma)*}$ пространства $F^{V(\Gamma)}$
биективно на себя.

\end{theorem}

\begin{proof}

Очевидно, что условие $(v)$ влечет условие $(iv)$, а условие $(vi)$
влечет условие $(v)$. Таким образом, для доказательства теоремы достаточно доказать равносильность условий
$(i)$, $(ii)$, $(iii)$, $(iv)$, $(vi)$, $(vii)$, $(viii)$, что и делается далее.

Предположим, что имеет место $(iv)$. Тогда у графа $\Gamma$ относительно каждой его вершины $u$
имеется $(F,\lambda)$-пропагатор $p_u$ с конечным носителем. Пусть ${\bf P}$ --- матрица из  $M^{cf}_{V(\Gamma)\times V(\Gamma)}(F)$,
у которой для любых $u_1, u_2 \in V(\Gamma)$ элемент на пересечении строки, соответствующей $u_1$, и столбца, соответствующего $u_2$,
равен $p_{u_2}(u_1)$. Тогда $({\bf A}_{\Gamma,F} - \lambda {\bf E}) {\bf P} =  {\bf E}$ и (с учетом симметричности матрицы ${\bf A}_{\Gamma,F}$)
также $ {\bf P}^t({\bf A}_{\Gamma,F} - \lambda {\bf E}) =  {\bf E}$. Отсюда (с учетом
${\bf A}_{\Gamma,F} - \lambda {\bf E} \in M^{rcf}_{V(\Gamma)\times V(\Gamma)}(F)$) следует
$${\bf P} = ({\bf P}^t({\bf A}_{\Gamma,F} - \lambda {\bf E})) {\bf P} =  {\bf P}^t(({\bf A}_{\Gamma,F} - \lambda {\bf E}){\bf P}) =  {\bf P}^t.$$
Таким образом, ${\bf P}$ --- симметрическая обратная в алгебре $M^{rcf}_{V(\Gamma)\times V(\Gamma)}(F)$ к ${\bf A}_{\Gamma,F} - \lambda {\bf E}$
матрица. Итак, из $(iv)$ следует  $(vi)$. Поскольку условие $(vi)$ влечет очевидным образом условие  $(iv)$,
заключаем, что условия  $(iv)$ и $(vi)$ равносильны.

Далее, очевидно, что условие $(iv)$ равносильно условию  $(vii)$, которое следует из условия $(viii)$.
Очевидно также, что условие $(vi)$ влечет условия  $(ii)$, $(iii)$, $(viii)$, а каждое из условий  $(ii)$ и $(iii)$
влечет условие $(i)$.

Таким образом, для завершения доказательства теоремы~\ref{t5.3} достаточно показать, что из условия $(i)$ следует условие $(iv)$.
Но это вытекает из теоремы~\ref{t4.2}.

\end{proof}

\begin{remark}
\label{r5.3}
При выполнении эквивалентных условий $(i)-(viii)$ обратный к
$A_{\Gamma,F} - \lambda E$ оператор
непрерывен. Это следует, например, из $(vi)$.

\end{remark}

Для локально конечного графа $\Gamma$, поля $F$  с абсолютным значением $|.|_{\rm v}$ и $\lambda \in F$ ослаблением условия инъективности оператора
$A_{\Gamma,F} - \lambda E$ является не только условие отсутствия у оператора $A_{\Gamma,F}$
собственной функции, соответствующей собственному значению $\lambda$, с {\it конечным} носителем (и равносильные ему условия $(i)-(v)$
из теоремы~\ref{t5.1}), но и условие отсутствия у оператора $A_{\Gamma,F}$ собственной функции, соответствующей собственному значению
$\lambda$, с {\it бесконечным} носителем (образно выражаясь, условие ``отсутствия коллапса'' функций
с бесконечным носителем). Как показывает следующая теорема, последнее условие равносильно
условию отсутствия в $F^{V(\Gamma)}$
функции с бесконечным носителем, образ которой под действием оператора $A_{\Gamma,F} - \lambda E$ имеет конечный
носитель, и оба эти условия равносильны условию, естественным образом формулируемому
в терминах $(F,\lambda)$-пропагаторов графа $\Gamma$.

\begin{theorem}
\label{t5.4}

Пусть $\Gamma$ --- локально конечный граф, $F$ --- поле  с абсолютным значением $|.|_{\rm v}$ и $\lambda \in F$.
Тогда следующие условия $(i)-(iii)$ равносильны$:$

$(i)$ У оператора $A_{\Gamma,F}$ отсутствуют собственные функции, соответствующие собственному значению $\lambda$,
с бесконечным носителем.

$(ii)$ Оператор $A_{\Gamma,F} - \lambda E$ отображает каждую функцию из $F^{V(\Gamma)}$
с бесконечным носителем в функцию также с бесконечным носителем.

$(iii)$ Множество $S_{F,\lambda}(\Gamma)$ конечно, а относительно некоторой $($или, что согласно следствию~$\ref{c5.1}$ равносильно, любой$)$ вершины из
$V(\Gamma) \setminus S_{F,\lambda}(\Gamma)$ у графа $\Gamma$ имеются лишь $(F,\lambda)$-пропагаторы с конечным носителем.

\end{theorem}

\begin{proof}
Очевидным образом $(ii)$ влечет $(i)$. Кроме того, из предложения~\ref{p5.2} и теоремы~\ref{t5.2}
следует, что $(i)$ влечет $(iii)$. Таким образом, для доказательства теоремы~\ref{t5.4} остается показать, что $(iii)$ влечет $(ii)$.

Предположим, что имеет место $(iii)$. Если граф $\Gamma$ конечен, то $(ii)$ выполняется тривиальным образом.
Предположим поэтому, что граф $\Gamma$ бесконечен. Для каждой вершины $u \in V(\Gamma) \setminus S_{F,\lambda}(\Gamma)$ пусть $p_u$ --- некоторый
$(F,\lambda)$-пропагатор с конечным носителем графа $\Gamma$ относительно $u$. Заменяя в случае необходимости $(F,\lambda)$-пропагатор
$p_u$ для $u \in V(\Gamma) \setminus S_{F,\lambda}(\Gamma)$ его ограничением на основную компоненту $\Gamma^2$-связности его расширенного
носителя (см. \S~\ref{s4}), будем, не теряя общности,
предполагать $\Gamma^2$-связность расширенного носителя $p_u$ для каждой вершины $u \in V(\Gamma) \setminus S_{F,\lambda}(\Gamma)$.
Из предложения~\ref{p4.6} следует, что каждая вершина графа
$\Gamma$ принадлежит лишь конечному числу расширенных носителей $(F,\lambda)$-пропагаторов
$p_u$, $u \in V(\Gamma) \setminus S_{F,\lambda}(\Gamma)$.
(Действительно, если некоторая вершина графа
$\Gamma$ принадлежит бесконечно многим расширенным носителям $(F,\lambda)$-пропагаторов $p_u$, $u \in V(\Gamma) \setminus S_{F,\lambda}(\Gamma)$,
то по предложению~\ref{p4.6} у оператора $A_{\Gamma,F}$ имеется собственная функция, соответствующая собственному значению $\lambda$,
с бесконечным носителем. Если $f$ --- такая функция, то для произвольной вершины $u \in V(\Gamma) \setminus S_{F,\lambda}(\Gamma)$
функция $p_u + f$ является $(F,\lambda)$-пропагатором с бесконечным носителем графа $\Gamma$ относительно $u$, что противоречит $(iii)$.)

Обозначим через ${\bf P}$ матрицу из $M^{cf}_{V(\Gamma)\times V(\Gamma)}(F)$,
у которой для $u_1, u_2 \in V(\Gamma)$ элемент на пересечении строки, соответствующей $u_1$, и столбца, соответствующего $u_2$,
равен $p_{u_2}(u_1)$ при $u_2 \in  V(\Gamma) \setminus S_{F,\lambda}(\Gamma)$ и равен $\delta_{u_1,u_2}$ при
$u_2 \in  S_{F,\lambda}(\Gamma)$ (где $\delta_{u_1,u_2} = 1$, если $u_1 = u_2$,
и $\delta_{u_1,u_2} = 0$, если $u_1 \not = u_2$). Поскольку каждая вершина графа
$\Gamma$ принадлежит лишь конечному числу носителей функций $p_u$, $u \in V(\Gamma) \setminus S_{F,\lambda}(\Gamma)$, то
${\bf P}$, а следовательно, и ${\bf P}^t$ принадлежат $M^{rcf}_{V(\Gamma)\times V(\Gamma)}(F)$.
Далее, у матрицы ${\bf P}^t ({\bf A}_{\Gamma,F} - \lambda {\bf E})$
ее $(V(\Gamma) \setminus S_{F,\lambda}(\Gamma)) \times V(\Gamma)$-подматрица совпадает с
$(V(\Gamma) \setminus S_{F,\lambda}(\Gamma)) \times V(\Gamma)$-подматрицей
единичной $V(\Gamma) \times V(\Gamma)$-матрицы, а ее
$S_{F,\lambda}(\Gamma) \times V(\Gamma)$-подматрица совпадает с
$S_{F,\lambda}(\Gamma) \times V(\Gamma)$-подматрицей
матрицы ${\bf A}_{\Gamma,F} - \lambda {\bf E}$. Поэтому
(с учетом конечности $S_{F,\lambda}(\Gamma)$) для любой функции $f \in F^{V(\Gamma)}$ с бесконечным носителем,
обозначая через ${\bf f}$ столбец (с элементами из $F$), строки которого соответствуют $V(\Gamma)$, а элемент
в строке, соответствующей $u \in V(\Gamma)$, равен $f(u)$, получаем, что столбец
$({\bf P}^t ({\bf A}_{\Gamma,F} - \lambda {\bf E})){\bf f}$ имеет (как и ${\bf f}$) бесконечно много ненулевых элементов, что в силу
${\bf P}^t \in M^{rcf}_{V(\Gamma)\times V(\Gamma)}(F)$ влечет бесконечность числа ненулевых
элементов также у столбца
$({\bf A}_{\Gamma,F} - \lambda {\bf E}){\bf f}$, т. е. бесконечность носителя функции
$(A_{\Gamma,F} - \lambda E)(f)$.
Таким образом, имеет место $(ii)$.

\end{proof}

\section{Пропагаторы и суммы по путям}
\label{s6}

В этом параграфе мы рассматриваем локально конечные {\it связные} графы.
Полученные результаты применимы к компонентам связности произвольных локально конечных графов.

В предыдущих параграфах было показано, что для локально конечного графа $\Gamma$ и поля $F$ с
абсолютным значением $|.|_{\rm v}$ ряд представляющих интерес свойств
оператора смежности $A_{\Gamma,F}$ (и, в частности, свойств его собственных функций)
естественным образом формулируется в терминах $(F,\lambda)$-пропагато\-ров графа $\Gamma$
для $\lambda \in F$. Целесообразность такой переформулировки находит подтверждение
в наличии для пропагаторов (в отличие от собственных функций) конструкции их построения,
использующей определенные суммы по путям графа, которая описывается (для связных графов) в настоящем параграфе.
Суммы по путям и соответствующие производящие функции широко применяются
в теории графов (как конечных, см., например,~\cite{God}, так и бесконечных, например,
при исследовании случайных блужданий на них, см.~\cite{Wo}).
В описываемой далее конструкции пропагаторы получаются в виде наборов производящих
функций определенных сумм по путям, а в случае $F = \mathbb{C}$ также наборов
аналитических продолжений этих производящих функций.
Прежде, чем перейти к описанию
конструкции, напомним, что если $\Gamma$ --- локально конечный связный граф, $F$ --- поле с
абсолютным значением $|.|_{\rm v}$ и $\lambda \in F$,
то согласно утверждению 2) предложения~\ref{p3.1} и теореме~\ref{t5.2}
(или теореме~\ref{t5.3}) для доказательства того, что $\lambda$ --- собственное значение оператора смежности $A_{\Gamma,F}$, достаточно
доказать существование у $\Gamma$ относительно какой-либо его вершины $(F',\lambda)$-пропагатора с бесконечным носителем
для какого-либо расширения $F'$ поля $F_0(\lambda)$, где $F_0$ --- простое подполе поля $F$, с произвольным
абсолютным значением $|.|_{\rm v'}$. Выбор содержащего $F_0(\lambda)$ поля $F'$ ``большим'' способно облегчить нахождение $(F',\lambda)$-пропагатора
с бесконечным носителем.

Перейдем к описанию конструкции пропагаторов на основе использования определенных сумм по путям графа.

Пусть $\Gamma$ --- локально конечный связный граф, $K$ --- поле и $x$ --- независимая переменная над $K$.
Как обычно, для произвольных $v, w \in V(\Gamma)$ производящей функцией над $K$ для числа путей графа $\Gamma$, начинающихся в $v$ и заканчивающихся в $w$,
называется определяемый следующим образом элемент $K$-алгебры $K[[x]]$ формальных степенных рядов от $x$ с коэффициентами из $K$:
\begin{equation}
\label{eq6.1}
W_{\Gamma,K}(x,v,w) := \sum_{n \in \mathbb{Z}_{\geq 0}} (N_{\Gamma}(v,w,n)\cdot 1_K)\ x^n,
\end{equation}
где, для каждого $n \in \mathbb{Z}_{\geq 0}$, $N_{\Gamma}(v,w,n)$ есть число путей длины $n$ графа $\Gamma$, начинающихся
в $v$ и заканчивающихся в $w$ (в дальнейшем, полагаясь на контекст, мы иногда будем опускать единичный элемент $1_K$ поля $K$ в записи $m \cdot 1_K$,
где $m \in \mathbb{Z}_{\geq 0}$).
Отметим, что
\begin{equation}
\label{eq6.2}
W_{\Gamma,K}(x,v,w) = W_{\Gamma,K}(x,w,v).
\end{equation}

$K$-алгебра $K[[x]]$ формальных степенных рядов от $x$ с коэффициентами из $K$ естественным образом
вкладывается в свое поле частных $K((x))$. Для произвольного $v \in V(\Gamma)$ определим элемент
$p_{\Gamma,K,x,v} \in K[[x]]^{V(\Gamma)} \subseteq K((x))^{V(\Gamma)}$, полагая
$$p_{\Gamma,K,x,v}(w) := - x W_{\Gamma,K}(x,v,w)\ \mbox{для всех}\ w \in V(\Gamma).$$
Для $v, w \in V(\Gamma)$ элемент $W_{\Gamma,K}(x,v,w)$ поля $K((x))$ можно
интерпретировать следующим образом: каждому пути графа $\Gamma$, начинающемуся в $v$ и заканчивающемся в $w$, приписывается вес из $K[[x]]$, равный
$x$ в степени длина этого пути, после чего $W_{\Gamma,K}(x,v,w)$
совпадает с суммой (в смысле~\cite[гл. IV, $\S\ 5, {\rm n}^\circ 4$]{Burb}) весов по всем путям графа $\Gamma$, начинающимся в $v$ и заканчивающимся в $w$.
С учетом тривиального наблюдения, что для произвольных $v, w \in V(\Gamma)$ и $n \in \mathbb{Z}_{\geq 1}$ в графе $\Gamma$
путь длины $n$ из $v$ в $w$ есть в точности произведение (т. е. последовательное прохождение) пути длины $n-1$ из $v$ в какую-нибудь вершину
из $\Gamma(w)$ и пути длины $1$ из этой последней вершины в вершину $w$, и потому
\begin{equation}
\label{eq6.3}
N_{\Gamma}(v,w,n) =
\sum_{w' \in \Gamma(w)} N_{\Gamma}(v,w',n-1),
\end{equation}
отсюда вытекает справедливость утверждения $1)$ следующего предложения~\ref{p6.1}. Справедливость утверждения $2)$  предложения~\ref{p6.1}
вытекает из~\eqref{eq6.2}, а справедливость утверждения $3)$ следует из связности графа $\Gamma$.

\begin{propos}
\label{p6.1}

Пусть $\Gamma$ --- локально конечный связный граф,  $K$ --- поле и $x$ --- независимая переменная над $K$. Снабдим $K((x))$ произвольным
абсолютным значением.
Тогда $($во введенных выше обозначениях$)$ справедливы следующие утверждения$:$

$1)$ Для $v \in V(\Gamma)$ функция $p_{\Gamma,K,x,v} \in K((x))^{V(\Gamma)}$  является
$(K((x)),1/x)$-пропа\-гатором графа $\Gamma$ относительно вершины $v$.

$2)$ Для $v, v' \in V(\Gamma)$  имеем $p_{\Gamma,K,x,v}(v') = p_{\Gamma,K,x,v'}(v)$.

$3)$ Если $K$ --- поле характеристики $0$, то для $v \in V(\Gamma)$ носитель $p_{\Gamma,K,x,v}$ совпадает с $V(\Gamma)$.

\end{propos}

\begin{remark}
\label{r6.1}
Утверждение 1) предложения~\ref{p6.1} можно рассматривать
как своего рода указание на специфичность случая $\lambda = 0$ при исследовании
$(F,\lambda)$-про\-пагаторов локально конечных связных графов $\Gamma$.

\end{remark}

\begin{remark}
\label{r6.2}

В случае {\it регулярного} степени $d \in \mathbb{Z}_{\geq 1}$ связного графа $\Gamma$
и поля частных
$K((x))$ (с произвольным абсолютным значением)
$K$-алгебры $K[[x]]$ формальных степенных рядов от $x$ с коэффициентами из произвольного поля $K$
в роли, аналогичной той, в которой выше для построения
значений пропагаторов
$p_{\Gamma,K,x,v}(w),$ $v,w \in V(\Gamma),$ использовалась производящая функция
$$W_{\Gamma,K}(x,v,w) = \sum_{n \in \mathbb{Z}_{\geq 0}} N_{\Gamma}(v,w,n)\cdot 1_K x^n,$$
(см.~\eqref{eq6.1}),
может быть использована производящая функция
$$\widehat{W}_{\Gamma,K}(x,v,w) := \sum_{n \in \mathbb{Z}_{\geq 0}} \widehat{N}_{\Gamma}(v,w,n)\cdot 1_K x^n \in K[[x]],$$
где $v, w \in V(\Gamma)$ и для $n \in \mathbb{Z}_{\geq 0}$ через $\widehat{N}_{\Gamma}(v,w,n)$ обозначается число таких
путей $v = u_0,...,u_n = w$ графа $\Gamma$, что $u_{i-1} \not = u_{i+1}$ для всех $0 < i < n$.
В качестве аналога~\eqref{eq6.3} выступает соотношение
$$\widehat{N}_{\Gamma}(v,w,n) = (\sum_{w' \in \Gamma(w)}\widehat{N}_{\Gamma}(v,w',n-1)) -
(d - 1) \widehat{N}_{\Gamma}(v,w,n-2),$$
справедливое для всех $v, w \in V(\Gamma)$ и $n \in \mathbb{Z}_{\geq 0}$,
кроме случая $v = w \in V(\Gamma)$, $n = 0$ и случая $v = w \in V(\Gamma)$, $n = 2$.
С использованием этого соотношения можно показать, что для $v \in V(\Gamma)$ функция $\widehat{p}_{\Gamma,K,x,v} \in K((x))^{V(\Gamma)}$, определяемая посредством
$$\widehat{p}_{\Gamma,K,x,v}(w) = - \frac{x}{1 - x^2} \widehat{W}_{\Gamma,K}(x,v,w)\ \mbox{для всех}\ w \in V(\Gamma),$$
является $(K((x)),(1+(d-1)x^2)/x)$-пропагатором графа $\Gamma$ относительно вершины $v$.

Однако даже в случае регулярного графа $\Gamma$ такое использование $\widehat{W}_{\Gamma,K}(x,v,w)$
вместо $W_{\Gamma,K}(x,v,w)$
не представляет интереса в принципиальном плане, поскольку
$$\widehat{W}_{\Gamma,K}(x,v,w) = \frac{1-x^2}{1+(d-1)x^2} W_{\Gamma,K}(\frac{x}{1+(d-1)x^2},v,w)$$
(см.~\cite[с. 72]{God},~\cite{Bar}).

\end{remark}

Из утверждения 3) предложения~\ref{p6.1} с учетом теоремы~\ref{t5.2} (или теоремы~\ref{t5.3}) и утверждения 2) предложения~\ref{p3.1} (а также равенства
$K(x) = K(1/x)$, где $1/x$ ---
также независимая переменная над $K$) вытекает справедливость следующего утверждения.
Если $\Gamma$ --- бесконечный локально конечный связный граф  и $x$ --- независимая переменная над $\mathbb{Q}$,
то $x$ является собственным значением оператора смежности $A_{\Gamma,\mathbb{Q}(x)}$ $($для произвольного
абсолютного значения на поле рациональных функций $\mathbb{Q}(x)$$)$.
Отсюда (с учетом утверждения 1) предложения~\ref{p3.1}) вытекает приводимая ниже теорема~\ref{t6.1}. Ее формулировке мы
предпошлем описание контекста, в котором ее естественно рассматривать.

Пусть $\Gamma$ --- бесконечный локально конечный связный граф, $F$ --- поле характеристики $0$ с произвольным абсолютным
значением и $\mathbb{Q}$ --- простое подполе поля $F$.
Как было только что доказано, у оператора смежности $A_{\Gamma,\mathbb{Q}(x)}$ имеется собственная функция
$f \in \mathbb{Q}(x)^{V(\Gamma)}$, соответствующая собственному значению $x$.
Для каждой вершины $w$ графа $\Gamma$ элемент $f(w)$ поля $\mathbb{Q}(x)$
однозначно записывается в виде $f_{1}(w)/f_{2}(w)$, где $f_{1}(w) \in \mathbb{Q}[x]$,
$0 \not = f_{2}(w) \in \mathbb{Q}[x]$, старший коэффициент
$f_{2}(w)$ единичен и при этом $\deg (f_{1}(w))$
и $\deg (f_{2}(w))$ минимально возможные для таких представлений $f(w)$.
Заменяя в случае необходимости $f$ ее произведением на
$\text{НОД}(f_2(w): w \in V(\Gamma))/\text{НОД}(f_1(w): w \in V(\Gamma)) \in \mathbb{Q}(x) \setminus \{0\}$,
можно при этом предполагать, что
$\text{НОД}(f_1(w): w \in V(\Gamma)) = 1$ и $\text{НОД}(f_2(w): w \in V(\Gamma)) = 1$.
Если $\lambda$ --- такой ненулевой элемент поля $F$,
что $(f_{1}(w))(\lambda) \not = 0$ для некоторой вершины $w \in V(\Gamma)$ и
$(f_{2}(w))(\lambda) \not = 0$ для всех вершин $w \in V(\Gamma)$
(что имеет место, например, в случае трансцендентного над $\mathbb{Q}$ элемента $\lambda$), то
функция из $\mathbb{Q}(\lambda)^{V(\Gamma)} \leq F^{V(\Gamma)}$, значение которой в каждой вершине $w$ графа $\Gamma$ равно
значению, принимаемому рациональной функцией $f_{1}(w)/f_{2}(w) \in \mathbb{Q}(x)$ при $x = \lambda$,
очевидным образом является собственной функцией оператора смежности $A_{\Gamma,\mathbb{Q}(\lambda)}$,
соответствующей собственному значению $\lambda$.
В частности, справедлива следующая

\begin{theorem}
\label{t6.1}
Пусть $\Gamma$ --- бесконечный локально конечный связный граф, $F$ --- поле характеристики $0$ с произвольным абсолютным значени\-ем
и $\lambda$ --- элемент поля $F$, трансцендентный над его простым подполем $\mathbb{Q}$. Тогда
$\lambda$ является собственным значением оператора смежности $A_{\Gamma,\mathbb{Q}(\lambda)}$ $($а следовательно,
и оператора смежности $A_{\Gamma,F}$$)$.

\end{theorem}

\begin{remark}
\label{r6.3} В предположениях теоремы~\ref{t6.1} ясно, что носители всех соответствующих $\lambda$ собственных функций оператора смежности $A_{\Gamma,F}$ бесконечны и, более того, имеют лишь бесконечные
связные компоненты. (Действительно, $\lambda$ является собственным значением
оператора
$A_{\langle X \rangle_{\Gamma},F}$ для любой связной компоненты $X$ носителя любой
собственной функции оператора $A_{\Gamma,F}$, и потому $X$ не может быть конечной
ввиду трансцендентности $\lambda$ над $\mathbb{Q}$.)

\end{remark}

\begin{remark}
\label{r6.4}

Вообще говоря, аналог теоремы~\ref{t6.1} для поля $F$ положительной
характеристики не имеет места, см. \S~\ref{s8}, раздел~\ref{s8.1}.

\end{remark}

Вернемся к исследованию пропагаторов локально конечных связных графов.
Из предложения~\ref{p6.1} и утверждения 2) предложения~\ref{p4.1} вытекает

\begin{propos}
\label{p6.2}

В обозначениях предложения~$\ref{p6.1}$ пусть $K_0$ --- простое подполе
поля $K$. Тогда у графа $\Gamma$ относительно произвольной его вершины $v$ существует
$(K_0(x),1/x)$-пропагатор с расширенным носителем, содержащимся в расширенном
носителе $(K((x)),1/x)$-пропага\-тора $p_{\Gamma,K,x,v}$ графа $\Gamma$ относительно $v$.

\end{propos}

Пусть $\Gamma$ --- локально конечный связный граф, $F$ --- поле с некоторым абсолютным значением, $F_0$ --- простое подполе поля $F$ и $x$ --- независимая переменная над $F$.
Согласно предложению~\ref{p6.2} для любого абсолютного значения на $F(x)$ у графа $\Gamma$ относительно произвольной его вершины $v$ имеется
$(F(x),1/x)$-пропагатор $p_v \in F_0(x)^{V(\Gamma)}$.
Для каждой вершины $w$ графа $\Gamma$ элемент $p_v(w)$ поля $F_0(x)$
однозначно записывается в виде $p_{v,1}(w)/p_{v,2}(w)$, где $p_{v,1}(w) \in F_0[x]$, $0 \not = p_{v,2}(w) \in F_0[x]$, старший коэффициент
$p_{v,2}(w)$ единичен и при этом $\deg (p_{v,1}(w))$
и $\deg (p_{v,2}(w))$ минимально возможные для таких представлений $p_v(w)$. Если теперь $\lambda$ --- такой ненулевой элемент поля $F$,
что $(p_{v,2}(w))(\lambda) \not = 0$ для всех $w \in V(\Gamma)$
(что имеет место, например, в случае трансцендентного над $F_0$ элемента $\lambda$), то
функция из $F_0(\lambda)^{V(\Gamma)} \leq F^{V(\Gamma)}$, значение которой в каждой вершине $w$ графа $\Gamma$ равно
значению, принимаемому рациональной функцией $p_{v,1}(w)/p_{v,2}(w) \in F_0(x)$ при $x = \lambda$,
очевидным образом является $(F_0(\lambda),1/{\lambda})$-пропагатором (а следовательно, и $(F,1/{\lambda})$-пропагатором)
графа $\Gamma$ относительно вершины $v$.
В частности, мы вновь (см. утверждение 1) следствия~\ref{c4.2}) получаем, что для $\lambda \in F$ отсутствие
у графа $\Gamma$ относительно какой-либо его вершины
$(F,\lambda)$-пропагатора возможно лишь в случае алгебраического над $F_0$
элемента $\lambda$.

\begin{remark}
\label{r6.5}

Из утверждения 1) предложения~\ref{p6.1} следует, что если $\Gamma$ --- локально конечный связный граф и $F$ --- поле
положительной характеристики, полное относительно (обязательно неархимедова) абсолютного значения $|.|_{\rm v}$, то для $v \in V(\Gamma)$ и $\lambda \in F$ c $|\lambda|_{\rm v} > 1$
функция $f: V(\Gamma) \to F$, определяемая посредством
$$f(w) = - \lambda^{-1} \sum_{n \in \mathbb{Z}_{\geq 0}} (N_{\Gamma}(v,w,n)\cdot 1_F)\ \lambda^{-n}\ \mbox {для всех}\ w \in V(\Gamma)$$
(где сумма справа обозначает очевидным образом существующий предел по $|.|_{\rm v}$ частичных сумм
$$\sum_{0 \leq n \leq n'} (N_{\Gamma}(v,w,n)\cdot 1_F)\ \lambda^{-n}$$
при $n' \to \infty$, $n' \in \mathbb{Z}_{\geq 0}$),
является $(F,\lambda)$-пропагатором графа $\Gamma$ относительно вершины $v$. Поскольку при сделанных предположениях элемент $\lambda$
трансцендентен над простым подполем поля $F$, то в рассматриваемом случае существование
$(F,\lambda)$-пропагатора графа $\Gamma$ относительно $v$ следует из утверждения 1) следствия~\ref{c4.2}, а также из предложения~\ref{p6.2}.
Назначение настоящего замечания --- предъявить ``явный вид'' одного из таких $(F,\lambda)$-пропагаторов.

\end{remark}

В оставшейся части настоящего параграфа мы разовьем изложенный подход
к построению пропагаторов локально конечных связных графов в случае поля $F = \mathbb{C}$
с обычным модулем $|.|$ комплексного числа в качестве абсолютного значения,
причем особое внимание будет уделено графам с ограниченными в совокупности
степенями вершин. Вначале мы приведем для этого случая некоторые простые следствия
предыдущих результатов параграфа.

Пусть $\Gamma$ --- локально конечный связный граф и $v \in V(\Gamma)$. Напомним, что для $w \in V(\Gamma)$
через $W_{\Gamma,\mathbb{C}}(x,v,w) \in \mathbb{C}((x))$ обозначается производящая функция над $\mathbb{C}$
для числа путей графа $\Gamma$, начинающихся в $v$ и заканчивающихся
в $w$ (см.~\eqref{eq6.1}).
Для $w \in V(\Gamma)$ и $0 \not = \lambda \in \mathbb{C}$ обозначим
через $W_{\Gamma,\mathbb{C}}(\lambda^{-1},v,w)$
сумму ряда, получающегося из $W_{\Gamma,\mathbb{C}}(x,v,w)$ подстановкой $\lambda^{-1}$ вместо $x$, {\it если
эта сумма определена} (т. е. получающийся ряд $\sum_{n \in \mathbb{Z}_{\geq 0}} N_{\Gamma}(v,w,n) \lambda^{-n}$
сходится), и {\it для таких} $\lambda$
положим
$$s_{\Gamma,\lambda,v}(w) := -\lambda^{-1} W_{\Gamma,\mathbb{C}}(\lambda^{-1},v,w).$$
Заметим, что для $w \in V(\Gamma)$ и $\lambda \in \mathbb{C}$ число $s_{\Gamma,\lambda,v}(w)$ определено тогда и только тогда,
когда определено число $s_{\Gamma,\lambda,w}(v)$, и для таких $w$ и $\lambda$ имеем
$$s_{\Gamma,\lambda,v}(w) = s_{\Gamma,\lambda,w}(v).$$
Обозначим через $D_{\Gamma,v}$ множество таких $0 \not = \lambda \in \mathbb{C}$, что
$s_{\Gamma,\lambda,v}(w)$ (или, эквивалентно, $W_{\Gamma,\mathbb{C}}(\lambda^{-1},v,w)$)
определено для всех $w \in V(\Gamma)$. Заметим, что в случае, когда степени вершин графа $\Gamma$
не превосходят некоторого натурального числа $d$, для произвольных $w \in V(\Gamma)$,
$n \in \mathbb{Z}_{\geq 0}$
имеем, очевидно, $N_{\Gamma}(v,w,n) \leq d^n$, и потому
$$\{c \in \mathbb{C} : |c| > d\} \subseteq D_{\Gamma,v}.$$
Более того, далее будет показано (см. замечание~\ref{r6.7}), что в случае такого графа $\Gamma$
для произвольных фиксированных $v \in V(\Gamma)$ и $\lambda \in \{c \in \mathbb{C} : |c| > d\}$
сумма
$$\sum_{w \in V(\Gamma)}|s_{\Gamma,\lambda,v}(w)|^2$$
конечна. Очевидно, что при $|V(\Gamma)| > 1$ из $\lambda \in D_{\Gamma,v}$ и
$\lambda \in \mathbb{R}_{> 0}$ следует $\lambda > 1$.

\medskip

В следующих предложениях~\ref{p6.3} и~\ref{p6.4} собраны несколько простых утверждений (некоторые из которых используются в дальнейшем)
относительно $s_{\Gamma,\lambda,v}(w)$
как функции аргумента $w \in V(\Gamma)$ для $\lambda \in D_{\Gamma,v}$.
Справедливость утверждения $1)$ предложения~\ref{p6.3} легко устанавливается с использованием~\eqref{eq6.3}
(ср. с доказательством
утверждения 1) предложения~\ref{p6.1}). Справедливость утверждений $2)$ и $3)$ предложения~\ref{p6.3}
легко следует из вещественности и неотрицательности коэффициентов $W_{\Gamma}(v,w,x)$,
а в случае утверждения 3) еще и выполнения неравенств $N_{\Gamma}(v,w,n) \leq d^n$
для всех $n \in \mathbb{Z}_{\geq 0}$.

\begin{propos}
\label{p6.3}

Пусть $\Gamma$ --- локально конечный связный граф, $v \in V(\Gamma)$ и
$\lambda \in D_{\Gamma,v}$ $($например,
это имеет место в случае, когда $\Gamma$ --- связный граф, степени вершин которого не превосходят некоторого натурального
числа $d$, $v \in V(\Gamma)$ и $|\lambda| > d$, см. выше$)$. Тогда справедливы следующие утверждения$:$

$1)$ Функция $s_{\Gamma,\lambda,v} \in \mathbb{C}^{V(\Gamma)}$ $($определенная выше$)$ является
$(\mathbb{C},\lambda)$-пропагато\-ром графа $\Gamma$ относительно вершины $v$.

$2)$ Если $w \in V(\Gamma)$ и $\lambda \in \mathbb{R}$, то $s_{\Gamma,\lambda,v}(w) \in \mathbb{R}$,
причем $s_{\Gamma,\lambda,v}(w) < 0$ при $\lambda > 0$.
Если, кроме того, $\lambda' \in \mathbb{R}$ и $\lambda' > \lambda > 0$, то $\lambda' \in D_{\Gamma,v}$
и $|s_{\Gamma,\lambda',v}(w)| < |s_{\Gamma,\lambda,v}(w)|$.

$3)$ Если степени вершин графа $\Gamma$ не превосходят некоторого натурального
числа $d$,
то при $|\lambda| > d$ имеем
$|s_{\Gamma,\lambda,v}(w)| \leq |s_{\Gamma,|\lambda|,v}(w)| \leq 1/(|\lambda| - d)$
для всех $w \in V(\Gamma)$.

\end{propos}

Пусть $\Gamma$ --- локально конечный связный граф и $v \in V(\Gamma)$. Тогда, если $\lambda \in D_{\Gamma,v}$ и
функция $s_{\Gamma,\lambda,v}$ имеет бесконечный носитель, то согласно утверждению 1) предложения~\ref{p6.3}
и теореме~\ref{t5.2} (или теореме ~\ref{t5.3}) число $\lambda$ является собственным значением оператора смежности графа $\Gamma$.
В силу утверждения 2) предложения~\ref{p6.3} для каждого положительного вещественного числа
$\lambda \in D_{\Gamma,v}$
(в частности, в случае когда степени вершин $\Gamma$ не превосходят некоторого натурального
числа $d$, для $\lambda \in \mathbb{R}_{> d}$) носитель функции $s_{\Gamma,\lambda,v}$ совпадает
с $V(\Gamma)$. Поэтому в случае бесконечного графа $\Gamma$
каждое такое вещественное число $\lambda$ является собственным значением оператора смежности графа $\Gamma$.
Мы уточним это заключение в утверждении 3) следующего предложения~\ref{p6.4}. Для его доказательства нам потребуются
утверждения 1) и 2) этого предложения.

\begin{propos}
\label{p6.4}

Пусть $\Gamma$ ---  локально конечный связный граф, $v \in V(\Gamma)$ и
$\lambda$ --- положительное вещественное число, содержащееся в $D_{\Gamma,v}$. $($Например,
это имеет место в случае, когда $\Gamma$ --- связный граф, степени вершин которого не превосходят некоторого натурального
числа $d$, $v \in V(\Gamma)$ и $\lambda \in \mathbb{R}_{> d}$.$)$
Тогда справедливы следующие утверждения$:$

$1)$ $\lambda \in D_{\Gamma,u}$ для всех $u \in V(\Gamma)$.

$2)$ $s_{\Gamma,\lambda,v}(w_1)/s_{\Gamma,\lambda,v}(w_2) \leq \lambda$ для всех $\{w_1,w_2\} \in E(\Gamma)$ $($напомним, что согласно
утверждению $2)$ предложения~$\ref{p6.3}$ для всех $w \in V(\Gamma)$ имеем
$s_{\Gamma,\lambda,v}(w) < 0$$)$.

$3)$ Если $\Gamma$ --- бесконечный граф, то существует такая собственная функция $f$ оператора смежности $A_{\Gamma, \mathbb{C}}$
графа $\Gamma$, соответствующая собственному значению $\lambda$, что $f(v) = 1$, $f(w) \in \mathbb{R}_{>0}$ для всех $w \in V(\Gamma)$
и $f(w_1)/f(w_2) \leq \lambda$ для всех $\{w_1,w_2\} \in E(\Gamma)$
$($и потому, в частности, $\lambda ^{-d_{\Gamma}(v,w)} \leq f(w) \leq \lambda ^{d_{\Gamma}(v,w)}$
для всех $w \in V(\Gamma)$$)$.

\end{propos}

\begin{proof}

Для доказательства утверждения 1) предложения достаточно заметить, что
для заданной $u \in V(\Gamma)$ и произвольных $w \in V(\Gamma)$ и $n \in \mathbb{Z}_{\geq d_{\Gamma}(u,v)}$
имеем $$N_{\Gamma}(u,w,n) \leq N_{\Gamma}(v,w,n + d_{\Gamma}(u,v)).$$

Для доказательства утверждения 2) предложения достаточно заметить, что в силу $N_{\Gamma}(v,w_1,n) \leq N_{\Gamma}(v,w_2,n+1)$, где
$n \in \mathbb{Z}_{\geq 0}$, имеем
$$\frac{\sum_{n \in \mathbb{Z}_{\geq 0}}N_{\Gamma}(v,w_1,n)\lambda^{-n}}{\sum_{n \in \mathbb{Z}_{\geq 0}}N_{\Gamma}(v,w_2,n)\lambda^{-n}} \leq
\frac{\lambda \sum_{n \in \mathbb{Z}_{\geq 0}}N_{\Gamma}(v,w_2,n+1)\lambda^{-n-1}}{\sum_{n \in \mathbb{Z}_{\geq 0}}N_{\Gamma}(v,w_2,n)\lambda^{-n}}
\leq \lambda.$$

Докажем утверждение 3) предложения.
Для произвольной вершины $u$ графа $\Gamma$ определим функцию $f_u \in \mathbb{C}^{V(\Gamma)}$,
полагая $f_u(w) = s_{\Gamma,\lambda,u}(w)/s_{\Gamma,\lambda,u}(v)$ для всех $w \in V(\Gamma)$.
Так как $\lambda \in D_{\Gamma,u}$ в силу уже доказанного утверждения 1) предложения~\ref{p6.4},
то из утверждений 2) и 1) предложения~\ref{p6.3} следует, что функция $f_u$
корректно определена, причем принимает только вещественные положительные значения (поскольку
функция $s_{\Gamma,\lambda,u}$ принимает вещественные отрицательные значения всюду на $V(\Gamma)$),
$f_u(v) = 1$ и $(A_{\Gamma,\mathbb{C}}(f_u))(w) = \lambda f_u(w)$ для всех
$w \in B_{\Gamma}(v,d_{\Gamma}(u,v) - 1)$.
Кроме того, в силу уже доказанного утверждения 2) предложения~\ref{p6.4}
имеем $f_u(w_1)/f_u(w_2) \leq \lambda$ для всех $\{w_1,w_2\} \in E(\Gamma)$.
Поскольку в силу бесконечности связного локально конечного графа $\Gamma$ вершина $u$ может
быть выбрана на расстоянии от $v$, превосходящем любое наперед заданное натуральное число,
то отсюда следует существование для произвольного натурального числа $r$ такой
функции $f_r \in \mathbb{C}^{V(\Gamma)}$, что
$f_r(w) \in \mathbb{R}_{>0}$ для всех $w \in V(\Gamma)$,
$f_r(v) = 1$, $(A_{\Gamma,\mathbb{C}}(f_r))(w) = \lambda f_r(w)$ для всех $w \in B_{\Gamma}(v,r)$
и, наконец, $f_r(w_1)/f_r(w_2) \leq \lambda$ для всех $\{w_1,w_2\} \in E(\Gamma)$.
Ясно, что последовательность функций $(f_r)_{r \in \mathbb{Z}_{\geq 1}}$ содержит подпоследовательность,
сходящуюся в топологии поточечной сходимости к некоторой функции $f$, которая очевидным образом
обладает всеми указанными в утверждении 3) предложения~\ref{p6.4} свойствами.

\end{proof}

Дополнительные углубленные результаты получаются с использованием, по существу, аналитических продолжений
ранее рассматриваемых производящих функций $s_{\Gamma,\lambda,v}$ в случае,
когда $\Gamma$ --- связный граф, степени вершин которого
не превосходят некоторого натурального числа $d$ (и, по-прежнему, $F = \mathbb{C}$ с обычным
абсолютным значением). Как будет показано, описываемая далее конструкция позволяет строить
{\it квадратично суммируемые $(\mathbb{C},\lambda)$-пропагаторы графа $\Gamma$ для всех
$\lambda \in \mathbb{C} \setminus (\mathbb{R}_{\geq - d} \cap \mathbb{R}_{\leq d})$.}

Итак, пусть $\Gamma$ --- связный граф, степени вершин которого
не превосходят некоторого натурального числа $d$ и $F = \mathbb{C}$ с обычным
абсолютным значением.
Пусть $l_2(V(\Gamma))$ --- содержащееся в $\mathbb{C}^{V(\Gamma)}$ гильбертово пространство
квадратично суммируемых функций со скалярным произведением $$(f_1,f_2) =
\sum_{u \in V(\Gamma)}f_1(u)\overline{f_2(u)}\ \mbox {для всех}\ f_1, f_2 \in l_2(V(\Gamma)),$$
где $\overline{f_2(u)}$ для $u \in V(\Gamma)$ --- число комплексно сопряженное с $f_2(u)$
(заметим, что топология на $l_2(V(\Gamma))$ не индуцируется
топологией топологического векторного пространства
$\mathbb{C}^{V(\Gamma)}$), и пусть ${\cal L}(l^2(V(\Gamma)))$ --- множество всюду определенных
ограниченных линейных операторов из гильбертова пространства $l_2(V(\Gamma))$ в себя.
Используемая далее терминология, касающаяся операторов из ${\cal L}(l^2(V(\Gamma)))$ и, в частности,
их спектров, совпадает с используемой в теории $C^*$-алгебр, к каковым ${\cal L}(l^2(V(\Gamma)))$
(с естественными операциями) относится. Так, для произвольного оператора $A \in {\cal L}(l^2(V(\Gamma)))$ его резольвентным
множеством $r(A)$ называется множество всех таких $\lambda \in \mathbb{C}$, что
$A - \lambda E$ есть биекция $l_2(V(\Gamma))$ на $l_2(V(\Gamma))$ (и тогда
$(A - \lambda E)^{-1} \in {\cal L}(l^2(V(\Gamma)))$ в силу теоремы Банаха об обратном операторе),
определенный для $\lambda \in r(A)$ оператор $(\lambda E - A)^{-1}$ называется резольвентой $A$ в $\lambda$,
а спектр $\sigma(A)$ есть $\mathbb{C} \setminus r(A)$.

Содержащееся в $\mathbb{C}^{V(\Gamma)}$ в качестве подпространства гильбертово пространство $l_2(V(\Gamma))$
(топология которого не индуцируется
топологией топологического векторного пространства
$\mathbb{C}^{V(\Gamma)}$) является $A^{({\rm alg})}_{\Gamma,\mathbb{C}}$-инвариант\-ным,
а ограничение $A^{({\rm alg})}_{\Gamma,\mathbb{C}}$ на $l_2(V(\Gamma))$, обозначаемое далее через $A^{(l_2)}_{\Gamma,\mathbb{C}}$,
есть оператор из ${\cal L}(l^2(V(\Gamma)))$.

\begin{remark}
\label{r6.6}

Как уже было сказано во введении, в случае связного графа $\Gamma$, степени вершин
которого ограничены в совокупности, оператор $A^{(l_2)}_{\Gamma,\mathbb{C}}$
совпадает с оператором, рассматриваемым в~\cite{BM},~\cite{MW} и называемым там оператором смежности.
В~\cite{BM},~\cite{MW} оператор, называемый там оператором смежности, определяется также в более общем случае счетного локально конечного графа $\Gamma$
как оператор в $l_2(V(\Gamma))$, являющийся замыканием
ограничения $A^{({\rm alg})}_{\Gamma,\mathbb{C}}$
на пространство функций с конечными носителями.
\end{remark}

Опрератор $A^{(l_2)}_{\Gamma,\mathbb{C}}$ обладает следующими несложно устанавливаемыми свойствами (ср.~\cite{BM},~\cite{MW}):

a) Норма $A^{(l_2)}_{\Gamma,\mathbb{C}}$ не превосходит $d$ (см.,
например,~\cite[теорема 6.12-А]{Tay}).

b) $A^{(l_2)}_{\Gamma,\mathbb{C}}$ самосопряжен (см., например,~\cite[с. 329--330]{Tay}).

c) $\sigma(A^{(l_2)}_{\Gamma,\mathbb{C}})$ содержится в
$\mathbb{R}_{\geq - d} \cap \mathbb{R}_{\leq d}$ (что следует из а), b) и, например,~\cite[теорема 2.2.5]{Br}).

Оператор  $A^{(l_2)}_{\Gamma,\mathbb{C}}$ используется нами для доказательства следующей
теоремы, первое утверждение которой непосредственно вытекает из указанного свойства c).

\begin{theorem}
\label{t6.2}

Пусть $\Gamma$ --- связный граф, степени вершин которого не превосходят некоторого натурального
числа $d$. Тогда справедливы следующие утверждения$:$

$1)$ Множество $\mathbb{C} \setminus (\mathbb{R}_{\geq - d} \cap \mathbb{R}_{\leq d})$
содержится в резольвентном множестве оператора $A^{(l_2)}_{\Gamma,\mathbb{C}}$. Таким образом,
для каждого $\lambda \in \mathbb{C} \setminus (\mathbb{R}_{\geq - d} \cap \mathbb{R}_{\leq d})$
определен линейный оператор
$R_{\Gamma,d,\lambda} := (A^{(l_2)}_{\Gamma,\mathbb{C}} - \lambda E)^{-1} \in {\cal L}(l^2(V(\Gamma)))$.

$2)$ Для произвольных $v \in V(\Gamma)$ и $\lambda \in \mathbb{C} \setminus (\mathbb{R}_{\geq - d} \cap \mathbb{R}_{\leq d})$
функция $R_{\Gamma,d,\lambda}(\delta_v) \in \mathbb{C}^{V(\Gamma)}$
является $(\mathbb{C},\lambda)$-пропагатором графа $\Gamma$ относительно вершины $v$,
причем содержащимся в $l_2(V(\Gamma))$ и единственным с этим свойством.

$3)$ Для произвольных $v, w \in V(\Gamma)$ и
$\lambda \in \mathbb{C} \setminus (\mathbb{R}_{\geq - d} \cap \mathbb{R}_{\leq d})$
положим $$\tilde s_{\Gamma,d,\lambda,v}(w) := (R_{\Gamma,d,\lambda}(\delta_v))(w) = (R_{\Gamma,d,\lambda}(\delta_v),\delta_w)$$
$($склярное произведение $R_{\Gamma,d,\lambda}(\delta_v)$ и $\delta_w$ в $l_2(V(\Gamma))$$)$.
Тогда для произвольных фиксированных $v, w \in V(\Gamma)$ определенная на
$\mathbb{C} \setminus (\mathbb{R}_{\geq - d} \cap \mathbb{R}_{\leq d})$
функция $\tilde s_{\Gamma,d,\lambda,v}(w)$ аргумента $\lambda$
является аналитической функцией, совпадающей на множестве $\{\lambda \in \mathbb{C} : |\lambda| > d\}$
с $s_{\Gamma,\lambda,v}(w)$, рассматриваемой как функция аргумента $\lambda$.

\end{theorem}

\begin{proof}
Из утверждения 1) теоремы (вытекающего из свойства с) оператора
$A^{(l_2)}_{\Gamma,\mathbb{C}}$) следует, что
для произвольных $v \in V(\Gamma)$ и $\lambda \in \mathbb{C} \setminus (\mathbb{R}_{\geq - d} \cap \mathbb{R}_{\leq d})$
функция $R_{\Gamma,d,\lambda}(\delta_v) \in \mathbb{C}^{V(\Gamma)}$
является $(\mathbb{C},\lambda)$-пропагатором графа $\Gamma$ относительно вершины $v$,
содержащимся в $l_2(V(\Gamma))$, причем единственным
$(\mathbb{C},\lambda)$-про\-пагатором графа $\Gamma$ относительно вершины $v$,
содержащимся в $l_2(V(\Gamma))$, поскольку разность различных $(\mathbb{C},\lambda)$-пропагаторов графа $\Gamma$ относительно вершины $v$,
содержащихся в $l_2(V(\Gamma))$, была бы собственной функцией оператора $A^{(l_2)}_{\Gamma,\mathbb{C}}$,
соответствующей собственному значению $\lambda$, а это невозможно в силу $\lambda \not \in \sigma(A^{(l_2)}_{\Gamma,\mathbb{C}})$) (см. с)).

Остается доказать утверждение 3) теоремы.
Его справедливость в части аналитичности функции
$\tilde s_{\Gamma,d,\lambda,v}(w)$ аргумента $\lambda \in \mathbb{C} \setminus (\mathbb{R}_{\geq - d} \cap \mathbb{R}_{\leq d})$ следует, например, из~\cite[теорема 4.12]{Stone}.
Покажем, наконец,
что для произвольных $v, w \in V(\Gamma)$ и
$\lambda \in \mathbb{C}$ с $|\lambda| > d$ имеет место равенство
$\tilde s_{\Gamma,d,\lambda,v}(w) = s_{\Gamma,\lambda,v}(w)$.
При $|\lambda| > d$ ввиду свойства a) оператора $A^{(l_2)}_{\Gamma,\mathbb{C}}$
норма оператора $\lambda^{-1} A^{(l_2)}_{\Gamma,\mathbb{C}}$ меньше $1$,
и следовательно, при $|\lambda| > d$ имеем (в смысле сходимости суммы по равномерной операторной топологии
на ${\cal L}(l^2(V(\Gamma)))$)
$$R_{\Gamma,d,\lambda} = (A^{(l_2)}_{\Gamma,\mathbb{C}} - \lambda E)^{-1}
= -\lambda^{-1}(E - \lambda^{-1} A^{(l_2)}_{\Gamma,\mathbb{C}})^{-1}$$
$$= -\lambda^{-1}(E + \sum_{n \in \mathbb{Z}_{\geq 1}} (\lambda^{-1} A^{(l_2)}_{\Gamma,\mathbb{C}})^n),$$
что для произвольных $v, w \in V(\Gamma)$ влечет
$$\tilde s_{\Gamma,d,\lambda,v}(w)
= (R_{\Gamma,d,\lambda}(\delta_v),\delta_w) =
-\lambda^{-1}((\delta_v,\delta_w) +
(\sum_{n \in \mathbb{Z}_{\geq 1}} \lambda^{-n}(A^{(l_2)}_{\Gamma,\mathbb{C}})^n(\delta_v),\delta_w))$$
$$= -\lambda^{-1}((\delta_v,\delta_w) +
\sum_{n \in \mathbb{Z}_{\geq 1}}(\lambda^{-n}(A^{(l_2)}_{\Gamma,\mathbb{C}})^n(\delta_v),\delta_w))
= s_{\Gamma,\lambda,v}(w).$$

\end{proof}

\begin{remark}
\label{r6.7}

Если выполнены условия теоремы~\ref{t6.2}, то согласно утверждениям $2)$ и $3)$ этой теоремы
для произвольных фиксированных $v \in V(\Gamma)$ и $\lambda \in \{c \in \mathbb{C} : |c| > d\}$
содержащаяся в $\mathbb{C}^{V(\Gamma)}$ функция $s_{\Gamma,\lambda,v}(w)$ аргумента $w \in V(\Gamma)$
лежит в $l_2(V(\Gamma))$.

\end{remark}

\begin{remark}
\label{r6.8}

Подчеркнем, что из утверждения $2)$ теоремы~\ref{t6.2} следует {\it существование
и единственность} для произвольного связного графа $\Gamma$, степени вершин которого ограничены в совокупности, и произвольных $v \in V(\Gamma)$, $\lambda \in \mathbb{C} \setminus (\mathbb{R}_{\geq - d} \cap \mathbb{R}_{\leq d})$,
где $d$ --- максимум степеней вершин графа $\Gamma$,
квадратично суммируемого $(\mathbb{C},\lambda)$-пропагатора $\tilde s_{\Gamma,d,\lambda,v}$
графа $\Gamma$ относительно вершины $v$ (своего рода {\it канонического
$(\mathbb{C},\lambda)$-пропагатора $\Gamma$ относительно} $v$).
(Для указанных $\Gamma$  и $\lambda$
у $A_{\Gamma,\mathbb{C}}$ отсутствуют, как следствие,
собственные квадратично суммируемые функции, соответствующие
собственному значению $\lambda$.) Хотя в части существования пропагаторов это утверждение, даже в случае
связного графа $\Gamma$ с не превосходящими некоторого натурального числа $d$ степенями вершин и $F = \mathbb{C}$ c обычным
абсолютным значением, слабее утверждения 2) следствия~\ref{c4.2}, оно при $\lambda \in  \mathbb{C} \setminus (\mathbb{R}_{\geq - d} \cap \mathbb{R}_{\leq d})$ имеет то преимущество, что
доставляет ``явный вид'' пропагатора относительно произвольной заданной вершины, причем единственного содержащегося
в  $l_2(V(\Gamma))$. Заманчивым выглядит изучение аналитических свойств
этих канонических пропагаторов $\tilde s_{\Gamma,d,\lambda,v}$ (как вектор-функций аргумента
$\lambda \in  \mathbb{C} \setminus (\mathbb{R}_{\geq - d} \cap \mathbb{R}_{\leq d})$)
с использованием, например, соотношения Гильберта для $R_{\Gamma,d,\lambda}$ и равенства $R_{\Gamma,d,\lambda}^2$ производной $R_{\Gamma,d,\lambda}$ по $\lambda$.
Для заданного $w \in V(\Gamma)$
в связи с предыдущими результатами работы значительный интерес представляют нули
$\tilde s_{\Gamma,d,\lambda,v}(w)$, рассматриваемой как функции аргумента
$\lambda \in  \mathbb{C} \setminus (\mathbb{R}_{\geq - d} \cap \mathbb{R}_{\leq d})$.

\end{remark}

\begin{remark}
\label{r6.9}

Обратим внимание еще на один аспект утверждения $3)$ теоремы~\ref{t6.2}. Фиксируя $d \in \mathbb{Z}_{\geq 2}$
и меняя связные графы $\Gamma$, степени вершин которых не превосходят $d$,  а также вершины $v, w \in V(\Gamma)$,
мы получаем широкий класс голоморфных на  $\{c \in \mathbb{C}: |c| > d\}$ функций
$s_{\Gamma,\lambda,v}(w)$ аргумента $\lambda$, допускающих аналитическое продолжение до
голоморфных на  $\mathbb{C} \setminus  (\mathbb{R}_{\geq - d} \cap \mathbb{R}_{\leq d})$
функций $\tilde s_{\Gamma,d,\lambda,v}(w)$ аргумента $\lambda$.
На эту ситуацию можно смотреть и со следующих позиций: для связного графа $\Gamma$, степени вершин которого не превосходят $d$,
и для произвольных $v, w \in V(\Gamma)$ производящий ряд
$$W_{\Gamma,\mathbb{C}}(x,v,w) = \sum_{n \in \mathbb{Z}_{\geq 0}} N_{\Gamma}(v,w,n) x^n,$$
абсолютно сходящийся при $x \in \{c \in \mathbb{C} : |c| < 1/d\}$, можно естественным образом
``просуммировать'' для любого
$x \in \{c \in \mathbb{C} : |c| \geq 1/d\} \setminus (\mathbb{R}_{\leq - 1/d} \cup \mathbb{R}_{\geq 1/d})$,
приписав ему значение $-x^{-1} \tilde s_{\Gamma,d,x^{-1},v}(w)$.

\end{remark}

\begin{remark}
\label{r6.10}

В связи с утверждением $3)$ теоремы~\ref{t6.2} обратим внимание на \S~\ref{s8}, раздел~\ref{s8.3},
где строится бесконечный кубический связный граф $\Gamma$ такой, что для комплексного невещественного корня $\lambda_0$
многочлена $x^3 + x^2 - 1 \in \mathbb{C}[x]$ (любого из двух) оператор $A_{\Gamma,\mathbb{C}} - \lambda_0 E$ инъективен, т. е. согласно теореме~\ref{t5.3} относительно
каждой вершины $v$ у $\Gamma$ имеется
единственный $(\mathbb{C},\lambda_0)$-пропагатор $p_v$, причем $p_v$ имеет конечный носитель. Но тогда
согласно утверждениям $2)$ и $3)$ теоремы~\ref{t6.2} для произвольных $v, w \in V(\Gamma)$ имеем
$p_v(w) = \tilde s_{\Gamma,3,\lambda_0,v}(w)$. Таким образом, для произвольной вершины $v \in V(\Gamma)$ для
всех, за исключением конечного числа, вершин $w \in V(\Gamma)$ аналитическая на
$\mathbb{C} \setminus (\mathbb{R}_{\geq - 3} \cap \mathbb{R}_{\leq 3})$ согласно
утверждению $3)$ теоремы~\ref{t6.2}
функция $\tilde s_{\Gamma,3,\lambda,v}(w)$ аргумента $\lambda$ обращается в $0$ в точке $\lambda_0$.
При этом согласно утверждению 2) предложения~\ref{p6.3} для произвольных $v, w \in V(\Gamma)$ аналитическая на
$\mathbb{C} \setminus (\mathbb{R}_{\geq - 3} \cap \mathbb{R}_{\leq 3})$
функция $\tilde s_{\Gamma,3,\lambda,v}(w)$ аргумента $\lambda$ всюду на $\mathbb{R}_{\geq 3}$ принимает
вещественные отрицательные значения.

\end{remark}

\section{``Типичность'' сюръективности и ``исключительность'' (в случае бесконечного связного
графа $\Gamma$ и поля $F$ нулевой характеристики)
инъективности для операторов $A_{\Gamma,F} - \lambda E$}
\label{s7}

Всюду в этом параграфе $\Gamma$ --- {\it локально конечный граф и} $F$ --- {\it поле с некоторым абсолютным значением} $|.|_{\rm v}$.
Напомним (см. теоремы~\ref{t5.1} и~\ref{t5.3}), что для $\lambda \in F$ условие сюръективности оператора $A_{\Gamma, F} - \lambda E$
эквивалентно условию наличия $(F,\lambda)$-пропага\-торов графа $\Gamma$ относительно всех его вершин (т. е. условию $S_{F,\lambda}(\Gamma) = \emptyset$), а
условие инъективности оператора $A_{\Gamma, F} - \lambda E$
эквивалентно условию наличия $(F,\lambda)$-про\-пагаторов с конечными носителями графа $\Gamma$ относительно всех его вершин
(т. е. условию $L_{F,\lambda}(\Gamma) = V(\Gamma)$).

Ясно, что в случае конечного графа $\Gamma$ условия сюръективности и инъективности оператора $A_{\Gamma, F} - \lambda E$
равносильны и ``типичны'' в том смысле, что не выполняются лишь для $\lambda$, являющихся
собственными значениями матрицы смежности графа $\Gamma$ над $F$.
Здесь и далее, называя некоторое свойство элементов $\lambda$ поля $F$ ``типичным'' или, напротив,
``исключительным'', мы имеем в виду (проецируя ситуацию на случай,
когда $F$ есть $\mathbb{C}$
или, более общо, есть поле, не являющееся алгебраическим расширением своего простого подполя),
что свойство может не выполняться или, соответственно, выполняться лишь для определенных алгебраических над простым подполем элементов $\lambda$.

Однако, как будет показано в этом
параграфе, из предыдущих результатов настоящей работы следует, что ситуация меняется в случае бесконечного
графа $\Gamma$. Хотя для $\lambda \in F$ условие наличия $(F,\lambda)$-пропагаторов бесконечного графа $\Gamma$ относительно всех его вершин
(т. е. условие сюръективности $A_{\Gamma, F} - \lambda E$) согласно следствию~\ref{c4.2} по-прежнему является  ``типичным'', более сильное условие
наличия $(F,\lambda)$-пропагаторов с конечными носителями бесконечного связного графа $\Gamma$ относительно всех его вершин
(т. е. условие  инъективности $A_{\Gamma, F} - \lambda E$) является ``исключительным'' для $\lambda$ в случае
поля $F$
характеристики $0$
(но может не быть  ``исключительным'' для $\lambda$ в случае поля $F$ положительной характеристики, см. \S~\ref{s8}, раздел~\ref{s8.1.1}).
Как будет показано ниже (см. предложения~\ref{p7.3} и~\ref{p7.2}), в случае бесконечного связного
графа $\Gamma$ и поля $F$ характеристики $0$ ``исключительным'' для $\lambda \in F$ является даже условие
наличия $(F,\lambda)$-пропагатора с конечным носителем графа $\Gamma$ относительно {\it какой-либо} вершины
(т. е. условие $L_{F,\lambda}(\Gamma) \not = \emptyset$).
Таким образом, будет показано, что для бесконечного локально конечного связного графа $\Gamma$ и поля $F$ характеристики $0$ ``типичной'' для $\lambda \in F$ является ситуация, когда $S_{F,\lambda}(\Gamma) = \emptyset$,
$L_{F,\lambda}(\Gamma) = \emptyset$ и (см. теоремы~\ref{t4.1} и~\ref{t4.2}) носители всех собственных функций
оператора $A_{\Gamma,F}$, соответствующих собственному значению $\lambda$, бесконечны и в объединении
дают все множество $V(\Gamma)$.
Отметим, однако, что имеются такие бесконечные локально конечные связные графы $\Gamma$ и поля $F$
характеристики $0$ с абсолютными значениями, что для некоторых $\lambda \in F$
имеет место равенство $S_{F,\lambda}(\Gamma) = V(\Gamma)$ (см. \S~\ref{s8}, раздел~\ref{s8.5}), и
что имеются такие графы $\Gamma$ и поля $F$ с теми же свойствами, что для некоторых $\lambda \in F$
имеет место равенство $L_{F,\lambda}(\Gamma) = V(\Gamma)$ или, другими словами, $\lambda$ не является
собственным значением оператора $A_{\Gamma, F}$ (см. \S~\ref{s8}, разделы~\ref{s8.2},~\ref{s8.3} и~\ref{s8.4}).

Иная ситуация вполне может иметь место в случае поля $F$ положительной характеристики:
для каждого простого $p \in \mathbb{Z}_{\geq 1}$
несложно привести пример бесконечного локально конечного связного вершинно-симметрического графа $\Gamma$ такого,
что для любого поля $F$ характеристики $p$ с произвольным абсолютным значением и любого
$\lambda \in F \setminus \{0\}$
имеет место равенство $L_{F,\lambda}(\Gamma) = V(\Gamma)$ или, другими словами, $\lambda$ не является
собственным значением оператора $A_{\Gamma, F}$ (см. \S~\ref{s8}, раздел~\ref{s8.1.1}).

Итак, из следствия~\ref{c4.2} теоремы~\ref{t4.1} вытекает,  что для локально конечного графа $\Gamma$ и поля $F$ с абсолютным значением $|.|_{\rm v}$ ``типичной'' для $\lambda \in F$ является ситуация, когда у графа $\Gamma$ относительно каждой его вершины
имеется $(F,\lambda)$-пропагатор, или, другими словами, для $\lambda \in F$ ``исключительным'' является условие отсутствия
у графа $\Gamma$ относительно какой-либо его вершины
$(F,\lambda)$-пропагатора (т. е. условие $S_{F,\lambda}(\Gamma) \not = \emptyset$). Вернемся к выводу из теоремы~\ref{t4.1} ``исключительности'' для $\lambda \in F$ условия $S_{F,\lambda}(\Gamma) \not = \emptyset$,
дополнив его одним уточняющим наблюдением.
Итак, согласно теореме~\ref{t4.1} для локально конечного графа $\Gamma$, поля $F$ и $\lambda \in F$ условие
$S_{F,\lambda}(\Gamma) \not = \emptyset$ эквивалентно существованию
у оператора $A_{\Gamma, F}$ собственной функции с конечным носителем, соответствующей
собственному значению $\lambda$. Отсюда следует ``исключительность'' для $\lambda \in F$ этого свойства,
поскольку, если $X$ --- конечный носитель некоторой собственной функции
оператора $A_{\Gamma, F}$, соответствующей собственному значению $\lambda$,
то для любого содержащего $X$ конечного подмножества $X'$ множества $V(\Gamma)$ элемент
$\lambda$ {\it является
собственным значением матрицы смежности над $F$ индуцированного $X'$ конечного
подграфа $\langle X' \rangle_{\Gamma}$
графа $\Gamma$.} Добавим, что, как это следует из очевидного предложения~\ref{p7.1}, условие существования
у оператора $A_{\Gamma, F}$ собственной функции с конечным носителем, соответствующей
собственному значению $\lambda$, вообще говоря, более ограничительно, чем условие принадлежности
$\lambda$ множеству собственных значений матрицы смежности над $F$ некоторого индуцированного конечного
подграфа графа $\Gamma$.

\begin{propos}
\label{p7.1}

Пусть $\Gamma$ --- локально конечный граф, $F$ --- поле с некоторым абсолютным значением $|.|_{\rm v}$ и $\lambda \in F$.
Тогда для $\emptyset \not = X \subseteq V(\Gamma)$
следующие условия равносильны$:$

$(i)$ У оператора $A_{\Gamma, F}$ существует собственная функция, соответствующая
собственному значению $\lambda$, носитель которой содержится в $X$.

$(ii)$ У оператора $A_{\langle X \rangle_{\Gamma}, F}$ существует такая собственная
функция $f$, соответствующая
собственному значению $\lambda$, что для любой вершины
$u \in \Gamma (X) \setminus X$ $($или, эквивалентно, для любой вершины $u \in V(\Gamma) \setminus X$$)$
справедливо равенство $$\sum_{u' \in \Gamma(u)\cap X}f(u') = 0.$$

\end{propos}

\begin{proof}[{\rm очевидно}] Если имеет место $(i)$ и $\tilde f$ --- собственная функция
оператора $A_{\Gamma, F}$, соответствующая
собственному значению $\lambda$, носитель которой содержится в $X$, то $(ii)$ имеет место, поскольку очевидно, что в качестве
$f$ можно взять $\tilde f|_X$. Обратно, из справедливости $(ii)$ следует справедливость
$(i)$, поскольку $\tilde f \in F^{V(\Gamma)}$, совпадающая с $f$ на $X$ и равная $0$ всюду
на $V(\Gamma) \setminus X$, очевидным образом является собственной функцией
оператора $A_{\Gamma, F}$, соответствующей
собственному значению $\lambda$, носитель которой содержится в $X$.

\end{proof}

\begin{remark}
\label{r7.1}
Для пропагаторов справедлив следующий аналогичный результат (с аналогичным доказательством).
В предположениях
предложения~\ref{p7.1} пусть $v \in X$. Тогда
следующие условия равносильны:

$(i)$ У графа $\Gamma$ имеется $(F,\lambda)$-пропагатор относительно вершины $v$, носитель
которого содержится в $X$.

$(ii)$ У графа $\langle X \rangle_{\Gamma}$ имеется такой $(F,\lambda)$-пропагатор $p_v$
относительно вершины $v$, что для любой вершины
$u \in \Gamma (X) \setminus X$ $($или, эквивалентно, для любой вершины $u \in V(\Gamma) \setminus X$$)$
справедливо равенство $\sum_{u' \in \Gamma(u)\cap X}p_v(u') = 0.$

\end{remark}

Мы докажем теперь, что если $\Gamma$ --- {\it бесконечный} локально конечный {\it связный} граф и $F$ --- поле {\it характеристики} $0$
с некоторым абсолютным значением $|.|_{\rm v}$,
то для $\lambda \in F$, в отличие от случая конечного графа, ``исключительным'' является уже наличие $(F,\lambda)$-пропагатора
с {\it конечным} носителем
графа $\Gamma$ относительно {\it какой-либо} его вершины, что согласно теореме~\ref{t5.3} влечет
также ``исключительность'' инъективности
оператора $A_{\Gamma, F} - \lambda E$. Впрочем, ``исключительность'' инъективности
оператора $A_{\Gamma, F} - \lambda E$ в том смысле, что она возможна лишь для алгебраических
над $\mathbb{Q} \subseteq F$ элементов $\lambda$, вытекает из теоремы~\ref{t6.1}.

\medskip

Предварительно докажем следующее

\begin{propos}
\label{p7.2}
Пусть $F$ --- поле характеристики $0$ и $F[x]$ --- алгебра многочленов от независимой переменной $x$ с коэффициентами
из $F$. Пусть $\Delta$ --- конечный связный граф с $|V(\Delta)| > 1$ и ${\bf A}_{\Delta,F} - x {\bf E}$ --- характеристическая матрица его матрицы смежности
$($над $F$$)$. Тогда произвольный минор порядка $|V(\Delta)| - 1$ матрицы ${\bf A}_{\Delta,F} - x {\bf E}$,
рассматриваемый как элемент $F[x]$, отличен от нуля $($и принадлежит $\mathbb{Z}[x] \subseteq F[x]$$).$
Как следствие, матрица ${\bf A}_{\Delta,F} - x {\bf E}$,
рассматриваемая над полем $\mathbb{Q}(x)$, имеет обратную, все элементы которой ненулевые.

\end{propos}

\begin{proof}
Не теряя общности, будем считать, что $F = \mathbb{Q}$.
Предположим, что нулевым элементом $\mathbb{Q}[x]$ является определитель матрицы ${\bf A}'$, получающейся из
${\bf A}_{\Delta,\mathbb{Q}} - x {\bf E}$
удалением строки, соответствующей вершине $u$ графа $\Delta$, и столбца, соответствующего вершине $u'$ графа $\Delta$.
Для каждого $i \in \{1,...,|V(\Delta)|\}$ пусть $v_i$ --- вершина графа $\Delta$, соответствующая $i$-ой строке
(и $i$-ому столбцу) матрицы ${\bf A}_{\Delta,\mathbb{Q}}$.
Поскольку $u \not = u'$ (ясно, что все главные миноры матрицы ${\bf A}_{\Delta,\mathbb{Q}} - x {\bf E}$ ненулевые),
можно, не теряя общности, предполагать, что $u = v_1$ и $u' = v_{|V(\Delta)|}$.
Пусть ${\boldsymbol \delta}_1 := (1,0,...,0)^t$ (матрица размера $|V(\Delta)| \times 1$ над $\mathbb{Q}[x]$).
Тогда в силу $\det ({\bf A}_{\Delta,\mathbb{Q}} - x {\bf E}) \not = 0$ и $\det({\bf A}') = 0$
система линейных уравнений
$$({\bf A}_{\Delta,\mathbb{Q}} - x {\bf E})(x_1,...,x_{|V(\Delta)|})^t = {\boldsymbol \delta}_1$$
над полем $\mathbb{Q}(x)$ имеет единственное решение $(x'_1,...,x'_{|V(\Delta)|})^t$,
причем (по правилу Крамера) $x'_{|V(\Delta)|} = 0$. Но тогда для любого $\lambda \in \mathbb{R}$ со свойством
$\det ({\bf A}_{\Delta,\mathbb{Q}} - \lambda {\bf E}) \not = 0$  (матрица ${\bf A}_{\Delta,\mathbb{Q}} - \lambda {\bf E}$
рассматривается как матрица над $\mathbb{R}$) система линейных уравнений
\begin{equation}
\label{eq7.1}
({\bf A}_{\Delta,\mathbb{Q}} - \lambda {\bf E})(x_1,...,x_{|V(\Delta)|})^t = {\boldsymbol \delta}_1,
\end{equation}
рассматриваемая как система линейных уравнений над полем $\mathbb{R}$, также имеет единственное решение $(x'_{\lambda,1},...,x'_{\lambda,|V(\Delta)|})^t$,
причем $x'_{\lambda,|V(\Delta)|} = 0$.
Пусть $d$ --- максимум степеней вершин графа $\Delta$ и $\lambda_0 \in  \mathbb{R}_{> d}$.
Тогда  $\det ({\bf A}_{\Delta,\mathbb{Q}} - \lambda_0 {\bf E}) \not = 0$. Кроме того,
согласно утверждению $1)$ предложения~\ref{p6.3} решением~\eqref{eq7.1}, рассматриваемой как система линейных уравнений
над $\mathbb{R}$, при $\lambda = \lambda_0$ является $(s_1,...,s_{|V(\Delta)|})^t$, где $s_i := s_{\Delta,\lambda_0,v_1}(v_i)$ для каждого $i \in \{1,...,|V(\Delta)|\}$.
Следовательно, $s_i = x'_{\lambda_0,i}$ для всех $i \in \{1,...,|V(\Delta)|\}$.
Однако $s_{|V(\Delta)|} = s_{\Delta,\lambda_0,v_1}(v_{|V(\Delta)|}) < 0$ в силу
утверждения $2)$ предложения~\ref{p6.3},
а $x'_{\lambda_0,|V(\Delta)|}=0$. Полученное противоречие доказывает предложение~\ref{p7.2}.

\end{proof}

\begin{remark}
\label{r7.2}
Для любого поля $F$ характеристики $p > 0$ аналог предложения~\ref{p7.2}, вообще говоря, не имеет места,
как показывает пример подграфа, порожденного шаром радиуса $2$, графа из раздела~\ref{s8.1.1} параграфа~\ref{s8} (с тем же самым $p$).

\end{remark}

\begin{remark}
\label{r7.3}
При дополнительном предположении, что $\Delta$ регулярен, заключение предложения~\ref{p7.2}
является следствием матричной теоремы Кирхгофа о деревьях.

\end{remark}

\begin{cor}
\label{c7.1}

Пусть $\Delta$ --- конечный связный граф с $|V(\Delta)| > 1$.
Пусть $F$ --- поле характеристики $0$ и $\lambda$ --- элемент $F$, не являющийся собственным значением
матрицы смежности графа $\Delta$ $($над $F$$)$, но такой, что $(F,\lambda)$-пропага\-тор графа $\Delta$
относительно некоторой вершины графа $\Delta$ принимает нулевое значение в какой-то вершине
графа $\Delta$. Тогда $\lambda$ --- алгебраический над $\mathbb{Q} \subseteq F$ элемент степени $< |V(\Delta)|$.

\end{cor}
\begin{proof}
Поскольку $\lambda$ не является собственным значением
матрицы смежности графа $\Delta$, условие, что $(F,\lambda)$-пропагатор графа $\Delta$ относительно некоторой
вершины $u$ принимает нулевое значение в какой-то вершине $u'$, эквивалентно равенству нулю
определителя матрицы, получающейся из
${\bf A}_{\Delta,F} - \lambda {\bf E}$
удалением строки, соответствующей вершине $u$ графа $\Delta$, и столбца, соответствующего вершине
$u'$ графа $\Delta$.
С учетом этого требуемое утверждение вытекает из предложения~\ref{p7.2}.

\end{proof}

Теперь мы можем аргументированно повторить сказанное ранее.
Пусть $\Gamma$ --- бесконечный локально конечный связный граф и $F$ --- поле характеристики $0$
с некоторым абсолютным значением $|.|_{\rm v}$. Из приводимого ниже очевидного предложения~\ref{p7.3}
и предложения~\ref{p7.2} следует ``исключительность'' для $\lambda \in F$ условия наличия у графа $\Gamma$
относительно {\it некоторой} его вершины $(F,\lambda)$-пропагатора с {\it конечным} носителем
(т. е. условия $L_{F,\lambda}(\Gamma) \not = \emptyset$).
Еще более ``исключительным'' является для $\lambda \in F$ условие
инъективности
оператора $A_{\Gamma, F} - \lambda E$, которое согласно теореме~\ref{t5.3} равносильно условию наличия у графа $\Gamma$
относительно {\it каждой} его вершины
$(F,\lambda)$-пропагатора с {\it конечным} носителем (т. е. условию $L_{F,\lambda}(\Gamma) = V(\Gamma)$).

\begin{propos}
\label{p7.3}
Пусть $\Gamma$ --- локально конечный граф, $F$ --- поле с некоторым абсолютным значением $|.|_{\rm v}$ и $v \in V(\Gamma)$.
Тогда для $\lambda \in F$ из
существования $(F,\lambda)$-пропагатора с конечным носителем у $\Gamma$ относительно $v$
$($т. е. из $v \in L_{F,\lambda}(\Gamma)$$)$
следует существование такого содержащего $v$ конечного подмножества $X$ множества $V(\Gamma)$,
что для любого содержащего $X$ конечного подмножества $X'$ множества $V(\Gamma)$
равен нулю определитель матрицы ${\bf A}_{\langle X'\rangle_{\Gamma},F} - \lambda {\bf E}$
или равны нулю все миноры этой матрицы,
получающиеся из нее удалением строки, соответствующей вершине $v$, и столбца, соответствующего
любой из вершин множества $X' \setminus X$.

\end{propos}

\section{Некоторые примеры}
\label{s8}

\subsection{}
\label{s8.1}
Пусть $F$ --- произвольное {\it поле характеристики} $p > 0$ с произвольным абсолютным значением.

\subsubsection{}
\label{s8.1.1}
Несложно привести {\it пример бесконечного локально конечного связного
вершинно-симметрического графа $\Gamma$ такого, что никакое $\lambda \in F \setminus \{0\}$ не является собственным значением
оператора смежности $A_{\Gamma,F}$}.

Действительно, пусть $\Delta$ --- произвольный бесконечный локально конечный связный
вершинно-симметрический граф. Определим граф $\Gamma$, полагая
$$V(\Gamma) = \{ v_{u,k} : u \in V(\Delta), k \in \{1,...,p\}\},$$
$$E(\Gamma) = \{\{v_{u',k'},v_{u'',k''}\} : \{u',u''\} \in E(\Delta), k',k'' \in \{1,...,p\}\}.$$
Покажем, что никакое $\lambda \in F \setminus \{0\}$ не является собственным значением
оператора смежности $A_{\Gamma,F}$.

Пусть $\lambda \in F \setminus \{0\}$, и пусть $f \in F^{V(\Gamma)}$ такова, что $$\lambda f(v) = \sum_{w \in \Gamma(v)}f(w)$$
для всех $v \in V(\Gamma)$. Тогда для произвольных $u' \in V(\Delta), k_1,k_2 \in \{1,...,p\}$
в силу $\Gamma(v_{u',k_1}) = \Gamma(v_{u',k_2})$ имеем $\lambda f(v_{u',k_1}) = \lambda f(v_{u',k_2})$,
и потому (с учетом $\lambda \not = 0$) $f(v_{u',k_1}) = f(v_{u',k_2})$.
Следовательно, для произвольной вершины $v_{u,k}$, $u \in V(\Delta), k \in \{1,...,p\}$, графа $\Gamma$
имеем
$$f(v_{u,k}) = \frac{1}{\lambda} \sum_{u' \in \Delta(u)} \sum_{k' \in \{1,...,p\}}f(v_{u',k'}) =
\frac{1}{\lambda} \sum_{u' \in \Delta(u)} p\cdot f(v_{u',1}) = 0.$$
Таким образом, $f$ --- функция, тождественно равная 0 на $V(\Gamma)$. Как следствие, $\lambda$ не является собственным значением
оператора смежности $A_{\Gamma,F}$.

\subsubsection{}
\label{s8.1.2}
Несложно привести также {\it пример бесконечного локально конечного связного
вершинно-симметрического графа $\Gamma$ такого, что $0 \in F$ не является собственным значением
оператора смежности $A_{\Gamma,F}$}.

Определим граф $\Gamma$, полагая
$$V(\Gamma) = \{ v_{i,k} : i \in \mathbb{Z}, k \in \{1,...,p\}\},$$
$$E(\Gamma) = \{\{v_{2i',k'},v_{2i'+1,k''}\} : i' \in \mathbb{Z}, k',k'' \in \{1,...,p\}\} \cup$$
$$\cup \{\{v_{2i',k'},v_{2i'-1,k'}\} : i' \in \mathbb{Z}, k' \in \{1,...,p\}\}.$$
Покажем, что $0 \in F$ не является собственным значением
оператора смежности $A_{\Gamma,F}$.

Пусть $f \in F^{V(\Gamma)}$ такова, что $$(0 \cdot f(v) =)\ 0 = \sum_{w \in \Gamma(v)}f(w)$$
для всех $v \in V(\Gamma)$. Покажем, прежде всего, что для произвольных $i \in \mathbb{Z}$ и $k_1,k_2 \in \{1,...,p\}$
справедливо равенство $f(v_{i,k_1}) = f(v_{i,k_2})$. Действительно, при четном $i$ имеем
$$0 = \sum_{w \in \Gamma(v_{i-1,k_1})}f(w) = f(v_{i,k_1}) + \sum_{k \in \{1,...,p\}}f(v_{i-2,k})$$
и
$$0 = \sum_{w \in \Gamma(v_{i-1,k_2})}f(w) = f(v_{i,k_2}) + \sum_{k \in \{1,...,p\}}f(v_{i-2,k}),$$
что влечет $f(v_{i,k_1}) = f(v_{i,k_2})$, а при нечетном $i$ имеем
$$0 = \sum_{w \in \Gamma(v_{i+1,k_1})}f(w) = f(v_{i,k_1}) + \sum_{k \in \{1,...,p\}}f(v_{i+2,k})$$
и
$$0 = \sum_{w \in \Gamma(v_{i+1,k_2})}f(w) = f(v_{i,k_2}) + \sum_{k \in \{1,...,p\}}f(v_{i+2,k}),$$
что влечет $f(v_{i,k_1}) = f(v_{i,k_2})$ и в этом случае.
Теперь для произвольного четного $i \in \mathbb{Z}$ и произвольного $k \in \{1,...,p\}$ имеем
$$0 = \sum_{w \in \Gamma(v_{i-1,k})}f(w) = f(v_{i,k}) + \sum_{k' \in \{1,...,p\}}f(v_{i-2,k'}) = f(v_{i,k}) + p \cdot f(v_{i-2,1}) = f(v_{i,k}),$$
а для произвольного нечетного $i \in \mathbb{Z}$ и произвольного $k \in \{1,...,p\}$ имеем
$$0 = \sum_{w \in \Gamma(v_{i+1,k})}f(w) = f(v_{i,k}) + \sum_{k' \in \{1,...,p\}}f(v_{i+2,k'}) = f(v_{i,k}) + p \cdot f(v_{i+2,1}) = f(v_{i,k}).$$
Таким образом, $f$ --- функция, тождественно равная 0 на $V(\Gamma)$. Как следствие, $0$ не является собственным значением
оператора смежности $A_{\Gamma,F}$.

\subsection{}
\label{s8.2}
Приведем {\it пример бесконечного кубического
связного графа $\Gamma$, для которого корни уравнения $x^3 - x^2 - 6x + 2 = 0$ в поле  $\mathbb C$ $($для
определенности с обычным абсолютным значением$)$ не являются
собственными значениями оператора смежности $A_{\Gamma, \mathbb C}$}.

Положим
$$V(\Gamma) = \{u_{i,1},...,u_{i,5},v_i : i \in \mathbb{Z}\},$$
$$E(\Gamma) = \{\{u_{i,1},u_{i,2}\}, \{u_{i,1},u_{i,3}\}, \{u_{i,2},u_{i,4}\}, \{u_{i,2},u_{i,5}\},
\{u_{i,3},u_{i,4}\}, \{u_{i,3},u_{i,5}\},$$
$$\{u_{i,4},u_{i,5}\}:
i \in \mathbb{Z}\}\cup \{\{u_{i,1},v_{i}\} : i \in \mathbb{Z}\} \cup \{\{v_i,v_{i+1}\} : i \in \mathbb{Z}\}.$$
Тогда $\Gamma$ --- бесконечный кубический
связный граф. Пусть $\lambda$ --- произвольный комплексный корень уравнения $x^3 - x^2 - 6x + 2 = 0$, и пусть
$f \in \mathbb{C}^{V(\Gamma)}$ такова, что
$$\lambda f(w) = \sum_{w' \in \Gamma(w)} f(w')$$
для каждой вершины $w$ графа $\Gamma$.

Предположим, что $f(v_j) \not = 0$ для некоторого $j \in \mathbb{Z}$. Тогда без потери общности можно предполагать,
что $f(v_j) = 1$ и, кроме того,
что $f(u_{j,2}) = f(u_{j,3})$ (в силу $\Gamma (u_{j,2}) = \Gamma (u_{j,3})$ и $\lambda \not = 0$) и
$f(u_{j,4}) = f(u_{j,5})$ (в силу наличия у $\Gamma$ автоморфизма, меняющего местами $u_{j,4}$ с $u_{j,5}$ и
стабилизирующего остальные вершины графа $\Gamma$). Но тогда
$$\lambda f(u_{j,1}) = 1 + 2 f(u_{j,2}),$$ $$\lambda f(u_{j,2}) = f(u_{j,1}) + 2 f(u_{j,4}),$$ $$\lambda f(u_{j,4}) = 2 f(u_{j,2}) + f(u_{j,4}),$$
что влечет $(\lambda^3 - \lambda^2 - 6 \lambda + 2) f(u_{j,4}) = 2$. Противоречие с выбором $\lambda$
доказывает, что $f(v_i) = 0$ для всех $i \in \mathbb{Z}$.

Для произвольного $i \in \mathbb{Z}$ из $f(v_{i-1}) = f(v_i) = f(v_{i+1}) = 0$ и $\lambda f(v_i) = f(v_{i-1}) + f(v_{i+1}) + f(u_{i,1})$
следует $f(u_{i,1}) = 0$, что с учетом
$\lambda f(u_{i,1}) =  f(v_i) + f(u_{i,2}) + f(u_{i,3})$ и равенства
$f(u_{i,2}) = f(u_{i,3})$ (справедливого в силу $\Gamma (u_{i,2}) = \Gamma (u_{i,3})$ и $\lambda \not = 0$)
влечет $f(u_{i,2}) = f(u_{i,3}) = 0$. Наконец, из $\lambda f(u_{i,4}) = f(u_{i,2}) + f(u_{i,3}) + f(u_{i,5}) = f(u_{i,5})$
и $\lambda f(u_{i,5}) = f(u_{i,2}) + f(u_{i,3}) + f(u_{i,4}) = f(u_{i,4})$ с учетом $\lambda^2 \not = 1$
следует, что $f(u_{i,4}) = f(u_{i,5}) = 0$. Таким образом, $f = 0$ и $\lambda$ не является
собственным значением оператора смежности $A_{\Gamma,\mathbb C}$
графа $\Gamma$.

Отметим, что группа $Aut(\Gamma)$ автоморфизмов графа $\Gamma$ имеет 4 орбиты на $V(\Gamma)$.

\begin{remark}
\label{r8.1}
Использованная при построении этого примера схема имеет следующий вид.
Пусть $\Delta$ --- конечный связный граф, $u \in V(\Delta)$ и $\lambda \in \mathbb C$.
Предположим, что

$(1)$ у графа $\Delta$ отсутствуют $(\mathbb C, \lambda)$-пропагаторы относительно вершины $u$.

Ясно, что при выполнении условия $(1)$ число $\lambda$
является собственным значением матрицы смежности графа $\Delta$ (рассматриваемой над $\mathbb C$).
В частности, $\lambda$ --- вполне вещественное целое алгебраическое число. Предположим дополнительно, что

$(2)$ $\lambda$ --- простое (т. е. не кратное) собственное значение матрицы смежности графа $\Delta$.

Согласно теореме~\ref{t4.1} при выполнении условий (1) и (2) соответствующий $\lambda$ собственный вектор матрицы
смежности графа $\Delta$,
рассматриваемый как функция из $\mathbb C^{V(\Delta)}$, принимает ненулевое значение в вершине $u$.

(Для приводимого выше примера $V(\Delta) = \{u = u_{1},...,u_{5}\},$
$E(\Delta) = \{\{u_{1},u_{2}\},$ $\{u_{1},u_{3}\}, \{u_{2},u_{4}\}, \{u_{2},u_{5}\},
\{u_{3},u_{4}\}, \{u_{3},u_{5}\},\{u_{4},u_{5}\}$, $\lambda$ --- любой комплексный
корень уравнения $x^3 - x^2 - 6x + 2 = 0$.)

Для каждого $i \in \mathbb{Z}$ пусть $\Delta _i$ --- граф, на который имеется изоморфизм $\varphi_i$
графа $\Delta$. Определим граф $\Gamma$, полагая
$V(\Gamma) = (\cup_{i \in \mathbb{Z}} V(\Delta_i)) \cup \{v_i : i \in \mathbb{Z}\}$
(предполагается, что множества $V(\Delta_i)$, $i \in \mathbb{Z}$, попарно дизъюнктны и
дизъюнктны с множеством $\{v_i : i \in \mathbb{Z}\}$),
$E(\Gamma) = (\cup_{i \in \mathbb{Z}} E(\Delta_i)) \cup \{\{\varphi_i(u),v_i\} : i \in \mathbb{Z}\} \cup \{\{v_i,v_{i+1}\} : i \in \mathbb{Z}\}$.
Тогда $\Gamma$ --- бесконечный локально конечный связный граф такой, что $\lambda$
не являются собственным значением оператора смежности $A_{\Gamma,\mathbb C}$.
Действительно, если $f \in \mathbb C ^{V(\Gamma)}$ --- собственная функция оператора смежности
$A_{\Gamma,\mathbb{C}}$,
соответствующая собственному значению $\lambda$, то с учетом $(1)$ имеем $f(v_i) = 0$ для всех
$i \in \mathbb{Z}$ (из $f(v_i) \not = 0$ для некоторого
$i \in \mathbb{Z}$ следовало бы наличие $(\mathbb C, \lambda)$-пропагатора графа $\Delta_i$ относительно
$\varphi_i(u)$, каковым являлось бы ограничение на $V(\Delta_i)$ функции $-f/f(v_i)$),
что с учетом $(2)$ влечет равенство $f$ нулю всюду на $V(\Gamma)$ (наличие
ненулевого ограничения $f$ на $V(\Delta_i)$ для $i \in \mathbb{Z}$
ввиду $f(v_i) = 0$,
$f(\varphi_i(u)) = \lambda f(v_i) - f(v_{i-1}) - f(v_{i+1}) = 0$ и замечания, сделанного сразу после
формулировки условия (2), противоречило бы условию (2)); противоречие.

В заключение еше раз обратим внимание на то, что используемое в этой схеме построения число $\lambda$
является вещественным.

\end{remark}

\subsection{}
\label{s8.3}
В связи с замечанием~\ref{r8.1} и замечанием~\ref{r6.10} не лишен интереса вопрос о существовании
бесконечного локально конечного связного графа $\Gamma$ такого, что некоторое
$\lambda \in \mathbb{C} \setminus \mathbb{R}$ не является
собственным значением оператора смежности $A_{\Gamma, \mathbb C}$ (с обычным абсолютным
значением на $\mathbb C$ для определенности).
Положительный ответ на него дает приводимый здесь {\it пример $($бесконечного$)$ кубического
связного графа $\Gamma$ такого, что любой комплексный невещественный корень $\lambda$
уравнения $x^3 + x^2 - 1 = 0$} (которое имеет два комплексно сопряженных невещественных корня и один вещественный
корень) {\it не является
собственным значением оператора смежности $A_{\Gamma, \mathbb C}$}. (В силу неприводимости многочлена
$x^3 + x^2 - 1 = 0$ над $\mathbb{Q}$ отсюда будет следовать (см. предложение~\ref{p3.1}), что и вещественный корень
уравнения $x^3 + x^2 - 1 = 0$ не является
собственным значением $A_{\Gamma, \mathbb C}$.)

Пусть $\Delta$ --- граф с множеством вершин $$V(\Delta) = \{u_j : 1 \leq j \leq 6\}$$
и множеством ребер $$E(\Delta) = \{\{u_1,u_2\},\{u_1,u_3\},\{u_2,u_4\},\{u_2,u_5\},
\{u_3,u_4\},$$ $$\{u_3,u_6\},\{u_4,u_5\},\{u_5,u_6\}\}.$$
Граф $\Delta$ обладает тем (не совсем обычным для конечных графов) свойством, что для невещественного $\lambda$
его $(\mathbb C, \lambda)$-пропагатор относительно одной вершины, а именно вершины $u_1$ или $u_6$,
принимает нулевое значение в другой вершине, а именно соответственно в вершине $u_6$ или $u_1$.

Пусть $\check \Delta$ --- граф с множеством вершин $V(\check \Delta) = V(\Delta) \cup \{v,v'\}$
и множеством ребер $E(\check \Delta) = E(\Delta) \cup \{\{v,u_1\},\{v',u_6\}\}$.
Обозначим через $U$ подпространство векторного пространства $\mathbb C^{V(\check \Delta)}$
комплекснозначных функций на $V(\check \Delta)$, состоящее из всех таких $f \in \mathbb C^{V(\check \Delta)}$,
что $\lambda f(u_j) = \sum _{u \in \check \Delta (u_j)} f(u)$ для всех $1 \leq j \leq 6$.
Подпространство $U$ имеет размерность $2$. В качестве базиса $U$ могут быть взяты функции $f_1$ и $f_2$,
определяемые следующим образом.

$f_1(v) = 1-2\lambda, f_1(u_1) = -2-\lambda, f_1(u_2) = \lambda, f_1(u_3)=-(1+\lambda)/\lambda,$
$f_1(u_4) = 1, f_1(u_5) = (1+\lambda)/\lambda, f_1(u_6) = 0, f_1(v') = 0.$

$f_2(w) = f_1(g(w))$ для всех $w \in V(\check \Delta)$, где $g$ --- автоморфизм графа $\check \Delta$,
стабилизирующий вершины $u_3$, $u_4$ и меняющий местами вершины $v$ и $v'$, $u_1$ и $u_6$, $u_2$ и $u_5$.

Для дальнейшего важно заметить, что для произвольного $c \in \mathbb C$ функции $(c/(1-2\lambda))f_1$
и $(c/(1-2\lambda))f_2$ из $U$
таковы, что
$(c/(1-2\lambda))f_1(v) = c$, $(c/(1-2\lambda))f_1(u_6) = 0 = (c/(1-2\lambda))f_1(v')$,
$(c/(1-2\lambda))f_2(v') = c$, $(c/(1-2\lambda))f_2(u_1) = 0 = (c/(1-2\lambda))f_2(v)$
и
$\lambda (c/(1-2\lambda))f_k(u_j) = \sum _{u \in \check \Delta (u_j)} (c/(1-2\lambda))f_k(u)$ для всех $1 \leq k \leq 2$, $1 \leq j \leq 6$.

Пусть теперь $\Gamma$ --- граф с множеством вершин $$V(\Gamma) = \{u_{i,j} : i \in \mathbb{Z}, 1 \leq j \leq 6\}$$
и множеством ребер $$E(\Gamma) = \{\{u_{i,1},u_{i,2}\},\{u_{i,1},u_{i,3}\},\{u_{i,2},u_{i,4}\},\{u_{i,2},u_{i,5}\},$$
$$\{u_{i,3},u_{i,4}\},\{u_{i,3},u_{i,6}\},\{u_{i,4},u_{i,5}\},\{u_{i,5},u_{i,6}\},\{u_{i,6},u_{i+1,1}\}: i \in \mathbb{Z}\}.$$
Граф $\Gamma$ --- бесконечный кубический связный граф. Покажем, что число $\lambda$ не является
собственным значением оператора смежности $A_{\Gamma, \mathbb C}$.

Предположим противное, и пусть $\tilde f \in \mathbb C^{V(\Gamma)}$ --- собственная функция оператора смежности $A_{\Gamma, \mathbb C}$,
соответствующая собственному значению $\lambda$.
В силу невещественности $\lambda$ подграф графа $\Gamma$, порожденный носителем функции $\tilde f$
не имеет конечных связных компонент (поскольку ограничение $\tilde f$ на такую компоненту являлось бы собственным вектором матрицы смежности
конечного графа, соответствующим $\lambda$). В частности, найдется такое $m \in \mathbb{Z}$, что
$\tilde f(u_{i,1}) \not = 0 \not = \tilde f(u_{i,6})$
для всех $i \in \mathbb{Z}_{\geq m}$ или для всех $i \in \mathbb{Z}_{\leq m}$. В силу наличия у $\Gamma$ автоморфизма $u_{i,1} \mapsto u_{-i,6},
u_{i,2} \mapsto u_{-i,5}, u_{i,3} \mapsto u_{-i,3}, u_{i,4} \mapsto u_{-i,4}, u_{i,5} \mapsto u_{-i,2},
u_{i,6} \mapsto u_{-i,1}$, $i \in \mathbb{Z},$
и автоморфизма $u_{i,j} \mapsto u_{i+1,j}$, $i \in \mathbb{Z}, 1 \leq j \leq 6,$
будем, не теряя общности, предполагать, что $m = 0$ и
$\tilde f(u_{i,1}) \not = 0 \not = \tilde f(u_{i,6})$ для всех $i \in \mathbb{Z}_{\geq 0}$. (Для дальнейшего важно, собственно, лишь то,
что $\tilde f$ может быть выбрана с отличными от $0$ значениями $\tilde f(u_{0,6})$ и $\tilde f(u_{1,1})$.)

Положим $X_0 = \{u_{-1,6},u_{0,1},...,u_{0,6},u_{1,1}\}$, $X_1 = \{u_{0,6},u_{1,1},...,u_{1,6},u_{2,1}\}$.
Каждый из подграфов $\langle X_0 \rangle_{\Gamma}$ и $\langle X_1 \rangle_{\Gamma}$ графа $\Gamma$
очевидным образом изоморфен графу $\check \Delta$, что с учетом наличия функций $f_1$ и $f_2$ (см. выше)
влечет наличие функций $f^- \in \mathbb C^{X_0}$ и $f^+ \in \mathbb C^{X_1}$ таких,
что $f^-(u_{-1,6}) = \tilde f(u_{-1,6})$, $f^-(u_{0,6}) = 0 = f^-(u_{1,1})$,
$\lambda f^-(u_{0,j}) = \sum _{u \in \Gamma (u_{0,j})} f^-(u)$ для всех $1 \leq j \leq 6$,
$f^+(u_{2,1}) = \tilde f(u_{2,1})$, $f^+(u_{0,6}) = 0 = f^+(u_{1,1})$,
$\lambda f^+(u_{1,j}) = \sum _{u \in \Gamma (u_{1,j})} f^+(u)$ для всех $1 \leq j \leq 6$.
Продолжим функции $f^-$, $f^+$ до функций $f_{ext}^-$, $f_{ext}^+$ из $\mathbb C^{X_0 \cup X_1}$,
полагая $f_{ext}^-(w) = 0$ для всех $w \in X_1$ и $f_{ext}^+(w) = 0$ для всех $w \in X_0$.
Пусть $\hat f \in \mathbb C^{X_0 \cup X_1}$, $\hat f (w) = \tilde f(w) - f_{ext}^-(w) - f_{ext}^+(w)$
для всех $w \in X_0 \cup X_1$. Тогда $\hat f(u_{-1,6}) = 0 = \hat f(u_{2,1})$,
$\hat f(u_{0,6}) = \tilde f(u_{0,6}) \not = 0 \not = \tilde f(u_{1,1}) = \hat f(u_{1,1})$,
$\hat f (u_{i,j}) = \sum _{u \in \Gamma (u_{i,j})} \hat f(u)$ для всех $1\leq i \leq 2, 1 \leq j \leq 6$.
Но тогда ограничение функции $\hat f$ на множество $X := (X_0 \cup X_1) \setminus \{u_{-1,6},u_{2,1}\}$
есть собственная функция матрицы смежности конечного графа $\langle X \rangle_{\Gamma}$, соответствующая
собственному значению $\lambda$. Последнее, однако, невозможно в силу невещественности $\lambda$,
что завершает доказательство.

Отметим, что группа $Aut(\Gamma)$ автоморфизмов графа $\Gamma$ имеет 4 орбиты на $V(\Gamma)$.

\subsection{}
\label{s8.4}
Для графа из раздела~\ref{s8.2} и графов из замечания~\ref{r8.1} числа вершинной и реберной связности равны 1.
Несколько сложнее построить пример бесконечного локально конечного регулярного
графа, у которого числа вершинной и реберной связности больше 1,
но для которого не все комплексные числа являются собственными значениями
оператора смежности над $\mathbb C$ (с обычным абсолютным
значением для определенности). Ниже строится {\it пример бесконечного кубического
графа $\Gamma$, у которого числа вершинной и реберной связности равны $2$,
но для которого число $0$ не является собственным значением оператора смежности $A_{\Gamma, \mathbb C}$}.

Для каждого $i \in  \mathbb{Z}$ следующим образом определим граф $\Delta_i$:
$$V(\Delta_i) = \{u_{i,1},...,u_{i,16}\},$$
$$E(\Delta_i) = \{\{u_{i,1},u_{i,3}\}, \{u_{i,2},u_{i,4}\}, \{u_{i,3},u_{i,5}\}, \{u_{i,3},u_{i,6}\},
\{u_{i,4},u_{i,7}\}, \{u_{i,4},u_{i,8}\},$$
$$\{u_{i,5},u_{i,8}\}, \{u_{i,5},u_{i,9}\}, \{u_{i,6},u_{i,7}\},
\{u_{i,6},u_{i,9}\}, \{u_{i,7},u_{i,10}\}, \{u_{i,8},u_{i,11}\},$$
$$\{u_{i,9},u_{i,12}\}, \{u_{i,10},u_{i,13}\},
\{u_{i,10},u_{i,14}\}, \{u_{i,11},u_{i,12}\}, \{u_{i,11},u_{i,15}\},$$
$$\{u_{i,12},u_{i,13}\}, \{u_{i,13},u_{i,16}\}, \{u_{i,14},u_{i,15}\}, \{u_{i,14},u_{i,16}\}, \{u_{i,15},u_{i,16}\}\}.$$
(Предполагается, что $V(\Delta_{i'}) \cap V(\Delta_{i''}) = \emptyset$ для любых
различных $i', i'' \in \mathbb{Z}$.)

Пусть $i \in \mathbb Z$, и пусть $f \in \mathbb C^{V(\Delta_i)}$ такова, что
\begin{equation}
\label{eq8.1}
\sum_{u' \in \Delta_i(u)} f(u') = 0
\end{equation}
для всех $u \in V(\Delta_i)\setminus \{u_{i,1},u_{i,2}\}$.
Тогда, полагая $a := f(u_{i,10})$ и $b := f(u_{i,14})$
и используя~\eqref{eq8.1}, легко получаем
$$f(u_{i,1}) = f(u_{i,2}) = f(u_{i,7}) = f(u_{i,8}) = f(u_{i,12}) = f(u_{i,15}) = 0,$$
$$f(u_{i,3}) = a - 2b, f(u_{i,4}) = -a + \frac{1}{2}b, f(u_{i,5}) = \frac{1}{2}b, u_{i,6} = - \frac{1}{2}b,$$
$$f(u_{i,9}) = -a + 2b, f(u_{i,11}) = a - b, f(u_{i,13}) = -b, f(u_{i,16}) = -a.$$

Определим следующим образом граф $\Gamma$:
$$V(\Gamma) = \cup_{i \in \mathbb{Z}} V(\Delta_i) =  \{u_{i,1},...,u_{i,16} : i \in \mathbb{Z}\},$$
$$E(\Gamma) = (\cup_{i \in \mathbb{Z}} E(\Delta_i)) \cup \{\{u_{i,1},u_{i+1,1}\},
\{u_{i,2},u_{i+1,2}\}: i \in  \mathbb{Z}\}.$$
Тогда $\Gamma$ --- бесконечный кубический
граф, числа вершинной и реберной связности которого равны 2.
Покажем, что 0
не является собственным значением оператора смежности графа $\Gamma$.
Пусть $\tilde f \in \mathbb{C}^{V(\Gamma)}$ такова, что
\begin{equation}
\label{eq8.2}
\sum_{u' \in \Gamma(u)} \tilde f(u') = 0
\end{equation}
для всех $u \in V(\Gamma)$. Для произвольного $i \in \mathbb{Z}$,
применяя к ограничению функции $\tilde f$ на $V(\Delta_i)$ утверждение,
полученное в предыдущем абзаце, заключаем, что
\begin{equation}
\label{eq8.3}
\tilde f(u_{i,1}) = 0 = \tilde f(u_{i,2}).
\end{equation}
Далее, для произвольного $i \in \mathbb{Z}$ имеем
$$\Gamma(u_{i,1}) =  \{u_{i-1,1}, u_{i,3}, u_{i+1,1}\}, \Gamma(u_{i,2}) =  \{u_{i-1,2}, u_{i,4}, u_{i+1,2}\},$$
что с учетом~\eqref{eq8.2} и~\eqref{eq8.3} влечет $\tilde f(u_{i,3}) = 0 = \tilde f(u_{i,4})$ для всех $i \in \mathbb{Z}$.
Но тогда, вновь применяя к ограничению функции $\tilde f$ на $V(\Delta_i), i \in \mathbb{Z},$ утверждение,
полученное в предыдущем абзаце, заключаем, что $\tilde f(u) = 0$ для всех $u \in V(\Delta_i), i \in \mathbb{Z}.$
Таким образом, каждая функция $\tilde f \in \mathbb{C}^{V(\Gamma)}$ с условием~\eqref{eq8.2} тождественно равна 0 на $V(\Gamma)$.
Следовательно, 0 не является собственным значением оператора смежности графа $\Gamma$.

Отметим, что группа $Aut(\Gamma)$ автоморфизмов графа $\Gamma$ имеет конечное число орбит на $V(\Gamma)$.

\subsection{}
\label{s8.5}
Легко привести {\it примеры бесконечных локально конечных связных вер\-шинно-симметрических графов $\Gamma$ таких, что для подходящего
поля $F$ с произвольным абсолютным значением и подходящего элемента $\lambda \in F$ имеем $S_{F,\lambda}(\Gamma) = V(\Gamma)$}.

\subsubsection{}
\label{s8.5.1}
Пусть $F$ --- поле (с произвольным абсолютным значением) и $\lambda \in F$ является собственным значением матрицы смежности над $F$ конечного связного вершинно-симметрического графа
$\Lambda$. Укажем пример бесконечного локально конечного связного вершинно-симмет\-рического графа $\Gamma$, для которого $S_{F,\lambda}(\Gamma) = V(\Gamma)$.
Пусть $\Delta$ --- произвольный  бесконечный локально конечный связный вершинно-симметрический граф. Положим
$$V(\Gamma) = \{(v,w,j): v \in V(\Delta), w \in V(\Lambda), j \in \{1,2,3,4\}\},$$
$$E(\Gamma) = \{\{(v_1,w_1,j_1),(v_2,w_2,j_2)\} : \{v_1,v_2\} \in E(\Delta), w_1,w_2 \in V(\Lambda),$$
$$j_1,j_2 \in \{1,2,3,4\}\} \cup  \{\{(v,w_1,j_1),(v,w_2,j_2)\} : v \in V(\Delta), w_1,w_2 \in V(\Lambda),$$
$$j_1,j_2 \in \{1,2,3,4\}, |j_1-j_2| \in \{1, 3\}\} \cup  \{\{(v,w_1,j),(v,w_2,j)\} : v \in V(\Delta),$$ $$\{w_1,w_2\} \in E(\Lambda),  j \in \{1,2,3,4\}\}.$$
Ясно, что $\Gamma$ ---  бесконечный локально конечный связный вершинно-симметри\-ческий граф.
Далее, поскольку $\lambda \in F$ --- собственное значение матрицы смежности над $F$ графа
$\Lambda$, имеется собственная функция $f \in F^{V(\Lambda)}$ матрицы смежности над $F$ графа
$\Lambda$, соответствующая собственному значению $\lambda$. Пусть $v_0 \in V(\Delta)$.
Определим функцию $\tilde f \in F^{V(\Gamma)}$, полагая $\tilde f((v_0,w,1)) = f(w)$ для всех
$w \in V(\Lambda)$ и
$\tilde f((v_0,w,3)) = -f(w)$ для всех $w \in V(\Lambda)$, а во всех других
вершинах графа $\Gamma$ полагая $\tilde f$
равной $0$. Легко убедиться, что $\tilde f$ --- собственная функция оператора $A_{\Gamma,F}$, соответствующая
собственному значению $\lambda$, причем $\tilde f$ имеет конечный носитель. Согласно теореме~\ref{t4.1} это влечет
$S_{F,\lambda}(\Gamma) \not = \emptyset$, что в силу вершинной симметричности графа $\Gamma$ дает
требуемое равенство $S_{F,\lambda}(\Gamma) = V(\Gamma)$.

\subsubsection{}
\label{s8.5.2}
Если $\Gamma$ --- произвольный бесконечный локально конечный связный
вершинно-симметрический граф с тем свойством, что для $u \in V(\Gamma)$ найдется $v \in V(\Gamma)$, для которой
$\{u,v\} \in E(\Gamma)$ и $\Gamma(u) \setminus \{v\} = \Gamma(v) \setminus \{u\}$, то, как легко видеть,
для произвольного поля $F$ (с произвольным абсолютным значением) имеем $S_{F,-1}(\Gamma) = V(\Gamma)$.
Сходным образом, если $\Gamma$ --- произвольный бесконечный локально конечный связный
вершинно-симметрический граф с тем свойством, что для $u \in V(\Gamma)$ найдется $v \in V(\Gamma) \setminus \{u\}$, для которой
$\Gamma(u) = \Gamma(v)$, то для произвольного поля $F$ (с произвольным абсолютным значением)
имеем $S_{F,0}(\Gamma) = V(\Gamma)$.

\end{fulltext}

\end{document}